\documentclass[a4paper,oneside]{article}
\usepackage{amssymb}
\usepackage{amsfonts}
\usepackage{amsmath}

\input{tcilatex}

\begin{document}

\title{Schubert calculus and cohomology of Lie groups. Part II. Compact Lie
groups}
\author{Haibao Duan\thanks{%
The author's research is supported by 973 Program 2011CB302400 and NSFC
11131008.} \\
Hua Loo-Keng Key Laboratory of Mathematics, \\
Institute of Mathematics, Chinese Academy of Sciences,\\
dhb@math.ac.cn}
\maketitle

\begin{abstract}
Let $G$ be a compact Lie group with a maximal torus $T$. Based on a Schubert
presentation on the integral cohomology $H^{\ast }(G/T)$ of the flag
manifold $G/T$ \cite{DZ1} we have presented in \cite{DZ2} an explicit and
unified construction of the integral cohomology $H^{\ast }(G)$ for all $1$%
--connected Lie groups $G$. In this paper we extend this construction to the
compact Lie groups.

\begin{description}
\item[2000 Mathematical Subject Classification: ] 55T10, 14M15

\item[Key words and phrases:] Lie groups; Cohomology, Schubert calculus
\end{description}
\end{abstract}

\section{Introduction}

The Lie group $G$ under consideration is compact and connected; the
coefficient for the cohomologies is either the ring $\mathbb{Z}$ of
integers, or one of the finite fields $\mathbb{F}_{p}$, unless otherwise
stated.

For a maximal torus $T$ on $G$ the homogeneous space $G/T$ is canonically a
projective variety, called \textsl{the complete flag manifold} of the Lie
group $G$. Based on a Schubert presentation on the integral cohomology ring
of $G/T$ \cite[Theorem 1.2]{DZ1} we have presented in \cite{DZ2} a unified
construction for the cohomologies of all $1$--connected Lie groups. In this
paper we extend this construction to the compact Lie groups. As applications
of our general approach the integral cohomologies of the adjoint Lie groups $%
PSU(n),PSp(n),PE_{6},PE_{7}$ are determined, see Theorem 4.7, Theorem 4.12,
as well as the historical remarks in Section 4.5.

The problem of computing the cohomologies of Lie groups was raised by Cartan
in 1929. It is a focus of algebraic topology for the fundamental roles of
Lie groups playing in geometry and topology, see \cite[Chapter VI]{D}, \cite%
{K,MT}. On the other hand, the classical Schubert calculus amounts to the
determination of the cohomology rings of the flag manifolds $G/T$ \cite[p.331%
]{W}. The present work, together with the companion ones \cite{DZ2,DZ3},
completes our project to determine the integral cohomologies of all compact
Lie groups $G$, as well as the structure of the $\func{mod}p$ cohomology $%
H^{\ast }(G;\mathbb{F}_{p})$ as a module over the Steenrod algebra $\mathcal{%
A}_{p}$, in the context of Schubert calculus.

We begin with a general procedure that reduces the construction and
computation with the cohomology of a compact Lie group to that of the $1$%
--connected ones. It is well known that all the $1$\textsl{--connected} 
\textsl{simple Lie groups} $G$, together with their centers $\mathcal{Z}(G)$%
, are classified by the types $\Phi _{G}$\ of their root systems showing in
the following table, where $e\in G$ is the group unit:

\begin{center}
\textbf{Table 1.1.} The types and centers of the $1$--connected simple Lie
groups

{\footnotesize 
\begin{tabular}{l|l|l|l|l|l|l|l|l|l}
\hline\hline
${\small G}$ & ${\small SU(n)}$ & ${\small Sp(n)}$ & ${\small Spin(2n+1)}$ & 
${\small Spin(2n)}$ & ${\small G}_{2}$ & ${\small F}_{4}$ & ${\small E}_{6}$
& ${\small E}_{7}$ & ${\small E}_{8}$ \\ \hline
${\small \Phi }_{{\small G}}$ & $A_{n-1}$ & $B_{n}$ & $C_{n}$ & $D_{n}$ & $%
{\small G}_{2}$ & ${\small F}_{4}$ & ${\small E}_{6}$ & ${\small E}_{7}$ & $%
{\small E}_{8}$ \\ \hline
$\mathcal{Z}{\small (G)}$ & $\mathbb{Z}_{n}$ & $\mathbb{Z}_{2}$ & $\mathbb{Z}%
_{2}$ & 
\begin{tabular}{l}
$\mathbb{Z}_{4}${\small ,\ }${\small n=2k+1}$ \\ 
$\mathbb{Z}_{2}{\small \oplus }\mathbb{Z}_{2}${\small , }${\small n=2k}$%
\end{tabular}
& ${\small \{e\}}$ & ${\small \{e\}}$ & $\mathbb{Z}_{3}$ & $\mathbb{Z}_{2}$
& ${\small \{e\}}$ \\ \hline
\end{tabular}
}.
\end{center}

\noindent In general, let $G$ be a compact Lie group with $\mathcal{Z}%
_{0}(G) $ the identity component of the center of $G$, and with $G^{\prime }$
the commutator subgroup of $G$. Then the intersection $F=G^{\prime }\cap 
\mathcal{Z}_{0}(G)$ is always a finite abelian group. The Cartan's
classification on compact Lie groups states that (\cite[Theorem 5.22]{S})

\bigskip

\noindent \textbf{Theorem 1.1.} \textsl{The isomorphism type of a compact
Lie group }$G$\textsl{\ is}

\begin{enumerate}
\item[(1.1)] $G=\left[ G^{\prime }\times \mathcal{Z}_{0}(G)\right] /F$,
\end{enumerate}

\noindent \textsl{where }$F$\textsl{\ is embedded in} \textsl{the numerator
group} $G^{\prime }\times \mathcal{Z}_{0}(G)$ \textsl{as} $\{(g,g^{-1})\mid
g\in F\}$.

\textsl{Moreover, the commutator subgroup }$G^{\prime }$\textsl{\ admits\ a
canonical presentation as}

\begin{enumerate}
\item[(1.2)] $G^{\prime }=\left[ G_{1}\times \cdots \times G_{k}\right] /K$%
\textsl{,}
\end{enumerate}

\noindent \textsl{where each }$G_{t}$\textsl{, }$1\leq t$\textsl{\ }$\leq k$%
\textsl{, is one of the }$1$\textsl{--connected simple Lie groups listed in
Table 1.1, and where }$K$\textsl{\ is a subgroup of} \textsl{the finite group%
} $\mathcal{Z}(G_{1})\times \cdots \times \mathcal{Z}(G_{k})$.$\square $

\bigskip

In views of (1.1) and (1.2) a Lie group $G$ is called \textsl{semi--simple}
if $\mathcal{Z}_{0}(G)=\{e\}$; \textsl{simple} if $\mathcal{Z}_{0}(G)=\{e\}$
and $k=1$. Since the commutator subgroup $G^{\prime }$ is always
semi--simple we shall call it \textsl{the semi--simple part} of $G$. Based
on Theorem 1.1 we give a diffeomorphism classification of compact Lie groups
in the following result. Let $T^{r}$ be the $r$\textsl{--}dimensional torus
group.

\bigskip

\noindent \textbf{Theorem 1.2.} \textsl{The diffeomorphism type of a compact
Lie group }$G$\textsl{\ with semi--simple part }$G^{\prime }$ \textsl{is}

\begin{enumerate}
\item[(1.3)] $G\cong $\textsl{\ }$G^{\prime }\times T^{r}$\textsl{,} $r=\dim 
\mathcal{Z}_{0}(G)$.
\end{enumerate}

\noindent \textbf{Proof. }Since $G^{\prime }$ is normal in $G$ the quotient
space $G/G^{\prime }=\mathcal{Z}_{0}(G)/F$ has the structure of an abelian
group isomorphic to the $r$\textsl{--}dimensional torus group\textsl{\ }$%
T^{r}$, $r=\dim \mathcal{Z}_{0}(G)$, and the quotient map $h:G\rightarrow
G/G^{\prime }=T^{r}$ is both a group homomorphism and a submersion with
fiber $G^{\prime }$.

Take a maximal torus $T^{\prime }$ on $G^{\prime }$. By (1.1) a maximal
torus on $G$ is $T=\left[ T^{\prime }\times \mathcal{Z}_{0}(G)\right] /F$.
Since the restriction $h\mid T:T\rightarrow T^{r}$ of $h$ on $T$ is a fiber
bundle in torus groups there is a monomorphism $\sigma :T^{r}\rightarrow T$
so that the composition $\left( h\mid T\right) \circ \sigma $ is the
identity on $T^{r}$. Since $\sigma $ can be viewed as a section of $h$, one
obtains the diffeomorphism (1.3) from the fact that the fiber of $h$ is a
group.$\square $

\bigskip

The unitary group $U(n)$ of order $n$ may serve as the first example
illustrating the subtle difference between the isomorphism and
diffeomorphism types of a Lie group. As a group it is isomorphic to $\left[
SU(n)\times S^{1}\right] /\mathbb{Z}_{n}$, while as a smooth manifold it is
diffeomorphic to the product space $SU(n)\times S^{1}$.

Let $G$ be a semi--simple Lie group whose center $\mathcal{Z}(G)$ contains
the cyclic group $\mathbb{Z}_{q}$ of order $q$. Consider the cyclic covering 
$c:G\rightarrow G/\mathbb{Z}_{q}$ of Lie groups. Since the classifying space 
$B\mathbb{Z}_{q}$ of the group $\mathbb{Z}_{q}$ is the Elienberg--MaClane
space $K(\mathbb{Z}_{q},1)$ the classifying map $f_{c}:$ $G/\mathbb{Z}%
_{q}\rightarrow B\mathbb{Z}_{q}$ of $c$ defines a cohomology class $\iota
\in H^{1}(G/\mathbb{Z}_{q};\mathbb{F}_{q})$, called \textsl{the
characteristic class of the covering}.

On the other hand let $\mathbb{Z}_{q}$ act on the circle group $S^{1}$ as
the anti--clockwise rotation through the angle $2\pi /q$. Then the obvious
group homomorphism

\begin{enumerate}
\item[(1.4)] $C:\left[ G\times S^{1}\right] /\mathbb{Z}_{q}\rightarrow G/%
\mathbb{Z}_{q}$
\end{enumerate}

\noindent is an oriented cycle bundle on the quotient group $G/\mathbb{Z}%
_{q} $ with Euler class $\omega =\beta _{q}(\iota )\in H^{2}(G/\mathbb{Z}%
_{q})$, where $\beta _{q}$ is the Bockstein homomorphism. The Gysin sequence
of the circle bundle $C$ then yields the exact sequence \cite[P.143]{MS}

\begin{enumerate}
\item[(1.5)] $\cdots \rightarrow H^{r}(G/\mathbb{Z}_{q})\overset{C^{\ast }}{%
\rightarrow }H^{r}(G\times S^{1})\overset{\theta }{\rightarrow }H^{r-1}(G/%
\mathbb{Z}_{q})\overset{\omega }{\rightarrow }H^{r+1}(G/\mathbb{Z}_{q})%
\overset{C^{\ast }}{\rightarrow }\cdots $
\end{enumerate}

\noindent where $\omega $\ denotes the homomorphism of taking product with
the class $\omega $, and where the space $\left[ G\times S^{1}\right] /%
\mathbb{Z}_{q}$ has been replaced by its diffeomorphism type $G\times S^{1}$
by Theorem 1.2. Let $p:G\times S^{1}\rightarrow S^{1}$\ be the projection
onto the second factor, and $\varepsilon \in H^{1}(S^{1})$ the canonical
orientation on $S^{1}$. In what follows we set

\begin{enumerate}
\item[(1.6)] $\xi _{1}:=p^{\ast }(\varepsilon )\in H^{1}(G\times S^{1})$.
\end{enumerate}

Let $J(\omega )$ and $\left\langle \omega \right\rangle $ be respectively
the subring and the ideal of $H^{\ast }(G/\mathbb{Z}_{q})$ generated by $%
\omega $. Write $H^{\ast }(G/\mathbb{Z}_{q})_{\left\langle \omega
\right\rangle }$ for the quotient ring $H^{\ast }(G/\mathbb{Z}%
_{q})/\left\langle \omega \right\rangle $ with quotient map $g$. Then, in
addition to the short exact sequence of rings

\begin{enumerate}
\item[(1.7)] $0\rightarrow \left\langle \omega \right\rangle \rightarrow
H^{\ast }(G/\mathbb{Z}_{q})\overset{g}{\rightarrow }H^{\ast }(G/\mathbb{Z}%
_{q})_{\left\langle \omega \right\rangle }\rightarrow 0$,
\end{enumerate}

\noindent the cohomology $H^{\ast }(G/\mathbb{Z}_{q})$ can be regarded as a
module over its subring $J(\omega )$.

\bigskip

\noindent \textbf{Theorem\textbf{\ }1.3.} \textsl{The induced map }$C^{\ast
} $ \textsl{fits into} \textsl{the exact sequence}

\begin{enumerate}
\item[(1.8)] $0\rightarrow H^{\ast }(G/\mathbb{Z}_{q})_{\left\langle \omega
\right\rangle }\overset{C^{\ast }}{\rightarrow }H^{\ast }(G\times S^{1})%
\overset{\theta }{\rightarrow }H^{\ast }(G/\mathbb{Z}_{q})\overset{\omega }{%
\rightarrow }\left\langle \omega \right\rangle \rightarrow 0$
\end{enumerate}

\noindent \textsl{in which the homomorphism }$\theta $\textsl{\ has the
following} \textsl{properties}

\begin{quote}
\textsl{i)} $\theta (\xi _{1})=q\in H^{0}(G/\mathbb{Z}_{q})$\textsl{;}

\textsl{ii)} $\theta (x\cup C^{\ast }(y))=\theta (x)\cup y$ \textsl{for} $%
x\in H^{\ast }(G\times S^{1})$ \textsl{and} $y\in H^{\ast }(G/\mathbb{Z}%
_{q}) $\textsl{.}
\end{quote}

\textsl{Moreover, if the map }$g$ \textsl{in (1.7) admits a split
homomorphism} $j:H^{\ast }(G/\mathbb{Z}_{q})_{\left\langle \omega
\right\rangle }$ $\rightarrow $ $H^{\ast }(G/\mathbb{Z}_{q})$\textsl{, then
the map}

\begin{quote}
$h:J(\omega )\otimes H^{\ast }(G/\mathbb{Z}_{q})_{\left\langle \omega
\right\rangle }\rightarrow H^{\ast }(G/\mathbb{Z}_{q})$
\end{quote}

\noindent \textsl{by} $h(\omega ^{r}\otimes x)=\omega ^{r}\cup j(x)$ \textsl{%
induces an isomorphism of }$J(\omega )$\textsl{--modules}

\begin{enumerate}
\item[(1.9)] $H^{\ast }(G/\mathbb{Z}_{q})\cong \frac{J(\omega )\otimes
H^{\ast }(G/\mathbb{Z}_{q})_{\left\langle \omega \right\rangle }}{%
\left\langle \omega \cdot \func{Im}\theta \right\rangle }$\textsl{.}
\end{enumerate}

\noindent \textbf{Proof.} The sequence (1.8) is easily seen to be a compact
form of (1.5), while the properties i) and ii) about the map $\theta $ are
also known, see \cite[Lemma 1]{Ma}.

The group $J(\omega )$ has the basis $\{1,\omega ,\cdots ,\omega ^{r}\}$ for
some $r\geq 1$. Granted with the map $j$ the sequence (1.7) is applicable to
show that every $x\in H^{\ast }(G/\mathbb{Z}_{q})$ admits an expansion of
the form

\begin{quote}
$x=j(a_{0})+\omega \cup j(a_{1})+\cdots +\omega ^{r}\cup j(a_{r})$, $%
a_{i}\in H^{\ast }(G/\mathbb{Z}_{q})_{\left\langle \omega \right\rangle }$, $%
0\leq i\leq r$.
\end{quote}

\noindent That is the map $h$ is surjective. To show the formula (1.9)
consider an element $y=a_{0}+\omega \otimes a_{1}+\cdots +\omega ^{r}\otimes
a_{r}\in \ker h$. One infers from $g\circ h(y)=0$ that $a_{0}=0$. It follows
now from $0=h(y)=\omega \cdot (1\cup j(a_{1})+\cdots +\omega ^{r-1}\cup
j(a_{r}))$ that

\begin{quote}
$1\cup j(a_{1})+\cdots +\omega ^{r-1}\cup j(a_{r})\in \func{Im}\theta $
\end{quote}

\noindent by (1.8). That is $\ker h=\left\langle \omega \cdot \func{Im}%
\theta \right\rangle $.$\square $

\bigskip

Granted with Theorems 1.2 and 1.3 we can clarify the main ideas of our
approach. Firstly, Theorem 1.2 reduces the cohomology of a compact Lie group 
$G$ to that of its semi--simple part $G^{\prime }$ by the Kunneth formula $%
H^{\ast }(G)=$ $H^{\ast }(G^{\prime })\otimes H^{\ast }(T^{r})$. Next, for a
semi--simple Lie group $G$ its universal covering $c:$ $G_{0}\rightarrow G$
can always be decomposed into a sequence of cyclic coverings

\begin{quote}
$c:G_{0}\overset{c_{1}}{\rightarrow }G_{1}\overset{c_{2}}{\rightarrow }%
\cdots \overset{c_{k}}{\rightarrow }G_{k}=G$
\end{quote}

\noindent in which the cohomology of the $1$--connected Lie group $G_{0}$ is
known \cite[Theorem 1.9]{DZ2}. The cohomology $H^{\ast }(G)$ in question can
be calculated, in principle, from the known one $H^{\ast }(G_{0})$ by a
repeatedly application of the exact sequence (1.8).

The remaining sections of the paper are so arranged. In Sections 2 and 3 we
develop the formulae and the constructions required to implement above
procedure. As applications the integral cohomology of the adjoint Lie groups 
$PSU(n),PSp(n),PE_{6},PE_{7}$ are calculated in Section 4.

In addition, since the present work is a natural extension of the earlier
one \cite{DZ2}, we shall be free to adopt the notation and results developed
in \cite{DZ2}.

\section{Spectral sequence of the fibration $G\rightarrow G/T$}

Fix a maximal torus $T$ on $G$, and consider the Leray--Serre spectral
sequence $\{E_{r}^{\ast ,\ast }(G),d_{r}\}$ of the corresponding fibration

\begin{enumerate}
\item[(2.1)] $T\hookrightarrow G\overset{\pi }{\rightarrow }G/T$.
\end{enumerate}

\noindent It is well known that its second page is the Koszul complex with

\begin{enumerate}
\item[(2.2)] $E_{2}^{\ast ,\ast }(G)=H^{\ast }(G/T)\otimes H^{\ast }(T)$;

\item[(2.3)] $d_{2}(a\otimes t)=(\tau (t)\cup a)\otimes 1$ for $a\in H^{\ast
}(G/T)$, $t\in H^{1}(T)$
\end{enumerate}

\noindent where $\tau :$ $H^{1}(T)\rightarrow H^{2}(G/T)$ is the\textsl{\
Borel transgression} in the fibration (2.1).

Our calculation and construction with the exact sequence (1.8) actually take
place on the third page of the spectral sequence. In order to access $%
E_{3}^{\ast ,\ast }(G)$ we deduce a formula for the transgression $\tau $ in
Theorem 2.4, and give a concise characterization for the factor ring $%
H^{\ast }(G/T)$ of $E_{2}^{\ast ,\ast }(G)$ in Theorem 2.6. These results
are essential in Section 3 for us to construct explicit generators of the
ring $H^{\ast }(G)$ by certain polynomials in the Schubert classes on $G/T$.

\subsection{A formula for the Borel transgression $\protect\tau $}

In the diagram with top row the cohomology exact sequence of the pair $(G,T)$

\begin{quote}
\begin{tabular}{lllll}
$0\rightarrow H^{1}(G)\overset{i^{\ast }}{\rightarrow }$ & $H^{1}(T)$ & $%
\overset{\delta }{\rightarrow }$ & $H^{2}(G,T)$ & $\overset{j^{\ast }}{%
\rightarrow }H^{2}(G)\rightarrow \cdots $ \\ 
&  & $\searrow \tau $ & $\quad \cong \uparrow \pi ^{\ast }$ &  \\ 
&  &  & $H^{2}(G/T)$ & 
\end{tabular}%
,
\end{quote}

\noindent the induced map $\pi ^{\ast }$ is an isomorphism by the $1$%
--connectness of the pair $(G,T)$. The \textsl{Borel} \textsl{transgression}
in the fibration (2.1) is $\tau =$ $(\pi ^{\ast })^{-1}\circ \delta $ (\cite[%
p.185]{Mc}).

\bigskip

\noindent \textbf{Lemma 2.1.} \textsl{The diffeomorphism type of the flag
manifold }$G/T$\textsl{\ depends only on the semi--simple part }$G^{\prime }$
\textsl{of the group }$G$\textsl{\ as}

\begin{enumerate}
\item[(2.4)] $G/T\cong \frac{G_{1}}{T_{1}}\times \cdots \times \frac{G_{k}}{%
T_{k}}$\textsl{\ (see (1.2)),}
\end{enumerate}

\noindent \textsl{where }$T_{i}$\textsl{\ is a maximal torus of the }$1$%
\textsl{--connected simple Lie group }$G_{i}$\textsl{, }$1\leq i\leq k$%
\textsl{.}

\textsl{The transgression }$\tau $ \textsl{fits into the following exact
sequence, where }$Tor(A)$\textsl{\ denotes the torsion subgroup of an
abelian group }$A$\textsl{,}

\begin{enumerate}
\item[(2.5)] $0\rightarrow H^{1}(G)\overset{j^{\ast }}{\rightarrow }H^{1}(T)%
\overset{\tau }{\rightarrow }H^{2}(G/T)\overset{\pi ^{\ast }}{\rightarrow }%
TorH^{2}(G)\rightarrow 0$.
\end{enumerate}

\noindent \textbf{Proof.} Let $T^{\prime }$ be a maximal torus of the
semi--simple part $G^{\prime }$ of $G$. By (1.1) a maximal torus of $G$ is $%
T=\left[ T^{\prime }\times \mathcal{Z}_{0}(G)\right] /F$. The diffeomorphism
(2.4) comes from the obvious relation $G/T=G^{\prime }/T^{\prime }$ and
(1.2).

Since the second homotopy group of a Lie group is trivial, the homotopy
exact sequence of $\pi $ contains the free resolution of the fundamental
group $\pi _{1}(G)$

\begin{quote}
$0\rightarrow \pi _{2}(G/T)\rightarrow \pi _{1}(T)\rightarrow \pi
_{1}(G)\rightarrow 0$.
\end{quote}

\noindent Applying the co--functor $Hom(,\mathbb{Z})$ to this sequence, and
using the Huriwicz isomorphisms $\pi _{2}(G/T)=H_{2}(G/T)$, $\pi
_{1}(T)=H_{1}(T)$, $\pi _{1}(G)=H_{1}(G)$ to substitute for the relevant
groups, one obtains (2.5).$\square $

\bigskip

By (1.3) we have $H^{1}(G)=H^{1}(T^{r})$, $r=\dim \mathcal{Z}_{0}(G)$. By
(2.5) the transgression $\tau $ annihilates the direct summand $H^{1}(T^{r})$
of $H^{1}(T)$ and hence, depends only on the semi--simple part of $G$. For
this reason we can assume below that the Lie group $G$ under consideration
is semi--simple.

Equip the Lie algebra $L(G)$ of $G$ with an inner product $(,)$ so that the
adjoint representation acts as isometries on $L(G)$. Let $L(T)\subset L(G)$
be the \textsl{Cartan subalgebra} corresponding to the fixed maximal torus $%
T $ on $G$, and fix a set $\Delta =\{\alpha _{1},\cdots ,\alpha
_{n}\}\subset L(T)$ of simple roots of $G$, where $n=\dim T$.

The Euclidean space $L(T)$ contains three distinguished lattices. Firstly,
the set $\{\alpha _{1},\cdots ,\alpha _{n}\}$ of simple roots generates 
\textsl{the} \textsl{root lattice} $\Lambda _{r}$ of $G$. Next, the
pre--image of the exponential map $\exp :L(T)\rightarrow T$ at the group
unit $e\in T$ gives rise to the \textsl{unit lattice }$\Lambda _{e}:=\exp
^{-1}(e)$\textsl{\ }of $G$. Finally, using simple roots one defines the set $%
\Omega =\{\phi _{1},\cdots ,\phi _{n}\}$ of \textsl{fundamental dominant
weights} of $G$ by the formula $2(\phi _{i},\phi _{j})/(\alpha _{j},\alpha
_{j})=\delta _{i,j}$ that generates \textsl{the weight lattice} $\Lambda
_{\omega }$ of $G$.

Let $A=(b_{ij})_{n\times n}$,\ $b_{ij}=2(a_{i},\alpha _{j})/(\alpha
_{j},\alpha _{j})$, be the Cartan matrix of $G$, and let $A^{\tau }$\ be the
transpose of $A$. The following result can be found in \cite[(3.4)]{DL}.

\bigskip

\noindent \textbf{Lemma 2.2. }\textsl{On the space }$L(T)$ \textsl{one has }$%
\Lambda _{r}\subseteq \Lambda _{e}\subseteq \Lambda _{\omega }$\textsl{.} 
\textsl{In addition}

\textsl{i) the group }$G$\textsl{\ is }$1$\textsl{--connected if and only if 
}$\Lambda _{r}=\Lambda _{e}$\textsl{;}

\textsl{ii) the group }$G$\textsl{\ is adjoint if and only if }$\Lambda
_{e}=\Lambda _{\omega }$\textsl{;}

\textsl{iii) the basis }$\Delta $\textsl{\ on }$\Lambda _{r}$\textsl{\ can
be expressed by the basis }$\Omega $\textsl{\ on }$\Lambda _{\omega }$ 
\textsl{by the formula}

\begin{quote}
$\qquad \left( \alpha _{1},\cdots ,\alpha _{n}\right) =\left( \phi
_{1},\cdots ,\phi _{n}\right) \cdot A^{\tau }$.$\square $
\end{quote}

For a root $\alpha \in \Delta $ let $K(\alpha )\subset G$ be the subgroup
with Lie algebra $l_{\alpha }\oplus L_{\alpha }$, where $l_{\alpha }\subset
L(T)$ is the $1$--dimensional subspace spanned by $\alpha $, and $L_{\alpha
}\subset L(G)$ is the root space (viewed as an oriented real $2$--plane)
belonging to the root $\alpha $ (\cite[p.35]{H}). Then the circle subgroup $%
S^{1}=\exp (l_{\alpha })$ is a maximal torus on $K(\alpha )$, while quotient
space $K_{\alpha }/S^{1}$ is diffeomorphic to the $2$--dimensional sphere $%
S^{2}$. Moreover, the inclusion $(K_{\alpha },S^{1})\subset (G,T)$ induces
an embedding

\begin{enumerate}
\item[(2.6)] $s_{\alpha }:S^{2}=K_{\alpha }/S^{1}\rightarrow G/T$
\end{enumerate}

\noindent whose image is known as the \textsl{Schubert variety} associated
to the root $\alpha $ \cite{DZ0}. By the basis theorem of Chevalley \cite{Ch}
the maps $s_{\alpha }$ with $\alpha \in \Delta $ represent a basis of the
second homology $H_{2}(G/T)$. As a result if one lets $\omega _{i}\in
H^{2}(G/T)$\ be the Kronnecker dual of the homology class represented by the
map $s_{\alpha _{i}}$, then

\bigskip

\noindent \textbf{Lemma 2.3. }\textsl{The set }$\left\{ \omega _{1},\cdots
,\omega _{n}\right\} $\textsl{\ is a basis of the group }$H^{2}(G/T)$.$%
\square $

\bigskip

On the other hand let $\Theta =\{\theta _{1},\cdots ,\theta _{n}\}$ be a
basis for the unit lattice $\Lambda _{e}$. It defines $n$ oriented circle
subgroups on the maximal torus

\begin{enumerate}
\item[(2.7)] $\widetilde{\theta }_{i}:S^{1}=\mathbb{R}/\mathbb{Z}\rightarrow
T$, $\widetilde{\theta }_{i}(t):=\exp (t\theta _{i})$, $1\leq i\leq n$,
\end{enumerate}

\noindent that represent also a basis of the first homology $H_{1}(T)$. As
result if we let $t_{i}\in H^{1}(T)$ be the class Kronnecker dual to the map 
$\widetilde{\theta }_{i}$, then

\begin{enumerate}
\item[(2.8)] $H^{\ast }(T)=\Lambda (t_{1},\cdots ,t_{n})$ (i.e. the exterior
ring generated by $t_{1},\cdots ,t_{n}$).
\end{enumerate}

\noindent Let $C(\Theta )=\left( c_{i,j}\right) _{n\times n}$\ be the matrix
expressing the ordered basis $\Delta $ by the ordered basis $\Theta $ in
view of the inclusion $\Lambda _{r}\subseteq \Lambda _{e}$ by Lemma 2.2.
Namely, $\left( \alpha _{1},\cdots ,\alpha _{n}\right) =\left( \theta
_{1},\cdots ,\theta _{n}\right) C(\Theta )$.

\bigskip

\noindent \textbf{Theorem 2.4.}\textsl{\ With respect to the basis\ (2.6)
and (2.8) on the groups }$H^{2}(G/T)$\textsl{\ and } $H^{1}(T)$\textsl{,\
the transgression }$\tau $ \textsl{is given by the formula}

\begin{enumerate}
\item[(2.9)] $\left( \tau (t_{1}),\cdots ,\tau (t_{n})\right) =\left( \omega
_{1},\cdots ,\omega _{n}\right) C(\Theta )$\textsl{.}
\end{enumerate}

\noindent \textbf{Proof. }Assume firstly\textbf{\ }that the group $G$ is $1$%
--connected. Then a basis $\Theta $ of the unit lattice $\Lambda
_{e}=\Lambda _{r}$ can be taken to be $\Delta =\{\alpha _{1},\cdots ,\alpha
_{n}\}$. Since $C(\Theta )$ is then the identity matrix we are bound to show
that $\tau (t_{i})=\omega _{i}$, $1\leq i\leq n$.

For each root $\alpha _{i}\in $ $\Delta $ the inclusion $(K(\alpha
_{i}),S^{1})\subset (G,T)$ induces the following bundle map over $s_{\alpha
_{i}}$, where by the $1$--connectness of the group $G$ the group $K(\alpha
_{i})$ is isomorphic to the $3$--sphere $S^{3}$ and $\pi _{i}$ is the Hopf
fibration over $S^{2}$,

\begin{quote}
$\qquad \qquad 
\begin{array}{ccccc}
S^{1} & \rightarrow & K(\alpha _{i})=S^{3} & \overset{\pi _{i}}{\rightarrow }
& K(\alpha _{i})/S^{1}=S^{2} \\ 
\widetilde{\alpha }_{1}\downarrow \qquad &  & \downarrow &  & \downarrow
s_{\alpha _{i}} \\ 
T & \rightarrow & G & \overset{\pi }{\rightarrow } & G/T\qquad%
\end{array}%
$
\end{quote}

\noindent This indicates that in the homotopy exact sequence of $\pi _{i}$
the connecting homomorphism $\partial $ satisfies $\partial \left[ \iota _{2}%
\right] =\left[ \iota _{1}\right] $, where $\iota _{r}$ is the identity on
the $r$--sphere. By the naturality of $\partial $ one gets in the homotopy
exact sequence of $\pi $ that $\partial \left[ s_{\alpha _{i}}\right] =\left[
\widetilde{\alpha }_{i}\right] $. This shows $\tau (t_{i})=$ $\omega _{i}$
as the map $\tau $ is dual to $\partial $ by the proof of Lemma 2.1.

In a general case let $d:(G_{0},T_{0})\rightarrow (G,T)$ be the universal
covering of $G$ with $T_{0}$ the maximal torus on $G_{0}$ corresponding to $%
T $. Then

\begin{quote}
$\exp =d\circ \exp _{0}:(L(G_{0}),L(T_{0}))\rightarrow
(G_{0},T_{0})\rightarrow (G,T)$.
\end{quote}

\noindent where $\exp $ (resp. $\exp _{0}$) is the exponential map of the
group $G$ (resp. $G_{0}$). It follows that, if we let $p(\Lambda
_{r},\Lambda _{e}):$ $T_{0}=L(T_{0})/\Lambda _{r}$ $\rightarrow
T=L(T_{0})/\Lambda _{e}$ be the covering map induced by the inclusion $%
\Lambda _{r}\subset \Lambda _{e}$ of the lattices, then

\begin{enumerate}
\item[(2.10)] $d\mid T_{0}=p(\Lambda _{r},\Lambda _{e}):T_{0}\rightarrow T$,
\end{enumerate}

\noindent and the induced map $p(\Lambda _{r},\Lambda _{e})_{\ast }$ on $\pi
_{1}(T_{0})$ is

\begin{enumerate}
\item[(2.11)] $p(\Lambda _{r},\Lambda _{e})_{\ast }[\widetilde{\alpha }%
_{i}]=c_{i,1}[\widetilde{\theta }_{1}]+\cdots +c_{i,n}[\widetilde{\theta }%
_{n}]$ with $C(\Theta )=\left( c_{ij}\right) _{n\times n}$.
\end{enumerate}

On the other hand the restriction $d\mid T_{0}$ fits in the commutative
diagram

\begin{enumerate}
\item[(2.12)] $%
\begin{array}{ccc}
\pi _{2}(G_{0}/T_{0}) & \underset{\cong }{\overset{\partial _{0}}{%
\rightarrow }} & \pi _{1}(T_{0}) \\ 
\parallel &  & \qquad \downarrow (d\mid T_{0})_{\ast } \\ 
\pi _{2}(G/T) & \overset{\partial }{\rightarrow } & \pi _{1}(T)%
\end{array}%
$
\end{enumerate}

\noindent with $\partial _{0}$, $\partial $ the connecting homomorphisms in
the homotopy exact sequences of the bundles $G_{0}\rightarrow G_{0}/T_{0}$, $%
G\rightarrow G/T$, respectively, where the vertical identification on the
left comes from (2.1). It follows that, for each root $\alpha _{i}\in \Delta 
$,

\begin{quote}
$\partial \left[ s_{\alpha _{i}}\right] =(d\mid T_{0})_{\ast }\circ \partial
_{0}\left[ s_{\alpha _{i}}\right] $ (by the diagram (2.12))

$\qquad =(d\mid T_{0})_{\ast }\left[ \widetilde{\alpha }_{i}\right] $ (by
the proof of the previous case)

$\qquad =p(\Lambda _{r},\Lambda _{e})_{\ast }(\left[ \widetilde{\alpha }_{i}%
\right] )$ (by (2.10)).
\end{quote}

\noindent The proof is completed by (2.11), and by the fact that the map $%
\tau $ is dual to $\partial $.$\square $

\bigskip

\noindent \textbf{Example 2.5.} Formula (2.9) is ready to evaluate the
transgression $\tau $, hence the differential $d_{2}$ on $E_{2}^{\ast ,\ast
}(G)$. As examples we have by Lemma 2.2 that

i) if the group $G$ is $1$--connected one can take $\Delta =\{\alpha
_{1},\cdots ,\alpha _{n}\}$ as a basis for the unit lattice $\Lambda _{e}$,
and the transition matrix $C(\Theta )$ is the identity;

ii) if the group $G$ is of the adjoint type, then the set $\Omega =\{\phi
_{1},\cdots ,\phi _{n}\}$ of fundamental dominant weights is a basis of $%
\Lambda _{e}$, and the corresponding transition matrix $C(\Theta )$ from $%
\Lambda _{e}$ to $\Lambda _{r}$ is\textsl{\ }the transpose $A^{\tau }$ of
the Cartan matrix $A$.$\square $

\subsection{Schubert presentation of the ring $H^{\ast }(G/T)$}

Turning to a concise presentation of the factor subring $H^{\ast }(G/T)$ of $%
E_{2}^{\ast ,\ast }(G)$ assume that the rank of the semi--simple part of $G$
is $n$, and let $\{\omega _{1},\cdots ,\omega _{n}\}$ be the Schubert basis
on $H^{2}(G/T)$. It is shown in \cite[Theorem 1.2]{DZ1} that

\bigskip

\noindent \textbf{Theorem 2.6. }\textsl{For each Lie group }$G$\textsl{\
there exist a set }$\left\{ y_{1},\cdots ,y_{m}\right\} $\textsl{\ of
Schubert classes on }$G/T$\textsl{\ with }$\deg y_{i}$\textsl{\ }$>2$\textsl{%
,\ so that the set }$\{\omega _{1},\cdots ,\omega _{n},y_{1},\cdots ,y_{m}\}$%
\textsl{\ is a minimal system of generators of the integral cohomology ring }%
$H^{\ast }(G/T)$\textsl{.}

\textsl{With respect to these generators the ring }$H^{\ast }(G/T)$ \textsl{%
has the presentation}

\begin{enumerate}
\item[(2.13)] $H^{\ast }(G/T)=\mathbb{Z}[\omega _{1},\cdots ,\omega
_{n},y_{1},\cdots ,y_{m}]/\left\langle h_{i},f_{j},g_{j}\right\rangle
_{1\leq i\leq k;1\leq j\leq m}$\textsl{,}
\end{enumerate}

\noindent \textsl{in which}

\begin{quote}
\textsl{i) for each }$1\leq i\leq k$\textsl{, }$h_{i}\in \left\langle \omega
_{1},\cdots ,\omega _{n}\right\rangle $\textsl{;}

\textsl{ii) for each }$1\leq j\leq m$\textsl{, the pair }$(f_{j},g_{j})$%
\textsl{\ of polynomials is related to the Schubert class }$y_{j}$\textsl{\
in the fashion }

$\qquad f_{j}$\textsl{\ }$=$\textsl{\ }$p_{j}y_{j}+\alpha _{j}$\textsl{,
\quad }$g_{j}=y_{j}^{k_{j}}+\beta _{j}$\textsl{, }

\textsl{where }$p_{j}\in \{2,3,5\}$\textsl{\ and }$\alpha _{j},\beta _{j}\in
\left\langle \omega _{1},\cdots ,\omega _{n}\right\rangle $\textsl{.}$%
\square $
\end{quote}

\noindent \textbf{Example 2.7.} For the cases $G=SU(n),Sp(n)$, $E_{6}$ and $%
E_{7}$ we refer to Theorem 5.1 for the explicit form of the formula (2.13)
of the rings $H^{\ast }(G/T)$.$\square $

\section{Construction and computation in $E_{3}^{\ast ,\ast }(G)$}

By the fiber degrees the third page of the spectral sequence $\{E_{r}^{\ast
,\ast }(G),d_{r}\}$ has the decomposition $E_{3}^{\ast ,\ast
}(G)=E_{3}^{\ast ,0}(G)\oplus E_{3}^{\ast ,1}(G)\oplus \cdots \oplus
E_{3}^{\ast ,N}(G)$, $N=\dim T$. Based on Theorems 2.4 and 2.6 we single out
certain elements in the initial two summands $E_{3}^{\ast ,0}$ and $%
E_{3}^{\ast ,1}$, which will show to generate the ring $H^{\ast }(G)$
multiplicatively. Along the way we demonstrate that, with respect to our
explicit constructions of the generators on the ring $H^{\ast }(G)$, the
Bockstein homomorphism, Steenrod operators, as well as the homomorphism $%
\theta $ in (1.8), can be effectively evaluated by simple formulae, see
Section 3.4 and Theorem 3.10.

\subsection{The term $E_{3}^{\ast ,0}(G)$}

The formula (2.3) of the differential $d_{2}:E_{2}^{\ast ,\ast
}(G)\rightarrow E_{2}^{\ast ,\ast }(G)$ implies that $E_{3}^{\ast
,0}(G)=H^{\ast }(G/T)/\left\langle \func{Im}\tau \right\rangle $. From
(2.13) one gets

\bigskip

\noindent \textbf{Lemma 3.1. }$E_{3}^{\ast ,0}(G)=H^{\ast }(G/T)\mid _{\tau
(t_{1})=\cdots =\tau (t_{N})=0}$\textsl{.}$\square $

\subsection{Constructions in $E_{3}^{\ast ,1}(G)$}

For a $d_{2}$--cocycle $\gamma \in E_{2}^{\ast ,\ast }(G)$ write $[\gamma
]\in E_{3}^{\ast ,\ast }(G)$ for its cohomology class. Based on Theorems 2.4
and 2.6, we present two ways to construct elements in $E_{3}^{\ast ,1}(G)$.
The first one resorts to $\ker \tau $, while the second utilizes $\func{Im}%
\tau $.

For a $t\in \ker \tau $ the element $1\otimes t\in E_{2}^{0,1}$ is clear a $%
d_{2}$--cocycle. The class $\iota (t):=[1\otimes t]\in E_{3}^{0,1}(G)$ will
be called a \textsl{primary} $\QTR{sl}{1}$\textsl{--form} of $G$ with base
degree $0$. Moreover, in view of the exact sequence (2.5) one has

\bigskip

\noindent \textbf{Lemma 3.2. }\textsl{The map }$\iota :\ker \tau \rightarrow 
$\textsl{\ }$E_{3}^{0,1}(G)=H^{1}(G)$\textsl{\ is an isomorphism.}$\square $

\bigskip

For the second construction we take firstly the ring $\mathbb{Z}$ of
integers as coefficient for cohomologies. In view of the formula (2.13) one
has the surjective ring map

\begin{quote}
$f:\mathbb{Z}[\omega _{i},y_{j}]_{1\leq i\leq n,1\leq j\leq m}\rightarrow
H^{\ast }(G/T)$ with $\ker f=$ $\left\langle h_{i},f_{j},g_{j}\right\rangle $%
.
\end{quote}

\noindent Since $f$ is an isomorphism in degree $2$ the transgression $\tau $
has a unique lift $\widetilde{\tau }$ into the free polynomial ring $\mathbb{%
Z}[\omega _{i},y_{j}]$ subject to the relation $\tau =f\circ \widetilde{\tau 
}$. For a polynomial $P\in \left\langle \func{Im}\widetilde{\tau }%
\right\rangle $ we can write

\begin{enumerate}
\item[(3.1)] $P=p_{1}\cdot \widetilde{\tau }(t_{1})+\cdots +p_{N}\cdot 
\widetilde{\tau }(t_{N})$ with $p_{i}\in \mathbb{Z}[\omega _{i},y_{j}]$,
\end{enumerate}

\noindent and set $\widetilde{P}:=f(p_{1})\otimes t_{1}+\cdots +$ $%
f(p_{N})\otimes t_{N}\in E_{2}^{\ast ,1}(G)$. It is crucial for us to note
from $f(P)=d_{2}(\widetilde{P})$ (by (2.3)) that

\begin{quote}
"$P\in \left\langle \func{Im}\widetilde{\tau }\right\rangle \cap \ker f$ 
\textsl{implies} $\widetilde{P}\in \ker d_{2}$".
\end{quote}

\noindent Most importantly, one can show that

\bigskip

\noindent \textbf{Lemma 3.3.} \textsl{The map }$\varphi :\left\langle \func{%
Im}\widetilde{\tau }\right\rangle \cap \ker f\rightarrow E_{3}^{\ast ,1}(G)$%
\textsl{\ by }$\varphi (P)=[\widetilde{P}]$\textsl{\ is well defined. }

\bigskip

\noindent \textbf{Proof.}\textsl{\ }We are bound to show the class $[%
\widetilde{P}]\in E_{3}^{\ast ,1}(G)$ is independent of a choice of the
expansion (3.1). Assume in addition to (3.1) that one has a second expansion 
$P=h_{1}\cdot \widetilde{\tau }(t_{1})+\cdots +h_{N}\cdot \widetilde{\tau }%
(t_{N})$. Then the equation

\begin{quote}
$(p_{1}-h_{1})\cdot \widetilde{\tau }(t_{1})+\cdots +(p_{N}-h_{N})\cdot 
\widetilde{\tau }(t_{N})=0$
\end{quote}

\noindent holds in the ring $\mathbb{Z}[\omega _{i},y_{j}]$. We can assume
below that $p_{1}-h_{1}\neq 0$.

\noindent \textbf{Case 1.} The set $\{\widetilde{\tau }(t_{1}),\cdots ,%
\widetilde{\tau }(t_{N})\}$\ is a basis of $\func{Im}\widetilde{\tau }$:
Since $\{\widetilde{\tau }(t_{1}),\cdots ,\widetilde{\tau }(t_{N})\}\subset 
\mathbb{Z}[\omega _{i},y_{j}]$ is algebraically independent, the above
equation implies that all the differences $p_{i}-h_{i}$ with $i\neq 1$ are
divisible by $\widetilde{\tau }(t_{1})$. That is $p_{i}-h_{i}=q_{i}\cdot 
\widetilde{\tau }(t_{1})$ for some $q_{i}\in \mathbb{Z}[\omega _{i},y_{j}]$, 
$2\leq i\leq N$. The proof is done by the calculation

\begin{quote}
$d_{2}(f(q_{2})\otimes t_{1}t_{2}+\cdots +f(q_{N})\otimes t_{1}t_{N})$

$=f(p_{1}-h_{1})\otimes t_{1}+\cdots +f(p_{N}-h_{N})\otimes t_{N}$.
\end{quote}

\noindent \textbf{Case 2.} The set $\{\widetilde{\tau }(t_{1}),\cdots ,%
\widetilde{\tau }(t_{N})\}$\ is linearly dependent $\func{Im}\widetilde{\tau 
}$:\textbf{\ }Take a subset $\{\overline{t}_{1},\cdots ,\overline{t}%
_{n^{\prime }}\}\subset H^{1}(T)$ ($n^{\prime }\leq N$) so that its $%
\widetilde{\tau }$--image is a basis of the group $\func{Im}\widetilde{\tau }
$. Let $B=(b_{ij})_{n^{\prime }\times N}$ be the matrix expressing the
elements $\tau (t_{i})$ by the basis elements $\{\widetilde{\tau }(\overline{%
t}_{1}),\cdots ,\widetilde{\tau }(\overline{t}_{n^{\prime }})\}$. Denote by $%
B^{\tau }$ the transpose of $B$ and set

\begin{quote}
$(p_{1}^{\prime },\cdots ,p_{n^{\prime }}^{\prime })=(p_{1},\cdots
,p_{N})\cdot B^{\tau }$, $(h_{1}^{\prime },\cdots ,h_{n^{\prime }}^{\prime
})=(h_{1},\cdots ,h_{N})\cdot B^{\tau }$.
\end{quote}

\noindent Then, in addition to the next two expansions of $P$ in $\mathbb{Z}%
[\omega _{i},y_{j}]_{1\leq i\leq n,1\leq j\leq m}$

\begin{quote}
$P=p_{1}^{\prime }\cdot \widetilde{\tau }(\overline{t}_{1})+\cdots
+p_{n^{\prime }}^{\prime }\cdot \widetilde{\tau }(\overline{t}_{n^{\prime
}})=h_{1}^{\prime }\cdot \widetilde{\tau }(\overline{t}_{1})+\cdots
+h_{n^{\prime }}^{\prime }\cdot \widetilde{\tau }(\overline{t}_{n^{\prime
}}) $,
\end{quote}

\noindent one has the following relations in $E_{2}^{\ast ,1}(G)$

\begin{quote}
$f(p_{1})\otimes t_{1}+\cdots +$ $f(p_{N})\otimes t_{N}=f(p_{1}^{\prime
})\otimes \overline{t}_{1}+\cdots +$ $f(p_{n^{\prime }}^{\prime })\otimes 
\overline{t}_{n^{\prime }}$,

$f(h_{1})\otimes t_{1}+\cdots +$ $f(h_{N})\otimes t_{N}=f(h_{1}^{\prime
})\otimes \overline{t}_{1}+\cdots +$ $f(h_{n^{\prime }}^{\prime })\otimes 
\overline{t}_{n^{\prime }}$.
\end{quote}

\noindent The proof is done by the following computation in $E_{3}^{\ast
,1}(G)$

\begin{quote}
$\left[ f(p_{1})\otimes t_{1}+\cdots +f(p_{N})\otimes t_{N}\right] =\left[
f(p_{1}^{\prime })\otimes \overline{t}_{1}+\cdots +f(p_{n^{\prime }}^{\prime
})\otimes \overline{t}_{n^{\prime }}\right] $

$=\left[ f(h_{1}^{\prime })\otimes \overline{t}_{1}+\cdots +f(h_{n^{\prime
}}^{\prime })\otimes \overline{t}_{n^{\prime }}\right] =\left[
f(h_{1})\otimes t_{1}+\cdots +f(h_{N})\otimes t_{N}\right] $,
\end{quote}

\noindent where the second equality has been shown in Case 1.$\square $

\bigskip

The maps $\widetilde{\tau }$ and $\varphi $ above has its analogue for
cohomology over a finite field $\mathbb{F}_{p}$. Precisely, from Theorem 2.6
one can deduce a presentation of the ring $H^{\ast }(G/T;\mathbb{F}_{p})$ in
the following form (\cite[Lemma 2.3]{DZ3})

\begin{enumerate}
\item[(3.2)] $H^{\ast }(G/T;\mathbb{F}_{p})=\mathbb{F}_{p}[\omega
_{1},\ldots ,\omega _{n},y_{t}]/\left\langle \delta _{1},\cdots ,\delta
_{n},y_{t}^{k_{t}}+\sigma _{t}\right\rangle _{\text{ }t\in E(G,p)}$,
\end{enumerate}

\noindent where $\delta _{s}\in \mathbb{F}_{p}[\omega _{1},\ldots ,\omega
_{n}]$, $\sigma _{t}\in \left\langle \omega _{1},\ldots ,\omega
_{n}\right\rangle $, $E(G,p)=\{1\leq t\leq m;p_{t}=p\}$. Based on (3.2) one
formulates the $\mathbb{F}_{p}$--analogies of the maps $f$, $\widetilde{\tau 
}$ and $\varphi $ as

\begin{quote}
$f_{p}:\mathbb{F}_{p}[\omega _{1},\ldots ,\omega _{n},y_{t}]_{\text{ }t\in
E(G,p)}\rightarrow H^{\ast }(G/T;\mathbb{F}_{p})$,

$\widetilde{\tau }_{p}:H^{1}(T;\mathbb{F}_{p})\rightarrow \mathbb{F}%
_{p}[\omega _{1},\ldots ,\omega _{n},y_{t}]_{\text{ }t\in E(G,p)}$,

$\varphi _{p}:\left\langle \func{Im}\widetilde{\tau }_{p}\right\rangle \cap
\left\langle \ker f_{p}\right\rangle \rightarrow E_{3}^{\ast ,1}(G;\mathbb{F}%
_{p})$.
\end{quote}

\noindent The proof of Lemma 3.3 is applicable to show that

\bigskip

\noindent \textbf{Lemma 3.4.} \textsl{The correspondence }$\varphi _{p}$%
\textsl{\ is well defined. In particular}

\begin{quote}
$\varphi _{p}(P^{\prime }P)=0$\textsl{\ if }$P^{\prime }\in \left\langle 
\func{Im}\widetilde{\tau }_{p}\right\rangle $\textsl{\ and }$P\in \ker f_{p}$%
\textsl{.}$\square $
\end{quote}

\subsection{Extension from $E_{3}^{\ast ,\ast }(G)$ to $H^{\ast }(G)$}

In preparation to solve the extension problem from $E_{3}^{\ast ,\ast }(G)$
to $H^{\ast }(G)$ we let $F^{p}$ be the filtration on $H^{\ast }(G)$ defined
by the fibration (2.1). That is

\begin{center}
$0=F^{r+1}(H^{r}(G))\subseteq F^{r}(H^{r}(G))\subseteq \cdots \subseteq
F^{0}(H^{r}(G))=H^{r}(G)$
\end{center}

\noindent with $E_{\infty }^{p,q}(G)=F^{p}(H^{p+q}(G))/F^{p+1}(H^{p+q}(G))$.
The routine relation $d_{r}(E_{r}^{\ast ,0}(G))=0$ for $r\geq 2$ yields the
sequence of quotient maps

\begin{center}
$H^{r}(G/T)=E_{2}^{r,0}\rightarrow E_{3}^{r,0}\rightarrow \cdots \rightarrow
E_{\infty }^{r,0}=F^{r}(H^{r}(G))\subset H^{r}(G)$
\end{center}

\noindent whose composition agrees with the induces map $\pi ^{\ast
}:H^{\ast }(G/T)\rightarrow H^{\ast }(G)$ \cite[P.147]{Mc}. For this reason
we can reserve $\pi ^{\ast }$ also for the composition

\begin{enumerate}
\item[(3.3)] $\pi ^{\ast }:E_{3}^{\ast ,0}(G)\rightarrow \cdots \rightarrow
E_{\infty }^{\ast ,0}(G)=F^{r}(H^{r}(G))\subset H^{\ast }(G)$.
\end{enumerate}

The property $H^{odd}(G/T)=0$ by Theorem 3.1 indicates that $E_{r}^{odd,q}=0$
for $r,q\geq 0$. This implies that $%
F^{2k+1}(H^{2k+1}(G))=F^{2k+2}(H^{2k+1}(G))=0$ and that

\begin{quote}
$E_{\infty }^{2k,1}(G)=F^{2k}(H^{2k+1}(G))\subset H^{2k+1}(G)$.
\end{quote}

\noindent Combining this with $d_{r}(E_{r}^{\ast,1})=0$ for $r\geq3$ yields
the composition

\begin{enumerate}
\item[(3.4)] $\kappa :E_{3}^{\ast ,1}(G)\rightarrow \cdots \rightarrow
E_{\infty }^{\ast ,1}(G)\subset H^{2t+1}(G)$
\end{enumerate}

\noindent that interprets directly elements of $E_{3}^{\ast ,1}$ as
cohomology classes of the group $G$.

\bigskip

\noindent \textbf{Definition 3.5.} For a polynomial $P\in \ker f\cap
\left\langle \func{Im}\widetilde{\tau }\right\rangle $ (resp. $P\in \ker
f_{p}\cap \left\langle \func{Im}\widetilde{\tau }_{p}\right\rangle $) we
shall refer to the class $\kappa \varphi (P)\in H^{\ast }(G)$ (resp. $\kappa
\varphi _{p}(P)\in H^{\ast }(G;\mathbb{F}_{p})$) as \textsl{the primary }$1$%
\textsl{--form on }$G$\textsl{\ with} \textsl{characteristic polynomial} $P$.%
$\square $

\bigskip

\noindent \textbf{Example 3.6.} If the group $G$ is $1$--connected then $%
\left\langle \func{Im}\widetilde{\tau }\right\rangle =\left\langle \omega
_{1},\cdots ,\omega _{n}\right\rangle $ (resp. $\left\langle \func{Im}%
\widetilde{\tau }_{p}\right\rangle =\left\langle \omega _{1},\cdots ,\omega
_{n}\right\rangle $) by Theorem 2.4, and the presentation (2.13) implies
that the set of polynomials

\begin{quote}
$S(G)=\{h_{i},p_{j}\beta _{j}-y_{j}^{k_{j}}\alpha _{j}\mid 1\leq i\leq
k,1\leq j\leq m\}$

(resp. $S_{p}(G):=\left\{ \delta _{1},\cdots ,\delta _{n}\right\} $)
\end{quote}

\noindent belongs to $\ker f\cap \left\langle \func{Im}\widetilde{\tau }%
\right\rangle $ (resp. to $\ker f_{p}\cap \left\langle \func{Im}\widetilde{%
\tau }_{p}\right\rangle $). It will be called \textsl{a set of primary
characteristic polynomials of }$G$\textsl{\ over} $\mathbb{Z}$ (resp. 
\textsl{over} $\mathbb{F}_{p}$).

It has been shown in \cite{DZ2,DZ3} that

\textsl{i) the square free products of the primary }$1$\textsl{--forms }$%
\kappa \circ \varphi (P)$\textsl{\ with }$P\in S(G)$\textsl{\ is a basis of
the free part of the integral cohomology }$H^{\ast }(G)$\textsl{;}

\textsl{ii) the ring }$H^{\ast }(G;\mathbb{F}_{p})$\textsl{\ is generated by
the set }$\left\{ \kappa \circ \varphi _{p}(\delta )\mid \delta \in
S_{p}(G)\right\} $\textsl{\ of primary }$1$\textsl{--forms, together with }$%
\pi ^{\ast }(y_{t})$\textsl{, }$t\in E(G,p)$\textsl{.}

\noindent In Section 4 these results will act as the input for computing the
cohomologies of the adjoint Lie groups $PG$.$\square $

\subsection{Steenrod operations on $H^{\ast }(G;\mathbb{F}_{p})$}

Let $\mathcal{A}_{{p}}$ be the \textsl{Steenrod algebra} with $\mathcal{P}%
^{k}\in \mathcal{A}_{{p}}$, $k\geq 1$, the $k^{th}$ \textsl{reduced power}
(if $p=2$ it is also customary to write $Sq^{2k}$ instead of $\mathcal{P}%
^{k} $), and with $\delta _{p}=r_{p}\circ \beta _{p}\in \mathcal{A}_{{p}}$
the \textsl{Bockstein operator} \cite{SE}. We can reduce the $\mathcal{A}_{{p%
}}$ action on the primary $1$--forms $\kappa \varphi _{p}(P)$ to computation
with the characteristic polynomials $P$.

Since the Koszul complex $E_{2}^{\ast ,\ast }(G)$ is torsion free one has
for a prime $p$ the short exact sequence of complexes

\begin{quote}
$0\rightarrow $ $E_{2}^{\ast ,\ast }(G)\overset{\cdot p}{\rightarrow }%
E_{2}^{\ast ,\ast }(G)\overset{r_{p}}{\rightarrow }E_{2}^{\ast ,\ast }(G;%
\mathbb{F}_{p})\rightarrow 0$.
\end{quote}

\noindent With respect to the maps $\kappa $ and $\pi ^{\ast }$ the
connecting homomorphism $\widehat{\beta }_{p}$ of the associated cohomology
exact sequence clearly satisfies the commutative diagram

\begin{enumerate}
\item[(3.5)] 
\begin{tabular}{lll}
$E_{3}^{\ast ,1}(G;\mathbb{F}_{p})$ & $\overset{\widehat{\beta }_{p}}{%
\rightarrow }$ & $E_{3}^{\ast ,0}(G)$ \\ 
$\kappa \downarrow $ &  & $\pi ^{\ast }\downarrow $ \\ 
$H^{\ast }(G;\mathbb{F}_{p})$ & $\overset{\beta _{p}}{\rightarrow }$ & $%
H^{\ast }(G)$%
\end{tabular}%
.
\end{enumerate}

\noindent In addition, by Lemma 3.1 the quotient map $H^{\ast
}(G/T)\rightarrow E_{3}^{\ast ,0}(G)$ is

\begin{quote}
$f(P)\rightarrow f(P)\mid _{\tau (t_{1})=\cdots =\tau (t_{N})=0}$, $P\in 
\mathbb{Z}[\omega _{i},y_{j}]$.
\end{quote}

\noindent Let $P_{0}\in \left\langle \func{Im}\widetilde{\tau }\right\rangle 
$ be an integral lift of a polynomial $P\in \ker f_{p}\cap \left\langle 
\func{Im}\widetilde{\tau }_{p}\right\rangle $ (i.e. $P_{0}\equiv P\func{mod}%
p $).\textsl{\ }The diagram chasing

\begin{quote}
$%
\begin{array}{ccccc}
&  & \varphi (P_{0}) & \overset{r_{p}}{\rightarrow } & \varphi _{p}(P) \\ 
&  & d_{2}\downarrow \quad &  & d_{2}\downarrow \quad \\ 
\frac{1}{p}f(P_{0}) & \overset{\cdot p}{\longrightarrow } & f(P_{0}) &  & 0%
\end{array}%
$
\end{quote}

\noindent in above short exact sequence shows that

\bigskip

\noindent \textbf{Lemma 3.7.} $\beta _{p}(\kappa \circ \varphi _{p}(P))=\pi
^{\ast }(\frac{1}{p}f(P_{0})\mid _{\tau (t_{1})=\cdots =\tau (t_{N})=0})$.$%
\square $

\bigskip

Let $c:$ $(G_{0},T_{0})\rightarrow (G,T)$ be the universal covering of a
semi--simple Lie group $G$, and consider the fibration induced by the
inclusion $i:T_{0}\rightarrow G_{0}$

\begin{enumerate}
\item[(3.6)] $G/T\overset{\psi }{\hookrightarrow }BT_{0}\overset{Bi}{%
\rightarrow }BG_{0}$,
\end{enumerate}

\noindent where $BT_{0}$ (resp. $BG_{0}$) is the classifying space of the
group $T_{0}$ (resp. $G_{0}$). It is known that the compositions $\widetilde{%
s}_{\alpha }:=\psi \circ s_{\alpha }:S^{2}\rightarrow BT_{0}$ with $\alpha
\in \Delta $ represent a basis of the group $H_{2}(BT_{0})$, where $%
s_{\alpha }$ are the maps defined by (2.6). As a result we can also write $%
\{\omega _{1},\ldots ,\omega _{n}\}$ for the basis on $H^{2}(BT_{0})$
Kronnecker dual to the ordered basis $\{\widetilde{s}_{\alpha }\mid \alpha
_{i}\in \Delta \}$ on $H_{2}(BT_{0})$. In this sense $\psi ^{\ast }=%
\overline{f}_{p}$, where $\overline{f}_{p}$ denotes the restriction of $%
f_{p} $ on the subalgebra $H^{\ast }(BT_{0};\mathbb{F}_{p})=\mathbb{F}%
_{p}[\omega _{1},\ldots ,\omega _{n}]$. Since the lift $\widetilde{\tau }%
_{p} $ of the transgression $\tau $ takes values in $H^{\ast }(BT_{0};%
\mathbb{F}_{p})$ one has the subspace

\begin{quote}
$\ker \overline{f}_{p}\cap \left\langle \func{Im}\widetilde{\tau }%
_{p}\right\rangle \subset H^{\ast }(BT_{0};\mathbb{F}_{p})$
\end{quote}

\noindent which is clearly closed under the $\mathcal{A}_{{p}}$ action on $%
H^{\ast }(BT_{0};\mathbb{F}_{p})$. The proof of \cite[Lemma 3.2]{DZ3} is
applicable to show the following formula, that reduces the $\mathcal{P}^{k}$
action on $H^{\ast }(G;\mathbb{F}_{p})$ to that on the much simpler $%
\mathcal{A}_{{p}}$--algebra $H^{\ast }(BT_{0};\mathbb{F}_{p})$.

\bigskip

\noindent \textbf{Lemma 3.8.} \textsl{For} \textsl{a characteristic} \textsl{%
polynomial} $P\in \ker \overline{f}_{p}\cap \left\langle \func{Im}\widetilde{%
\tau }_{p}\right\rangle $ \textsl{one has}

\begin{enumerate}
\item[(3.7)] $\mathcal{P}^{k}(\kappa \circ \varphi _{p}(P))=\kappa \circ
\varphi _{p}(\mathcal{P}^{k}(P))$.$\square $
\end{enumerate}

\subsection{A refinement of the exact sequence (1.8)}

Returning to the situation concerned by Theorem 1.3 let $G$ be a
semi--simple Lie group $G$ whose center $\mathcal{Z}(G)$ contains the cyclic
group $\mathbb{Z}_{q}$. Then the circle bundle $C$ on the quotient group $G/%
\mathbb{Z}_{q}$ fits into the commutative diagram

\begin{enumerate}
\item[(3.8)] $%
\begin{array}{ccccc}
S^{1} & \hookrightarrow & \left[ T\times S^{1}\right] /\mathbb{Z}_{q} & 
\overset{C^{\prime }}{\rightarrow } & T^{\prime } \\ 
\parallel &  & \cap \quad &  & \cap \quad \\ 
S^{1} & \hookrightarrow & \left[ G\times S^{1}\right] /\mathbb{Z}_{q} & 
\overset{C}{\rightarrow } & G/\mathbb{Z}_{q} \\ 
&  & \pi ^{\prime }\downarrow \quad &  & \pi \downarrow \quad \\ 
&  & G/T & = & G/T%
\end{array}%
$,
\end{enumerate}

\noindent where $T\subset G$ is a fixed maximal torus on $G$, $T^{\prime
}:=T/\mathbb{Z}_{q}$, the vertical maps $\pi ^{\prime }$ and $\pi $ are the
obvious quotients by the maximal torus, and where $C^{\prime }$ denotes the
restriction of $C$ to $\left[ T\times S^{1}\right] /\mathbb{Z}_{q}$ . Since
the maximal torus $\left[ T\times S^{1}\right] /\mathbb{Z}_{q}$\ of $\left[
G\times S^{1}\right] /\mathbb{Z}_{q}$ has the factorization $T^{\prime
}\times S^{1}$ so that $C^{\prime }$\ is the projection onto the first
factor (by the proof of Theorem 1.2), one can take a basis $\Theta =\{\theta
_{1},\cdots ,\theta _{n},\theta _{0}\}$ ($n=\dim T^{\prime }$) for the unit
lattice of the group $\left[ G\times S^{1}\right] /\mathbb{Z}_{q}$, so that
the tangent map of $C^{\prime }$ at the group unit $e$ carries the subset $%
\{\theta _{1},\cdots ,\theta _{n}\}$ to a basis of the unit lattice of $G/%
\mathbb{Z}_{q}$. As a result if we let $\{t_{1},\cdots ,t_{n},t_{0}\}$ be
the basis of $H^{1}(\left[ T\times S^{1}\right] /\mathbb{Z}_{q})$
corresponding to $\Theta $ in the manner of (2.8), then

\begin{quote}
a) $C^{\prime \ast }$ maps $H^{\ast }(T^{\prime })$ isomorphically onto the
subring $\Lambda ^{\ast }(t_{1},\cdots ,t_{n})$ of $H^{\ast }(\left[ T\times
S^{1}\right] /\mathbb{Z}_{q})=\Lambda ^{\ast }(t_{1},\cdots ,t_{n},t_{0})$;

b) the transgression $\tau $ in $\pi $ is the restriction of $\tau ^{\prime
} $ on $H^{1}(T^{\prime })$.
\end{quote}

\noindent Summarizing, the bundle map $C$ from $\pi ^{\prime }$ to $\pi $
fits into the short exact sequence

\begin{enumerate}
\item[(3.9)] $0\rightarrow E_{2}^{\ast ,k}(G/\mathbb{Z}_{q})\overset{C^{\ast
}}{\rightarrow }E_{2}^{\ast ,k}(\left[ G\times S^{1}\right] /\mathbb{Z}_{q})%
\overset{\overline{\theta }}{\rightarrow }E_{2}^{\ast ,k-1}(G/\mathbb{Z}%
_{q})\rightarrow 0$,
\end{enumerate}

\noindent where the quotient map $\overline{\theta }$ has the following
description: if $x\in H^{\ast }(G/T)$, $y=y_{0}+t_{0}\cdot y_{1}\in H^{\ast
}(\left[ T\times S^{1}\right] /\mathbb{Z}_{q})$ with $y_{0},y_{1}\in \Lambda
^{\ast }(t_{1},\cdots ,t_{n})$, then

\begin{enumerate}
\item[(3.10)] $\overline{\theta }(x\otimes y)=x\otimes y_{1}$.
\end{enumerate}

\noindent It follows from (3.9) that one has an exact sequence of the form

\begin{enumerate}
\item[(3.11)] ${\small \cdots \rightarrow }E_{3}^{\ast ,r}(G/Z_{q})\overset{%
{\small C}^{\ast }}{{\small \rightarrow }}E_{3}^{\ast ,r}(\left[ G\times
S^{1}\right] /Z_{q})\overset{\overline{\theta }}{{\small \rightarrow }}%
E_{3}^{\ast ,r-1}(G/Z_{q})\overset{{\small \varpi }}{{\small \rightarrow }}%
E_{3}^{\ast ,r-1}(G/Z_{q})${\tiny \ }$\overset{{\small C}^{\ast }}{{\small %
\rightarrow }}{\small \cdots }$
\end{enumerate}

\noindent in which the map $\varpi $\ is induced by the endomorphism\textsl{%
\ }$x\otimes y\rightarrow (x\cup \varpi )\otimes y$ on $E_{2}^{\ast ,\ast
}(G/\mathbb{Z}_{q})$ with $\varpi :=\tau ^{\prime }(t_{0})$.

With respect to the multiplicative structure inherited from that on $%
E_{2}^{\ast ,\ast }$ the third page $E_{3}^{\ast ,\ast }$ is a bi-graded
ring \cite[P.668]{Wh}. Let $J(\varpi )$ and $\left\langle \varpi
\right\rangle $ be respectively the subring and the ideal of $E_{3}^{\ast
,\ast }(G/\mathbb{Z}_{q})$ generated by the class $\varpi \in E_{3}^{2,0}(G/%
\mathbb{Z}_{q})$. Write $E_{3}^{\ast ,\ast }(G/\mathbb{Z}_{q})_{\left\langle
\varpi \right\rangle }$ for the quotient ring $E_{3}^{\ast ,\ast }(G/\mathbb{%
Z}_{q})/\left\langle \varpi \right\rangle $ with quotient map $g$. Then, in
addition to the exact sequence

\begin{enumerate}
\item[(3.12)] $0\rightarrow \left\langle \varpi \right\rangle \rightarrow
E_{3}^{\ast ,\ast }(G/\mathbb{Z}_{q})\overset{g}{\rightarrow }E_{3}^{\ast
,\ast }(G/\mathbb{Z}_{q})_{\left\langle \varpi \right\rangle }\rightarrow 0$
\end{enumerate}

\noindent the proof of Theorem 1.3 is valid to show the following result. By
Lemma 3.2 the group $E_{3}^{0,1}(\left[ G\times S^{1}\right] /\mathbb{Z}%
_{q}) $ has an element that corresponds to the class $\xi _{1}$ defined by
(1.6), which we denote still by $\xi _{1}$.

\bigskip

\noindent \textbf{Theorem\textbf{\ }3.9.} \textsl{The induced map }$C^{\ast
} $ \textsl{fits in} \textsl{the exact sequence}

\begin{enumerate}
\item[(3.13)] $0\rightarrow E_{3}^{\ast ,\ast }(G/\mathbb{Z}%
_{q})_{\left\langle \varpi \right\rangle }\overset{C^{\ast }}{\rightarrow }%
E_{3}^{\ast ,\ast }([G\times S^{1}]/\mathbb{Z}_{q})\overset{\overline{\theta 
}}{\rightarrow }E_{3}^{\ast ,\ast }(G/\mathbb{Z}_{q})\overset{\varpi }{%
\rightarrow }\left\langle \varpi \right\rangle \rightarrow 0$
\end{enumerate}

\noindent \textsl{where}

\textsl{i)} $\overline{\theta }(\xi _{1})=q\in E_{3}^{0,0}(G/\mathbb{Z}_{q})$%
\textsl{;}

\textsl{ii)} $\overline{\theta }(x\cup C^{\ast }(y))=\overline{\theta }%
(x)\cup y$\textsl{,} $x\in E_{3}^{\ast ,\ast }(\left[ G\times S^{1}\right] /%
\mathbb{Z}_{q})$\textsl{,} $y\in E_{3}^{\ast ,\ast }(G/\mathbb{Z}_{q})$%
\textsl{,}

\textsl{iii)} \textsl{the class} $\varpi $ \textsl{satisfies }$\pi ^{\ast
}(\varpi )=\omega $ \textsl{(see (3.3)).}

\textsl{In addition, if} \textsl{the map }$g$\textsl{\ admits a split
homomorphism }$j:E_{3}^{\ast ,\ast }(G/\mathbb{Z}_{q})_{\left\langle \varpi
\right\rangle }$ $\rightarrow $ $E_{3}^{\ast ,\ast }(G/\mathbb{Z}_{q})$%
\textsl{, then the map }

\begin{quote}
$h:J(\varpi )\otimes E_{3}^{\ast ,\ast }(G/\mathbb{Z}_{q})_{\left\langle
\varpi \right\rangle }\rightarrow E_{3}^{\ast ,\ast }(G/\mathbb{Z}_{q})$
\end{quote}

\noindent \textsl{by} $h(\varpi ^{r}\otimes x)=\varpi ^{r}\cup j(x)$ \textsl{%
induces an isomorphism of }$J(\varpi )$\textsl{--modules}

\begin{enumerate}
\item[(3.14)] $E_{3}^{\ast ,\ast }(G/\mathbb{Z}_{q})\cong \frac{J(\varpi
)\otimes E_{3}^{\ast ,\ast }(G/\mathbb{Z}_{q})_{\left\langle \varpi
\right\rangle }}{\left\langle \varpi \cdot \func{Im}\overline{\theta }%
\right\rangle }$\textsl{.}$\square $
\end{enumerate}

\noindent \textbf{Proof. }It suffices to show the relations i), ii) and iii).%
\textbf{\ }By the choice of the class $\xi _{1}\in E_{3}^{0,1}(\left[
G\times S^{1}\right] /\mathbb{Z}_{q})$ and by the following commutative
diagram

\begin{quote}
\begin{tabular}{lll}
$E_{3}^{0,1}(\left[ G\times S^{1}\right] /\mathbb{Z}_{q})$ & $\overset{%
\overline{\theta }}{\rightarrow }$ & $E_{3}^{0,0}(G/\mathbb{Z}_{q})$ \\ 
$\quad \parallel $ &  & $\quad \parallel $ \\ 
$H^{1}(\left[ G\times S^{1}\right] /\mathbb{Z}_{q})$ & $\overset{\theta }{%
\rightarrow }$ & $H^{0}(G/\mathbb{Z}_{q})$%
\end{tabular}%
,
\end{quote}

\noindent property i) corresponds to i) of Theorem 1.3, while ii) comes
directly from the formula (3.10). Finally, the relation iii) is shown by the
commutative diagram induced by the map $\pi ^{\ast }$ in (3.3):

\begin{quote}
\begin{tabular}{lll}
$E_{3}^{0,0}(G/\mathbb{Z}_{q})$ & $\overset{\cup \varpi }{\rightarrow }$ & $%
E_{3}^{2,0}(G/\mathbb{Z}_{q})$ \\ 
$\quad \parallel $ &  & $\pi ^{\ast }\quad \downarrow $ \\ 
$H^{0}(G/\mathbb{Z}_{q})$ & $\overset{\cup \omega }{\rightarrow }$ & $%
H^{2}(G/\mathbb{Z}_{q})$%
\end{tabular}%
.$\square $
\end{quote}

The exact sequence (3.13) can be seen to be a refinement of the sequence
(1.8). In addition, the maps $\pi ^{\ast }$ and $\kappa $ in (3.3) and (3.4)
build up the obvious commutative diagram relating the operator $\theta $ in
(1.8) with the map $\overline{\theta }$ in (3.13)

\begin{enumerate}
\item[(3.15)] 
\begin{tabular}{lll}
$E_{3}^{2r,1}(\left[ G\times S^{1}\right] /\mathbb{Z}_{q})$ & $\overset{%
\overline{\theta }}{\rightarrow }$ & $E_{3}^{2r,0}(G/\mathbb{Z}_{q})$ \\ 
$\qquad \kappa \quad \downarrow $ &  & $\pi ^{\ast }\downarrow $ \\ 
$H^{2r+1}(G\times S^{1})$ & $\overset{\theta }{\rightarrow }$ & $H^{2r}(G/%
\mathbb{Z}_{q})$%
\end{tabular}%
.
\end{enumerate}

\noindent Moreover, the map $\overline{\theta }$ (hence $\theta $) admits a
simple formula we come to describe.

For\ a polynomial $P\in \left\langle \func{Im}\widetilde{\tau ^{\prime }}%
\right\rangle $\textsl{\ }we set $P_{0}=P\mid _{\widetilde{\tau }%
(t_{1})=\cdots =\widetilde{\tau }(t_{n})=0}$. Then $P-P_{0}\in \left\langle 
\func{Im}\widetilde{\tau }\right\rangle $ and $P_{0}$ is divisible by $%
\varpi =\widetilde{\tau ^{\prime }}(t_{0})$. This enables us to define the 
\textsl{derivation} $\partial P/\partial \varpi $ of $P$ with respect to $%
\varpi $ by the formula

\begin{enumerate}
\item[(3.16)] $\partial P/\partial \varpi :=P_{0}/\varpi $,
\end{enumerate}

\noindent and to get an expansion of $P$ in the form (note that $\widetilde{%
\tau ^{\prime }}(t_{i})=\widetilde{\tau }(t_{i})$, $i\leq n$, by b))

\begin{quote}
$P=p_{1}\cdot \widetilde{\tau ^{\prime }}(t_{1})+\cdots +$ $p_{n}\cdot 
\widetilde{\tau ^{\prime }}(t_{n})+\partial P/\partial \varpi \cdot 
\widetilde{\tau ^{\prime }}(t_{0})$.
\end{quote}

\noindent In term of Lemma 3.3 for a polynomial $P\in \ker f\cap
\left\langle \func{Im}\widetilde{\tau ^{\prime }}\right\rangle $ we have

\begin{quote}
$\varphi (P)=\left[ f(p_{1})\otimes t_{1}+\cdots +f(p_{n})\otimes
t_{n}+f(\partial P/\partial \varpi )\otimes t_{0}\right] $.
\end{quote}

\noindent The diagram (3.15), together with the formula (3.10), concludes
that

\bigskip

\noindent \textbf{Theorem 3.10.} \textsl{The map} $\overline{\theta }$ 
\textsl{in (3.13) (resp. the map} $\theta $ \textsl{in (1.8))} \textsl{%
satisfies that}

\begin{enumerate}
\item[(3.17)] $\overline{\theta }(\varphi (P))=f(\partial P/\partial \varpi
) $ \textsl{(resp.} $\theta (\kappa \circ \varphi (P))=\pi ^{\ast
}f(\partial P/\partial \varpi )$\textsl{),}
\end{enumerate}

\noindent \textsl{where} $P\in \ker f\cap \left\langle \func{Im}\widetilde{%
\tau ^{\prime }}\right\rangle $, $P_{0}=P\mid _{\widetilde{\tau }%
(t_{1})=\cdots =\widetilde{\tau }(t_{n})=0}$\textsl{.}$\square $

\bigskip

Let $P\in \ker f\cap \left\langle \func{Im}\widetilde{\tau ^{\prime }}%
\right\rangle $ be a polynomial with\textsl{\ }$\overline{\theta }(\varphi
(P))=0$. By the exact sequence (3.13) there exists a $1$--form $\eta \in
E_{3}^{\ast ,1}(G/\mathbb{Z}_{q})$ satisfying $C^{\ast }(\eta )=\varphi (P)$%
. The proof of the next result indicates an algorithm to construct from $P$
a characteristic polynomial $P^{\prime }\in \ker f\cap \left\langle \func{Im}%
\widetilde{\tau }\right\rangle $ for such a class $\eta $.

\bigskip

\noindent \textbf{Lemma 3.11. }\textsl{For a polynomial }$P\in \ker f\cap
\left\langle \func{Im}\widetilde{\tau ^{\prime }}\right\rangle $ \textsl{%
with\ }$\overline{\theta }(\varphi (P))=0$\textsl{,} \textsl{there exists a
polynomial }$P^{\prime }\in \ker f\cap \left\langle \func{Im}\widetilde{\tau 
}\right\rangle $ \textsl{such that }$C^{\ast }\varphi (P^{\prime })=\varphi
(P)$.

\bigskip

\noindent \textbf{Proof.} With $\overline{\theta }(\varphi (P))=0$ we have
by the exact sequence (3.13) that $d_{2}(\gamma )=f(\partial P/\partial
\varpi )$ for some $\gamma \in E_{2}^{\ast ,1}(G/\mathbb{Z}_{q})$. We can
assume further that $\gamma =\widetilde{H}$ for some $H=h_{1}\cdot 
\widetilde{\tau }(t_{1})+\cdots +h_{n}\cdot \widetilde{\tau }(t_{n})\in
\left\langle \func{Im}\widetilde{\tau }\right\rangle $. The desired
polynomial $P^{\prime }$ is given by $P^{\prime }:=P+(H-\partial P/\partial
\varpi )\cdot \varpi $.

Indeed, the obvious relation $H-\partial P/\partial \varpi \in \ker f$
implies that

\begin{quote}
$(H-\partial P/\partial \varpi )\cdot \varpi \in \left\langle \func{Im}%
\widetilde{\tau ^{\prime }}\right\rangle \cap \ker f$.
\end{quote}

\noindent Consequently, $P^{\prime }\in \left\langle \func{Im}\widetilde{%
\tau }\right\rangle \cap \ker f$. The relation $C^{\ast }(\varphi (P^{\prime
}))=\varphi (P)$ on $E_{3}^{\ast ,1}([G\times $ $S^{1}]/\mathbb{Z}_{q})$ is
verified by the following calculation in $E_{2}^{\ast ,1}(\left[ G\times
S^{1}\right] /\mathbb{Z}_{q})$

\begin{quote}
$d_{2}(\dsum\limits_{1\leq i\leq n}f(h_{i})\otimes t_{i}t_{0})=-f(\partial
P/\partial \varpi )\otimes t_{0}+\dsum\limits_{1\leq i\leq n}f(h_{i})\varpi
\otimes t_{i}$

$=\widetilde{P^{\prime }}-\widetilde{P}$.$\square $
\end{quote}

\section{The cohomology of adjoint Lie groups}

A Lie group $G$ is called \textsl{adjoint} if its center subgroup $\mathcal{Z%
}(G)$ is trivial. In particular, for any Lie group $G$ the quotient group $%
PG:=G/\mathcal{Z}(G)$ is adjoint. Granted with the constructions and
formulae developed in Section 3 we compute in this section the cohomologies
of the adjoint Lie groups $PG$ for $G=SU(n),Sp(n),E_{6},E_{7}$. In these
cases the quotient maps $c:G\rightarrow PG$ are cyclic. The corresponding
circle bundle over $PG$ is denoted by

\begin{quote}
$C:\left[ G\times S^{1}\right] /\mathbb{Z}_{q}\cong G\times S^{1}\rightarrow
PG$, $q=\left\vert \mathcal{Z}(G)\right\vert $,
\end{quote}

\noindent where the diffeomorphism $\cong $ follows from Theorem 1.2.

Briefly, the computation is carried out by three steps. Starting from the
formula (2.13) of the ring $H^{\ast }(G/T)$ we obtain firstly the subgroups $%
\func{Im}\pi ^{\ast }$ and $\func{Im}\kappa $. The exact sequence (3.13) is
then applied to formulate the additive cohomology $H^{\ast }(PG)$ by $\func{%
Im}\pi ^{\ast }$ and $\func{Im}\kappa $. Finally, the structure of $H^{\ast
}(PG)$ as a ring is determined by expressing the squares $x^{2}$ with $x\in 
\func{Im}\kappa $ as elements of $\func{Im}\pi ^{\ast }$, see \cite[Lemma 2.8%
]{DZ2}.

With the group $G$ being $1$--connected the cohomologies $H^{\ast }(G;%
\mathbb{F}_{p})$ are known \cite{DZ2}. The following result allows us to
exclude the cohomology $H^{\ast }(PG;\mathbb{F}_{p})$ with $(p,q)=1$ from
further consideration.

\bigskip

\noindent \textbf{Theorem 4.1.} \textsl{If }$(p,q)=1$\textsl{\ the map} $C$ 
\textsl{induces a ring isomorphism }

\begin{quote}
$C^{\ast }:H^{\ast }(PG;\mathbb{F}_{p})\cong H^{\ast }(G;\mathbb{F}_{p})$.
\end{quote}

\noindent \textbf{Proof. }With\textbf{\ }$p$ co--prime to $q$ one has $%
\omega \equiv 0\func{mod}p$. Therefore the exact sequence (1.7) becomes

\begin{quote}
$0\rightarrow H^{\ast }(PG;\mathbb{F}_{p})\overset{C^{\ast }}{\rightarrow }%
H^{\ast }(G\times S^{1};\mathbb{F}_{p})\overset{\theta }{\rightarrow }%
H^{\ast }(PG;\mathbb{F}_{p})\rightarrow 0$.
\end{quote}

\noindent Moreover, the relation $\theta (\xi _{1})=1$ by i) of Theorem 1.3
implies that, with respect to the decomposition $H^{\ast }(G;\mathbb{F}%
_{p})\oplus \xi _{1}\cdot H^{\ast }(G;\mathbb{F}_{p})$ on $H^{\ast }(G\times
S^{1};\mathbb{F}_{p})$ the map $C^{\ast }$ carries $H^{\ast }(PG;\mathbb{F}%
_{p})$ injectively into the first summand $H^{\ast }(G;\mathbb{F}_{p})$,
while by ii) of Theorem 1.3 the map $\theta $ maps the second summand $\xi
_{1}\cdot H^{\ast }(G;\mathbb{F}_{p})$ surjectively onto $H^{\ast }(PG;%
\mathbb{F}_{p})$. This establishes the isomorphism.$\square $

\subsection{The map $C^{\ast }:E_{3}^{\ast ,0}(G/\mathbb{Z}_{q})$ $%
\rightarrow E_{3}^{\ast ,0}(G\times S^{1})$ in (3.11)}

We begin by taking integers as coefficients for cohomology. Since the groups 
$G$ are $1$--connected $TorH^{2}(G\times S^{1})=0$ by Theorem 1.2.
Therefore, the transgression $\tau ^{\prime }$ in $\pi ^{\prime }$ (see
(3.8)) is surjective by (2.2). One gets from Theorem 5.1, as well the
formula $E_{3}^{\ast ,0}(G\times S^{1})=H^{\ast }(G/T)\mid _{\omega
_{1}=\cdots =\omega _{n}=0}$by Lemma 3.1, that

\begin{enumerate}
\item[(4.1)] 
\begin{tabular}{l}
$E_{3}^{\ast ,0}(SU(n)\times S^{{\small 1}})=E_{3}^{\ast ,0}(Sp(n)\times
S^{1})=\mathbb{Z}$; \\ 
$E_{3}^{\ast ,0}(E_{6}\times S^{1})=\frac{\mathbb{Z}[x_{6},x_{8}]}{%
\left\langle 2x_{6},3x_{8},x_{6}^{2},x_{8}^{3}\right\rangle }$; \\ 
$E_{3}^{\ast ,0}(E_{7}\times S^{1})=\frac{\mathbb{Z}%
[x_{6},x_{8},x_{10},x_{18}]}{\left\langle 2x_{6},3x_{8},2x_{10},2{x_{18},}%
x_{6}^{2},x_{8}^{3},x_{10}^{2},x_{18}^{2}\right\rangle }$.%
\end{tabular}
\end{enumerate}

\noindent where $x_{i}$'s are the special Schubert classes on $E_{n}/T$, $%
n=6,7$, defined by (5.1).

Similarly, for the groups $PG$ take a set $\Omega =\{\phi _{1},\cdots ,\phi
_{m}\}$ of fundamental dominant weights as a basis for the unit lattice $%
\Lambda _{e}$ of $PG$, and let $\{t_{1},\cdots ,t_{m}\}$ be the
corresponding basis on the group $H^{1}(T^{\prime })$ (see in (3.8)), where $%
m=n-1,n,6$ or $7$ in accordance to $G=SU(n),Sp(n),E_{6},E_{7}$. By iii) of
Lemma 2.2 the matrix $C(\Omega )$ expressing the basis $\Delta $\ of the
root lattice $\Lambda _{r}$ by $\Omega $ is the Cartan matrix of $G$.
Granted with the presentation of the rings $H^{\ast }(G/T)$ in Theorem 5.1,
as well as the results of Lemma 5.2, the formula $E_{3}^{\ast
,0}(PG)=H^{\ast }(G/T)\mid _{\tau (t_{1})=\cdots =\tau (t_{m})=0}$by Lemma
3.1 yields that

\begin{enumerate}
\item[(4.2)] 
\begin{tabular}{l}
$E_{3}^{\ast ,0}(PSU(n))=\frac{\mathbb{Z}[\omega _{1}]}{\left\langle
b_{n,r}\omega _{1}^{r},\text{ }1\leq r\leq n\right\rangle }$ with $%
b_{n,r}=g.c.d.\{\binom{n}{1},\cdots ,\binom{n}{r}\}$; \\ 
$E_{3}^{\ast ,0}(PSp(n))=\frac{\mathbb{Z}[\omega _{1}]}{\left\langle 2\omega
_{1},\omega _{1}^{2^{r+1}}\right\rangle }$, $n=2^{r}(2s+1)$; \\ 
$E_{3}^{\ast ,0}(PE_{6})=\frac{\mathbb{Z}[\omega _{1},x_{3}^{\prime },x_{4}]%
}{\left\langle 3\omega _{1},2x_{3}^{\prime },3x_{4},x_{3}^{\prime 2},\omega
_{1}^{9},x_{4}^{3}\right\rangle }$, $x_{3}^{\prime }=x_{3}+\omega _{1}^{3}$;
\\ 
$E_{3}^{\ast ,0}(PE_{7})=\frac{\mathbb{Z}[\omega
_{2},x_{3},x_{4},x_{6},x_{9}]}{\left\langle 2\omega _{2},\omega
_{2}^{2},2x_{3},3x_{4},2x_{5},2{x_{9},}%
x_{3}^{2},x_{4}^{3},x_{5}^{2},x_{9}^{2}\right\rangle }$,%
\end{tabular}
\end{enumerate}

\noindent where $\binom{n}{r}:=\frac{n!}{r!(n-r)!}$. Inputting (4.1) and
(4.2) into the section

\begin{quote}
$E_{3}^{\ast ,0}(PG)\overset{\varpi }{\rightarrow }E_{3}^{\ast ,0}(PG)%
\overset{C^{\ast }}{\rightarrow }E_{3}^{\ast ,0}(G\times S^{1})\rightarrow 0$
\end{quote}

\noindent of the exact sequence (3.11) one obtains that

\bigskip

\noindent \textbf{Lemma 4.2.} \textsl{In the order of} $G=SU(n),Sp(n),E_{6}$ 
\textsl{and} $E_{7}$\textsl{\ one has}

\begin{quote}
\textsl{i)} $J(\varpi )=\frac{\mathbb{Z}[\varpi ]}{\left\langle
b_{n,r}\varpi ^{r},\text{ }1\leq r\leq n\right\rangle },\frac{\mathbb{Z}%
[\varpi ]}{\left\langle 2\varpi ,\varpi ^{2^{r+1}}\right\rangle },\frac{%
\mathbb{Z}[\varpi ]}{\left\langle 3\varpi ,\varpi ^{9}\right\rangle },\frac{%
\mathbb{Z}[\varpi ]}{\left\langle 2\varpi ,\varpi ^{2}\right\rangle }$%
\textsl{;}

\textsl{ii) the map }$C^{\ast }:$\textsl{\ }$E_{3}^{\ast ,0}(PG)$ $%
\rightarrow E_{3}^{\ast ,0}(G\times S^{1})$ \textsl{is given by }

$\qquad C^{\ast }(x_{s})=x_{s}$\textsl{, }$C^{\ast }(x_{3}^{\prime })=x_{3}$%
\textsl{;} $C^{\ast }(\varpi )=0$\textsl{.}
\end{quote}

\noindent \textsl{where }$\varpi =\omega _{1}$ \textsl{for} $%
G=SU(n),Sp(n),E_{6}\QTR{sl}{;}$ $\varpi =\omega _{2}\QTR{sl}{\ }$\textsl{for 
}$G=E_{7}.\square $

\bigskip

Since $H^{\ast }(G/T)$ is torsion free $E_{3}^{\ast ,0}(X;\mathbb{F}_{p})$ $%
=E_{3}^{\ast ,0}(X)\otimes \mathbb{F}_{p}$ for both $X=G/\mathbb{Z}_{q}$ and 
$G\times S^{1}$. Let $J_{p}(\varpi )\subset E_{3}^{\ast ,0}(PG;\mathbb{F}%
_{p})$ be the subring generated by $\varpi $ and set $E_{3}^{\ast ,0}(PG;%
\mathbb{F}_{p})_{\left\langle \varpi \right\rangle }=E_{3}^{\ast ,0}(PG;%
\mathbb{F}_{p})\mid _{\varpi =0}$. Lemma 4.2 implies that

\bigskip

\noindent \textbf{Lemma 4.3.} \textsl{In the order of }$(G,p)=(SU(n),p)$ 
\textsl{with }$n=p^{r}n^{\prime }$\textsl{\ and }$(p,n\prime )=1$\textsl{, }$%
(Sp(n),2)$\textsl{\ with }$n=2^{r}(2d+1)$\textsl{, }$(E_{6},3)$\textsl{\ and 
}$(E_{7},2)$\textsl{, one has}

\begin{quote}
\textsl{i)} $J_{p}(\varpi )=\frac{\mathbb{F}_{p}[\varpi ]}{\left\langle
\varpi ^{h(G)}\right\rangle }$ \textsl{with }$h(G)=p^{r},2^{r+1}$\textsl{, }$%
9$\textsl{, }$2$\textsl{;}

\textsl{ii) }$E_{3}^{\ast ,0}(PG;\mathbb{F}_{p})_{\left\langle \varpi
\right\rangle }=\mathbb{F}_{p},\mathbb{F}_{2},\frac{\mathbb{F}_{3}[x_{4}]}{%
\left\langle x_{4}^{3}\right\rangle },\frac{\mathbb{F}_{2}[x_{3},x_{5},x_{9}]%
}{\left\langle x_{3}^{2},x_{5}^{2},x_{9}^{2}\right\rangle }$\textsl{;}

\textsl{iii) }$E_{3}^{\ast ,0}(PG;\mathbb{F}_{p})=J_{p}(\varpi )\otimes
E_{3}^{\ast ,0}(PG;\mathbb{F}_{p})_{\left\langle \varpi \right\rangle }$,

\textsl{iv)} \textsl{the map }$C^{\ast }$ \textsl{annihilates }$\varpi $%
\textsl{\ and carries} \textsl{the factor} $1\otimes E_{3}^{\ast ,0}(PG;%
\mathbb{F}_{p})_{\left\langle \varpi \right\rangle }$ \textsl{in iii)} 
\textsl{isomorphically onto} the \textsl{ring }$E_{3}^{\ast ,0}(G\times
S^{1};\mathbb{F}_{p})$\textsl{.}$\square $
\end{quote}

\subsection{The ring $H^{\ast }(PG;\mathbb{F}_{p})$}

By Theorem 4.1 we can assume that $(G,p)=(SU(n),p)$ with $p\mid n$; $%
(Sp(n),2)$; $(E_{6},3)$; $(E_{7},2)$. In view of the set $S_{p}(G)$ of
primary characteristic polynomials of $G$ over $\mathbb{F}_{p}$ presented in
Table 5.1 (see also ii) of Example 3.7) one obtains the degree set $D(G,p)$
of the polynomials in $S_{p}(G)$ as that tabulated below:

\begin{center}
\textbf{Table 4.1.} The degree set of the primary characteristic polynomials
over $\mathbb{F}_{p}$

\begin{tabular}{l|l|l|l|l}
\hline\hline
$(G,p)$ & $(SU(n),p)$ & $(Spin(n),2)$ & $(E_{6},3)$ & $(E_{7},2)$ \\ \hline
$D(G,p)$ & $\{2,3,\cdots ,n\}$ & $\{2,4,\cdots ,2n\}$ & $\{2,4,5,6,8,9\}$ & $%
\{2,3,5,8,9,12,14\}$ \\ \hline\hline
\end{tabular}%
.
\end{center}

For each $s\in D(G,p)$ let $\xi _{2s-1}:=\varphi _{p}(P)\in E_{3}^{\ast
,1}(G\times S^{1},\mathbb{F}_{p})$, where $P\in S_{p}(G)$ with $\deg P=s$.
With the groups $G$ being $1$--connected one has by \cite[Theorem 5.4]{DZ3}
the ring isomorphism

\begin{quote}
$E_{3}^{\ast ,\ast }(G\times S^{1},\mathbb{F}_{p})=E_{3}^{\ast ,0}(G\times
S^{1},\mathbb{F}_{p})\otimes \Lambda (\xi _{1},\xi _{2s-1})_{s\in D(G,p)}$,
\end{quote}

\noindent where $\xi _{1}\in E_{3}^{0,1}$\ is the class specified in Theorem
3.9. Let $D(PG,p)$ be the complement of the number $h(G)$ (see i) of Lemma
4.3) in $D(G,p)$, put $B=E_{3}^{\ast ,0}(G\times S^{1},\mathbb{F}%
_{p})\otimes \Lambda (\xi _{1},\xi _{2s-1})_{s\in D(PG,p)}$ and rewritten

\begin{enumerate}
\item[(4.3)] $E_{3}^{\ast ,\ast }(G\times S^{1},\mathbb{F}_{p})=B\oplus \xi
_{2h(G)-1}\cdot B$ .
\end{enumerate}

\noindent Then the exact sequence (3.13) takes the following useful form

\begin{enumerate}
\item[(4.4)] $0\rightarrow E_{3}^{\ast ,\ast }(PG;\mathbb{F}%
_{p})_{\left\langle \varpi \right\rangle }\overset{C^{\ast }}{\rightarrow }%
B\oplus \xi _{2h(G)-1}\cdot B\overset{\overline{\theta }}{\rightarrow }%
E_{3}^{\ast ,\ast }(PG;\mathbb{F}_{p})\overset{\varpi }{\rightarrow }%
\left\langle \varpi \right\rangle \rightarrow 0$.
\end{enumerate}

With $p\mid q$ one has $\overline{\theta }(\xi _{1})=0$ by i) of Theorem
3.9. By (4.4) there is a unique class $\iota ^{\prime }\in E_{3}^{0,1}(G/%
\mathbb{Z}_{q};\mathbb{F}_{p})=\mathbb{F}_{p}$ such that

\begin{quote}
$C^{\ast }(\iota ^{\prime })=\xi _{1}$ (with $\kappa (\iota ^{\prime })=$ $%
\iota $, $\widehat{\beta }_{p}(\iota ^{\prime })=\varpi $ by the diagram
(3.5)).
\end{quote}

\noindent Furthermore, with the set $S_{p}(G)$ of characteristic polynomials
for $1$--forms $\xi _{2s-1}$ being presented in Table 5.1, the formula
(3.17) is applicable to evaluate $\overline{\theta }(\xi _{2s-1})$, $s\in
D(G,p)$. The results recorded in Table 5.2 yield that

\begin{enumerate}
\item[(4.5)] $\overline{\theta }(\xi _{2s-1})=0$ \textsl{if} $s\in D(PG,p)%
\QTR{sl}{,}$ $\omega ^{h(G)-1}$ \textsl{if} $s=h(G)$.
\end{enumerate}

\noindent With $\overline{\theta }(\xi _{2s-1})=0$ for $s\in D(PG,p)$ Lemma
3.11 assures us a polynomial $P^{\prime }\in \left\langle \func{Im}%
\widetilde{\tau }\right\rangle _{p}\cap \ker f_{p}$ with $\deg P^{\prime }=s$
so that class $\zeta _{2s-1}=\varphi _{p}(P^{\prime })$ satisfies

\begin{quote}
$C^{\ast }(\zeta _{2s-1})=\xi _{2s-1}$.
\end{quote}

\noindent A set of such polynomials $P^{\prime }$ so obtained, denoted by $%
S_{p}(PG)$, are presented in Table 5.3. Taking iv) of Lemma 4.3 into
consideration the exactness of the sequence (4.4) forces the isomorphisms

\begin{quote}
$E_{3}^{\ast ,\ast }(PG;\mathbb{F}_{p})_{\left\langle \omega \right\rangle
}\cong E_{3}^{\ast ,0}(PG;\mathbb{F}_{p})_{\left\langle \varpi \right\rangle
}\otimes \Lambda (\iota ^{\prime },\zeta _{2s-1})_{s\in D(PG,p)}$,

$\xi _{2h(G)-1}\cdot B\overset{\overline{\theta }}{\underset{\cong }{%
\rightarrow }}\varpi ^{h(G)-1}\cdot E_{3}^{\ast ,\ast }(PG;\mathbb{F}%
_{p})_{\left\langle \varpi \right\rangle }$ (by ii) of Theorem 3.9)
\end{quote}

\noindent which imply that $\func{Im}\overline{\theta }=\left\langle 
\overline{\omega }^{h(G)-1}\right\rangle $. With the field $\mathbb{F}_{p}$
as coefficient the exact sequence (3.12) is splitable. The formula (3.14)
then yields the presentation

\begin{enumerate}
\item[(4.6)] $E_{3}^{\ast ,\ast }(PG;\mathbb{F}_{p})=E_{3}^{\ast ,0}(PG;%
\mathbb{F}_{p})\otimes \Lambda (\iota ^{\prime },\zeta _{2s-1})_{s\in
D(PG,p)}$.
\end{enumerate}

By (4.6) the ring $E_{3}^{\ast ,\ast }(PG;\mathbb{F}_{p})$ is generated by $%
E_{3}^{\ast ,0}$ and $E_{3}^{\ast ,1}$, hence the differentials $d_{r}$ act
trivially on $E_{r}^{\ast ,\ast }(PG;\mathbb{F}_{p})$, $r\geq 4$. We get $%
E_{3}^{\ast ,\ast }(PG;\mathbb{F}_{p})=E_{\infty }^{\ast ,\ast }(PG;\mathbb{F%
}_{p})$. In particular, the maps $\pi ^{\ast }$ and $\kappa $ in (3.3) and
(3.4) are all monomorphisms. Combining this with the isomorphism of $\mathbb{%
F}_{p}$--spaces

\begin{quote}
$E_{\infty }^{\ast ,\ast }(PG;\mathbb{F}_{p})=H^{\ast }(PG;\mathbb{F}_{p})$
\end{quote}

\noindent we get from (4.6) the following additive presentation of the
cohomology $H^{\ast }(PG;\mathbb{F}_{p})$

\begin{enumerate}
\item[(4.7)] $H^{\ast }(PG;\mathbb{F}_{p})=\pi ^{\ast }E_{3}^{\ast ,0}(PG;%
\mathbb{F}_{p})\otimes \Delta (\iota ,\zeta _{2s-1})_{s\in D(PG,p)}$,
\end{enumerate}

\noindent where for simplicity the notion $\zeta _{2s-1}$ is reserved for
the cohomology classes $\kappa (\zeta _{2s-1})\in H^{\ast }(G;\mathbb{F}%
_{p}) $, and where one needs to replace exterior ring $\Lambda (\iota
^{\prime },\zeta _{2s-1})$ by the $\mathbb{F}_{p}$--module $\Delta (\iota
,\zeta _{2s-1})$, since the properties $\iota ^{\prime 2},\zeta
_{2s-1}^{2}=0 $ on $E_{\infty }^{\ast ,\ast }(PG;\mathbb{F}_{p})$ may not
survive to $H^{\ast }(G;\mathbb{F}_{p})$, see \cite[Lemma 2.8]{DZ2}.

Since $\zeta _{2s-1}=\kappa \circ \varphi _{p}(P)$ with $P\in S_{p}(PG)$ and 
$\deg P=s$, the formulae in Lemmas 3.7 and 3.8 is functional to compute $%
\beta _{p}(\zeta _{2s-1})$ and $Sq^{2r}(\zeta _{2s-1})$ in term of $P$. The
calculation justifying the following results can be found in Section 5.

\bigskip

\noindent \textbf{Lemma 4.4. }\textsl{With respect to the presentation (4.7)
one has}

\begin{quote}
\textsl{a) for }$(G,p)=(SU(n),p)$\textsl{\ with }$n=p^{r}n^{\prime }$, $%
(n\prime ,p)=1$

$\qquad \beta _{p}(\zeta _{2s-1})=-p^{r-t-1}\omega ^{p^{t}}$ \textsl{if} $%
s=p^{t}$ \textsl{with} $t<r$\textsl{,} $0$ \textsl{otherwise}

\qquad \textsl{(resp.} $Sq^{2s-2}\zeta _{2s-1}=\zeta _{4s-3}$ \textsl{for }$%
2s-1\leq 2^{r-1}$\textsl{);}

\textsl{b) for }$(G,p)=(Sp(n),2)$\textsl{\ with }$n=2^{r}(2b+1)$\textsl{:}

$\qquad \beta _{2}(\zeta _{4s-1})=\omega ^{2^{r}}$ \textsl{if} $s=2^{r-1}$%
\textsl{,} $0$ \textsl{if} $s\neq 2^{r-1}$

\qquad \textsl{(resp. }$Sq^{4s-2}\zeta _{4s-1}=\zeta _{8s-3}$\textsl{\ for }$%
4s-1\leq 2^{r}$\textsl{);}

\textsl{c) for }$(G,p)=(E_{6},3)$\textsl{\ and in the order of }$s=2,4,5,6,8$

$\qquad \beta _{3}(\zeta _{2s-1})=0,-x_{4},0,0,-x_{4}^{2}$\textsl{;}

\textsl{d) for }$(G,p)=(E_{7},2)$\textsl{\ and in the order of }$%
s=3,5,8,9,12,14$

$\qquad \beta _{2}(\zeta
_{2s-1})=x_{3},x_{5},x_{3}x_{5},x_{9},x_{3}x_{9},x_{5}x_{9}$

\qquad \textsl{(resp. }$Sq^{2s-2}\zeta _{2s-1}=\zeta _{9},\zeta
_{17},0,0,0,0 $\textsl{).}$\square $
\end{quote}

Combining (4.7) with Lemma 4.4 we show that

\bigskip

\noindent \textbf{Theorem 4.5.} \textsl{The rings }$H^{\ast }(PG;\mathbb{F}%
_{p})$ \textsl{has the following presentations}

\begin{enumerate}
\item[i)] $H^{\ast }(PSU(n);\mathbb{F}_{2})=\frac{\mathbb{F}_{2}[\omega ]}{%
\left\langle \omega ^{2^{r}}\right\rangle }\otimes \Delta (\iota )\otimes
\Lambda _{\mathbb{F}_{2}}(\zeta _{3},\zeta _{5},\cdots ,\widehat{\zeta }%
_{2^{r+1}-1},\cdots ,\zeta _{2n-1})$\textsl{,}

\textsl{where }$n=2^{r}(2b+1)$,\textsl{\ }$\iota ^{2}=\omega $ \textsl{or }$%
0 $\textsl{\ in accordance to }$r=1$\textsl{\ or} $r>1$\textsl{;}

\item[ii)] $H^{\ast }(PSU(n);\mathbb{F}_{p})=\frac{\mathbb{F}_{p}[\omega ]}{%
\left\langle \omega ^{p^{r}}\right\rangle }\otimes \Lambda _{\mathbb{F}%
_{p}}(\iota ,\zeta _{3},\cdots ,\widehat{\zeta }_{2p^{r}-1},\cdots ,\zeta
_{2n-1})$\textsl{, }

\textsl{where }$p\neq 2$\textsl{, }$n=p^{r}n^{\prime }$\textsl{\ with }$%
(n^{\prime },p)=1$\textsl{;}

\item[iii)] $H^{\ast }(PSp(n);\mathbb{F}_{2})=\frac{\mathbb{F}_{2}[\omega ]}{%
\left\langle \omega ^{2^{r+1}}\right\rangle }\otimes \Delta (\iota )\otimes
\Lambda _{\mathbb{F}_{2}}(\zeta _{3},\zeta _{7},\cdots ,\widehat{\zeta }%
_{2^{r+2}-1},\cdots ,\zeta _{4n-1})$\textsl{, }

\textsl{where }$\iota ^{2}=\omega $\textsl{,} $n=2^{r}(2b+1)$\textsl{;}

\item[iv)] $H^{\ast }(PE_{6};\mathbb{F}_{3})=\frac{\mathbb{F}_{3}[\omega
,x_{4}]}{\left\langle \omega ^{9},x_{4}^{3}\right\rangle }\otimes \Lambda _{%
\mathbb{F}_{3}}(\iota ,\zeta _{3},\zeta _{7},\zeta _{9},\zeta _{11},\zeta
_{15})$\textsl{;}

\item[v)] $H^{\ast }(PE_{7};\mathbb{F}_{2})=\frac{\mathbb{F}_{2}[\omega
,x_{3},x_{5},x_{9}]}{\left\langle \omega
^{2},x_{3}^{2},x_{5}^{2},x_{9}^{2}\right\rangle }\otimes \Delta (\iota
,\zeta _{5},\zeta _{9})\otimes \Lambda _{\mathbb{F}_{2}}(\zeta _{15},\zeta
_{17},\zeta _{23},\zeta _{27})$\textsl{,}

\textsl{where }$\iota ^{2}=\omega ,\zeta _{5}^{2}=x_{5},\zeta _{9}^{2}=x_{9}$%
\textsl{.}
\end{enumerate}

\noindent \textbf{Proof. }In the presentations through i) to v) the first
factor is $\func{Im}\pi ^{\ast }\cong E_{3}^{\ast ,0}(PG;\mathbb{F}_{p})$,
see iii) of Lemma 4.3. It remains to justify the expressions of the squares $%
\iota ^{2},\zeta _{2s-1}^{2}$ as that indicated in the theorem.

The cases ii) and iv) are trivial, as in a characteristic $p\neq 2$ the
square of any odd degree cohomology class is zero. In the remaining cases
i), iii) and v) we have $p=2$. The relations $\iota ^{2}=\omega $ come from $%
\omega =\beta _{q}(\iota )\in H^{2}(PG)$ with $q=n,2,2$ in accordance to $%
G=SU(n),Sp(n),E_{7}$. To evaluate $\zeta _{2s-1}^{2}$ with $s\in D(PG,p)$ we
make use of the Steenrod operators $Sq^{2r}$ by which

\begin{quote}
$\zeta _{2s-1}^{2}=\delta _{2}\circ Sq^{2s-2}(\zeta _{2s-1})$ (see \cite{SE}%
).
\end{quote}

\noindent Results in Lemma 4.4 implies that $\zeta _{2s-1}^{2}=0$\ with the
only exceptions $\zeta _{5}^{2}=x_{5}$, $\zeta _{9}^{2}=x_{9}$ when $%
(G,p)=(E_{7},2)$. This completes the proof.$\square $

\subsection{The integral cohomology of $PSU(n),PSp(n)$}

The integral cohomology of a space $X$ admits a canonical decomposition

\begin{enumerate}
\item[(4.8)] $H^{\ast }(X)=\mathcal{F}(X)\underset{p}{\oplus }\sigma _{p}(X)$
with $\sigma _{p}(X):=\{x\in H^{\ast }(X)\mid p^{r}x=0,$ $r\geq 1\}$,
\end{enumerate}

\noindent where the summands $\mathcal{F}(X)$ and $\sigma _{p}(X)$ are a 
\textsl{free part }and the $p$--\textsl{primary component }of $H^{\ast }(X)$%
, and where the sum is over all primes $p$. Therefore, the determination of
the cohomology $H^{\ast }(X)$ essentially consists of two tasks:

a) express $\mathcal{F}(X)$ and $\sigma _{p}(X)$ by explicit generators of
the ring $H^{\ast }(X)$;

b) decide the actions $\mathcal{F}(X)\times \sigma _{p}(X)\rightarrow \sigma
_{p}(X)$ of the free part on $\sigma _{p}(X)$.

\noindent Keeping these in mind we compute the integral cohomology $H^{\ast
}(PG)$ for $G=SU(n)$ and $SP(n)$ by applying the exact sequence (1.8) (resp.
(3.13)).

With the groups $G=SU(n)$ and $Sp(n)$ being $1$--connected the integral
cohomologies of the groups $G\times S^{1}$ are well known. Indeed, with
respect to the degree set $D(G)$ of the set $S(G)$ of primary characteristic
polynomials for $G$ over $\mathbb{Z}$ given in Table 4.4 (see also Example
3.6) one has

\begin{enumerate}
\item[(4.9)] $E_{3}^{\ast ,\ast }(G\times S^{1})=\Lambda (\gamma _{1},\gamma
_{2s-1})_{s\in D(G)}$($\cong H^{\ast }(G\times S^{1})$ via $\kappa $ in
(3.4)).
\end{enumerate}

\noindent where $D(SU(n))=\left\{ 2,\cdots ,n\right\} $, $D(Sp(n))=\left\{
4,\cdots ,2n\right\} $, and where $\gamma _{1}:=\xi _{1}$, $\gamma
_{2s-1}=\varphi (P)\in E_{3}^{\ast ,1}(G\times S^{1})$ with $P\in S(G)$ and $%
\deg P=s$. Moreover, in addition to the relations by i) of Theorem 3.9

\begin{quote}
$\overline{\theta }(\gamma _{1})=n$ or $2$ for $G=SU(n)$ or $Sp(n)$,
\end{quote}

\noindent with $\gamma _{2s-1}=\varphi (P)$ the formula (3.17) is applicable
to evaluate $\overline{\theta }(\gamma _{2s-1})$ in term of $P$. Explicitly,
the computation recorded in Table 5.5 tells that

\begin{enumerate}
\item[(4.10)] $\overline{\theta }(\gamma _{2s-1})=\left\{ 
\begin{tabular}{l}
$\binom{n}{s}\varpi ^{s-1}$ for $SU(n)$; \\ 
$0$ if $s\neq 2^{r+1}$, $\varpi ^{2^{r+1}-1}$ if $s=2^{r+1}$, for $%
Sp(2^{r}(2b+1))$.%
\end{tabular}%
\right. $
\end{enumerate}

By \textsl{the prime factorization of an integer} $n\geq 2$ we mean the
unique expression $n=p_{1}^{r_{1}}\cdots p_{t}^{r_{t}}$ with $1<p_{1}<\cdots
<p_{t}$ the set of all prime factors of $n$. In term of this factorization
one defines the partition on the set $\left\{ 2,\cdots ,n\right\} $ by

\begin{quote}
$\left\{ 2,\cdots ,n\right\} =Q_{0}(n)\underset{1\leq i\leq t}{\amalg }%
Q_{p_{i}}(n)$ with $Q_{p_{i}}(n)=\{p_{i}^{r}\mid 1\leq r\leq r_{s}\}$.
\end{quote}

\noindent \textbf{Lemma 4.6.} \textsl{Assume that} $G=SU(n)$\textsl{\ with }$%
n=p_{1}^{r_{1}}\cdots p_{t}^{r_{t}}$\textsl{\ (resp. }$G=Sp(n)$\textsl{\
with }$n=2^{r}(2b+1)$\textsl{). For a }$s\in D(G)$ \textsl{let }$a_{s}\in 
\mathbb{Z}$\textsl{\ be the order of the class }$\overline{\theta }(\gamma
_{2s-1})\in E_{3}^{\ast ,0}(PG)$\textsl{.} \textsl{Then}

\begin{enumerate}
\item[i)] $a_{s}=p_{i}$ \textsl{or }$1$\textsl{\ for }$s\in Q_{{\small p}%
_{i}}(n)$\textsl{\ or} $s\in Q_{0}(n)$

\textsl{(resp. }$a_{s}=2$ \textsl{or }$1$\textsl{\ for }$s=2^{{\small r+1}}$%
\textsl{\ or} $s\neq 2^{{\small r+1}}$\textsl{);}

\item[ii)] \textsl{there exists a class }$\rho _{2s-1}\in E_{3}^{\ast
,1}(PG) $\textsl{\ satisfying }$C^{\ast }(\rho _{2s-1})=a_{s}\gamma _{2s-1}$%
\textsl{;}

\item[iii)] $\mathcal{F}(PG)=\Lambda (\rho _{2s-1}^{\prime })_{s\in D(G)}$%
\textsl{, where }$\rho _{2s-1}^{\prime }=\kappa (\rho _{2s-1})\in H^{\ast
}(PG)$\textsl{.}
\end{enumerate}

\noindent \textbf{Proof.} For $SU(n)$ (resp. for $Sp(n)$) the relation i)
comes from the formula (4.10) of $\overline{\theta }(\gamma _{2s-1})$, the
presentation (4.2) of the ring $E_{3}^{\ast ,0}(PSU(n))$ (resp. $E_{3}^{\ast
,0}(PSp(n))$), as well as \cite[Theorem 1.1]{DLin} (see also (5.5)).

Property ii) follows from $\overline{\theta }(a_{s}\gamma _{2s-1})=0$,
together with the exactness of the sequence (3.11).

To show iii) we set for a multi--index $I\subseteq D(G)$ that

\begin{quote}
$a_{I}=\underset{s\in I}{\Pi }a_{s}\in \mathbb{Z}$, $\gamma _{I}=\underset{%
s\in I}{\Pi }\gamma _{2s-1}\in H^{\ast }(G)$, $\rho _{I}=\underset{s\in I}{%
\Pi }\rho _{2s-1}^{\prime }\in H^{\ast }(PG)$.
\end{quote}

\noindent For $G=SU(n)$ or $Sp(n)$ the degree of the covering $%
c:G\rightarrow PG$ is $\deg c=n$ or $2$. On the other hand by (4.9) the
monomial $\gamma _{D(G)}$ is a generator of the top degree cohomology group $%
H^{m}(G)=\mathbb{Z}$, $m=\dim G$, while the relations $a_{D(G)}=\deg c$, $%
C^{\ast }(\rho _{D(G)})=a_{D(G)}\cdot \gamma _{D(G)}$ by i) and ii) implies
that the monomial $\rho _{D(G)}$ is a generator of the group $H^{m}(PG)$ $=%
\mathbb{Z}$. By \cite[Lemma 2.9]{DZ2} the set $\{1,\rho _{I}\mid $ $%
I\subseteq D(G)\}$ of monomials spans a direct summand of $\mathcal{F}(PG)$
with rank $2^{\left\vert D(G)\right\vert }$. It follows then from

\begin{quote}
$\dim (\mathcal{F}(PG)\otimes \mathbb{Q)=}\dim (\mathcal{F}(G)\otimes 
\mathbb{Q)=}2^{\left\vert D(G)\right\vert }$
\end{quote}

\noindent that the set $\{1,\rho _{I}\mid $ $I\subseteq D(G)\}$ of monomials
is a basis of $\mathcal{F}(PG)$.

The proof of iii) will be completed once we show that $\rho _{2s-1}^{\prime
2}=0$, $s\in D(G)$. Assume that $n=2^{r}(2b+1)$. Since the square of an odd
degree cohomology class belongs to $\func{Im}\beta _{2}$, and since the
relation $\rho _{2s-1}^{\prime }\in \func{Im}\kappa $ implies that $\rho
_{2s-1}^{2}\in \func{Im}\pi ^{\ast }$ by \cite[Lemma 2.8]{DZ2}, we have

\begin{quote}
$\rho _{2s-1}^{\prime 2}\in \func{Im}[\pi ^{\ast }\circ \widehat{\beta }%
_{2}:E_{3}^{\ast ,1}(PG;\mathbb{F}_{2})\rightarrow E_{3}^{\ast
,0}(PG)\rightarrow H^{\ast }(PG)]$.
\end{quote}

\noindent On the other hand Lemma 4.3 implies that, if $G=SU(n)$ (resp. $%
G=Sp(n)$),

\begin{quote}
$\func{Im}\pi ^{\ast }\circ \widehat{\beta }_{2}=\{0,2^{r-t-1}\omega
^{2^{t}}\mid 0\leq t\leq r-1\}$ (resp. $\{0,\omega ^{2^{r}}\}$).
\end{quote}

\noindent One obtains the desired relation $\rho _{2s-1}^{\prime 2}=0$ for
the degree reason.$\square $

\bigskip

To simplify notation we reserve $\rho _{2s-1}$ for $\rho _{2s-1}^{\prime }$.
Substituting the formula (4.9) for $H^{\ast }(G\times S^{1})$ into (1.8)
yields the following exact sequence, in which $H^{\ast }(PG)_{\left\langle
\omega \right\rangle }=\Lambda (\rho _{2s-1})_{s\in D(G)}$ by the relations
ii) and iii) of Lemma 4.6,

\begin{quote}
$0\rightarrow \Lambda (\rho _{2s-1})_{s\in D(G)}\overset{C^{\ast }}{%
\rightarrow }\Lambda (\gamma _{1},\gamma _{2s-1})_{s\in D(G)}\overset{\theta 
}{\rightarrow }H^{\ast }(PG)\overset{\omega }{\rightarrow }\left\langle
\omega \right\rangle \rightarrow 0$.
\end{quote}

\noindent Since the quotient group $H^{\ast }(PG)_{\left\langle \omega
\right\rangle }$ is free the map $g$ in the exact sequence (1.7) has a
splitting homomorphism. Therefore, the formula (3.14) becomes functional to
yield the presentation

\begin{enumerate}
\item[(4.11)] $H^{\ast }(PG)=\frac{J(\omega )\otimes \Lambda (\rho
_{2s-1})_{s\in D(G)}}{\left\langle \omega \cdot \func{Im}\theta
\right\rangle }$, where

\item[(4.12)] $J(\omega )=\left\{ 
\begin{tabular}{l}
$\frac{\mathbb{Z}[\omega ]}{\left\langle b_{n,r}\omega ^{r},\text{ }1\leq
r\leq n\right\rangle }$ for $G=SU(n)$; \\ 
$\mathbb{Z}\oplus \frac{\mathbb{Z}[\omega ]^{+}}{\left\langle 2\omega
,\omega ^{2^{r+1}}\right\rangle }$ for $G=Sp(n)$ with $2^{r}(2b+1)$%
\end{tabular}%
\right. $ by (4.2).
\end{enumerate}

\noindent \textbf{Theorem 4.7. }\textsl{The integral cohomologies of the
groups }$PSp(n)$\textsl{\ (}$n=2^{r}(2b+1)$\textsl{) and }$PSU(n)$\textsl{\ (%
}$n=p_{1}^{r_{1}}\cdots p_{t}^{r_{t}}$\textsl{) are}

\begin{enumerate}
\item[i)] $H^{\ast }(PSp(n))=\Lambda (\rho _{4s-1})_{s\in \{1,\cdots
,n\}}\oplus \sigma _{2}(PSp(n))$, \textsl{where}

$\sigma _{2}(PSp(n))=\frac{\mathbb{F}_{2}[\omega ]^{+}\otimes \Lambda (\rho
_{3},\rho _{7},\cdots ,\rho _{4n-1})}{\left\langle \omega ^{2^{r+1}},\text{ }%
\omega \cdot \rho _{2^{r+2}-1}\right\rangle }$\textsl{.}

\item[ii)] $H^{\ast }(PSU(n))=\Lambda (\rho _{2s-1})_{s\in \{2,\cdots ,n\}}%
\underset{1\leq s\leq t}{\oplus }\sigma _{p_{s}}(PSU(n))$\textsl{, where}

$\sigma _{p_{s}}(PSU(n))=\frac{\mathbb{Z}[\omega ]^{+}\otimes \Lambda (\rho
_{2s-1})_{s\in \{2,\cdots ,n\}}}{\left\langle \omega \theta (\gamma
_{I})\mid I\subseteq \{1\}\amalg Q_{p_{s}}(n)\right\rangle }$\textsl{.}
\end{enumerate}

\noindent \textbf{Proof.} For both $G=Sp(n)$ and $SU(n)$ the presentation of
the free part $\mathcal{F}(PG)$ stated in the theorem has been shown by iii)
of Lemma 4.6. It remains to establish the formulae for the ideals $\sigma
_{p}(PG)$.

If $G=Sp(n)$ the formulae (4.11) and (4.12) imply, in addition to $\sigma
_{p}(PSp(n))$ $=0$ when $p\neq 2$, that the map $h$ in Theorem 1.3 restricts
a surjection

\begin{quote}
$h_{2}:\mathbb{Z}[\omega ]^{+}\otimes \Lambda (\rho _{3},\rho _{7},\cdots
,\rho _{4n-1})\rightarrow \sigma _{2}(PSp(n))$
\end{quote}

\noindent with $\ker h_{2}=\left\langle \omega \cdot \func{Im}\theta (\gamma
_{I}),\omega \cdot \func{Im}\theta (\gamma _{1}\gamma _{I})\mid I\subseteq
D(Sp(n))\right\rangle $. By formula (4.10) and the property ii) of Theorem
1.3

\begin{quote}
$\theta (\gamma _{I})=0$ or $\omega ^{2^{r}}\rho _{I\backslash \{2^{r+1}\}}$
if $2^{r+1}\notin I$ or $2^{r+1}\in I$;

$\theta (\gamma _{1}\gamma _{I})=2\rho _{I}$ or $\rho _{I}$ if $%
2^{r+1}\notin I$ or $2^{r+1}\in I$.
\end{quote}

\noindent Summarizing $\ker h_{2}=\left\langle 2\omega ,\omega
^{2^{r+1}},\omega \cdot \rho _{2^{r+2}-1}\right\rangle $, completing the
proof of i).

For $G=SU(n)$\textsl{\ }one can show that\textsl{\ }the ideal $\ker
h=\left\langle \omega \cdot \func{Im}\theta \right\rangle $ admits the
simplification

\begin{enumerate}
\item[a)] $\ker h=\left\langle \omega \theta (\gamma _{I})\mid I\subseteq
\{1\}\amalg Q_{p_{s}}(n)\text{, }1\leq s\leq t\right\rangle $.
\end{enumerate}

\noindent To show this we note that each multi--index $K\subseteq \{2,\cdots
,n\}$ has the partition

\begin{quote}
$K=K_{0}\underset{1\leq s\leq t}{\sqcup }K_{s}$ with $K_{0}=K\cap Q_{0}(n)$, 
$K_{s}=K\cap Q_{p_{s}}(n)$.
\end{quote}

\noindent Let $b_{s}=p_{1}^{\left\vert K_{1}\right\vert }\cdots \widehat{%
p_{s}^{\left\vert K_{s}\right\vert }}\cdots p_{t}^{\left\vert
K_{t}\right\vert }$ with $\left\vert K_{s}\right\vert $ the cardinality of $%
K_{s}$. Since the set $\{b_{1},\cdots ,b_{t}\}$ of integers is co--prime
there is a set $\{q_{1},\cdots ,q_{t}\}$ of integers satisfying $\Sigma
q_{s}b_{s}=1$. One obtains a) from the following calculation based on ii) of
Theorem 1.3 and ii) of Lemma 4.6, where $H_{s}$ is the complement of $K_{s}$
in $K$:

\begin{quote}
$\theta (\gamma _{K})=\Sigma q_{s}\theta (b_{s}\gamma _{K})=\Sigma q_{s}\rho
_{H_{s}}\cup \theta (\gamma _{K_{s}})$;

$\theta (\gamma _{1}\gamma _{K})=\Sigma q_{s}\theta (b_{s}\gamma _{1}\gamma
_{K})=\Sigma q_{s}\rho _{H_{s}}\cup \theta (\gamma _{1}\gamma _{K_{s}})$.
\end{quote}

On the other hand, by (4.12) the ring $J(\omega )$ admits the decomposition

\begin{enumerate}
\item[b)] $J(\omega )=\mathbb{Z}\underset{1\leq s\leq t}{\oplus }$ $%
J_{s}(\omega )$ with $J_{s}(\omega )=\frac{\mathbb{Z}[\omega ]^{+}}{%
\left\langle p_{s}^{r_{s}}\omega ,p_{s}^{r_{s}-1}\omega
^{p_{s}},p_{s}^{r_{s}-2}\omega ^{p_{s}^{2}},\cdots ,\omega
^{p_{s}^{r_{s}}}\right\rangle }$.
\end{enumerate}

\noindent Substituting a) and b) into (4.11) and notice by ii) of Theorem
1.3 that

\begin{quote}
$\theta (\gamma _{I})\in \sigma _{p_{s}}(PSU(n))$ for all $I\subseteq
Q_{p_{s}}(n)$,
\end{quote}

\noindent one gets, in addtion to $\sigma _{p}(PSU(n))$ $=0$ for $p\notin
\{p_{1},\cdots ,p_{t}\}$, the isomorphisms

\begin{enumerate}
\item[c)] $\frac{J_{s}(\omega )\otimes \Lambda (\rho _{2s-1})_{s\in
\{1,\cdots ,n\}}}{\left\langle \omega \theta (\gamma _{I})\mid I\subseteq
\{1\}\amalg Q_{p_{s}}(n)\right\rangle }\cong \sigma _{p_{s}}(PSU(n))$, $%
1\leq s\leq t$.$\square $
\end{enumerate}

In view the formula c) for the ideal $\sigma _{p_{s}}(PSU(n))$ a complete
description of the ring $H^{\ast }(PSU(n))$ asks for a formula expressing
the classes $\theta (\gamma _{I})$ with $I\subseteq \{1\}\amalg Q_{p_{s}}(n)$
by the generators $\omega $ and $\rho _{2s-1}$ of the ring $H^{\ast
}(PSU(n)) $. We emphasize at this point that presently the classes $\rho
_{2s-1}$ are only determined modulo the ideal $\left\langle \omega
\right\rangle $, while formulae for $\theta (\gamma _{I})$ may vary with
respect to different choices of the classes $\rho _{2s-1}$. With an
appropriate choice of the primary $1$--forms $\rho _{2s-1}$ by their
characteristic polynomials (see formula (5.8)) the following result will be
established in Section 5.4.

\bigskip

\noindent \textbf{Proposition 4.8.} \textsl{If }$n=p^{r}$\textsl{\ the
classes }$\rho _{2p^{s}-1}$\textsl{, }$1\leq s\leq r$\textsl{, may be chosen
so that, for an }$I=\left\{ p^{i_{1}},\cdots ,p^{i_{k}}\right\} \in
\{1\}\sqcup Q_{p}(n)$\textsl{\ with\ }$0\leq i_{1}<\cdots <i_{k}\leq r$%
\textsl{,}

\begin{enumerate}
\item[i)] \textsl{the class }$\theta (\gamma _{I})\in H^{\ast }(PSU(n))$%
\textsl{\ is divisible by }$p$\textsl{\ whenever }$i_{k}<r$\textsl{;}

\item[ii)] $\theta (\gamma _{2p^{s}-1})=p^{r-s}$\textsl{\ if }$I=\left\{
p^{s}\right\} $ \textsl{is a singleton;}

\item[iii)] $\theta (\gamma _{I})=(\frac{1}{p}\theta (\gamma _{I^{e}}))\cdot
\rho _{2p^{i_{k}}-1}+(\frac{1}{p}\theta (\gamma _{I^{\partial }}))\cdot
\omega ^{p^{i_{k}}-p^{i_{k}-1}}$
\end{enumerate}

\noindent \textsl{where }$I^{e}=\left\{ p^{i_{1}},\cdots
,p^{i_{k-1}}\right\} $\textsl{, }$I^{\partial }=\left\{ p^{i_{1}},\cdots
,p^{i_{k}-1}\right\} $\textsl{.}$\square $

\bigskip

In view of the formula c) the assumption $n=p^{r}$ in Proposition 4.8 may
not be necessary for the following reason. If $n=p^{r}n^{\prime }$\textsl{\ }%
with $(p,n^{\prime })=1$ then

\begin{quote}
$Q_{p}(n)=Q_{p}(p^{r})$.
\end{quote}

\noindent It indicates that, with respect to the cohomology ring map induced
by the inclusion $PSU(p^{r})\rightarrow PSU(n)$, the set $\left\langle
\omega \theta (\gamma _{I})\mid I\subseteq \{1\}\amalg Q_{p}(n)\right\rangle 
$ of relations on $\sigma _{p}(PSU(n))$ are in one to one correspondence
with that on $\sigma _{p}(PSU(p^{r}))$, where the latter are handled by
Proposition 4.8.

\bigskip

\noindent \textbf{Example 4.9.} In Proposition 4.8 we notice that

\begin{quote}
a) if $i_{k}=i_{k-1}+1$ then $\gamma _{I^{\partial }}=0$ by the relation $%
\gamma _{2r-1}^{2}=0$ on $H^{\ast }(U(n))$;

b) the classes $\theta (\gamma _{I^{e}})$ and $\theta (\gamma _{I^{\partial
}})$ are always divisible by\textsl{\ }$p$ by ii).
\end{quote}

\noindent Therefore, formula iii) gives an effective recurrence to evaluate $%
\theta (\gamma _{I})$. As examples when $n=2^{3}$ we get that

\begin{quote}
$\theta (\gamma _{\{1,2,4\}})=2\rho _{3}\rho _{7}$;

$\theta (\gamma _{\{1,2,8\}})=2\rho _{3}\rho _{15}+\omega ^{4}\rho _{3}\rho
_{7}$;

$\theta (\gamma _{\{1,4,8\}})=2\rho _{7}\rho _{15}+\omega ^{2}\rho _{3}\rho
_{15}$;

$\theta (\gamma _{\{2,4,8\}})=\omega \rho _{7}\rho _{15}$;

$\theta (\gamma _{\{1,2,4,8\}})=\rho _{3}\rho _{7}\rho _{15}$.$\square $
\end{quote}

\subsection{The integral cohomology of the groups $PE_{6}$ and $PE_{7}$}

For a space $X$ and a prime $p$ the pair $\{H^{\ast }(X;\mathbb{F}%
_{p});\delta _{p}\}$ with $\delta _{p}=r_{p}\circ \beta _{p}$ is a cochain
complex whose cohomology $\overline{H}^{\ast }(X;\mathbb{F}_{p})$ is the $%
\func{mod}p$ \textsl{Bockstein cohomology }of $X$. For $(G,p)=(E_{6},3)$ and 
$(E_{7},2)$ the complexes $\{H^{\ast }(PG;\mathbb{F}_{p});\delta _{p}\}$
have been decided by Lemma 4.4 and Theorem 4.5. Explicitly we have

\begin{enumerate}
\item[(4.13a)] $H^{\ast }(PE_{6};\mathbb{F}_{3})=[\func{Im}\pi ^{\ast
}\otimes \Lambda _{\mathbb{F}_{3}}(\varsigma _{1},\varsigma _{7})]\otimes
\Lambda _{\mathbb{F}_{3}}(\varsigma _{2s-1})_{s\in \{2,5,6,8\}}$ with

\quad i) $\func{Im}\pi ^{\ast }=\frac{\mathbb{F}_{3}[\omega ,x_{4}]}{%
\left\langle \omega ^{9},x_{4}^{3}\right\rangle }$; ii) $\delta
_{3}(\varsigma _{1})=\omega $, $\delta _{3}(\varsigma _{7})=x_{4}$, $\delta
_{2}(\varsigma _{2s-1})=0$,

where in view of iv) of Theorem 4.5 and in the order of $r=1,2,4,5,6,8$

$\quad \varsigma _{2r-1}:=\iota ,\zeta _{3},\zeta _{7},\zeta _{9},\zeta
_{11},\zeta _{15}-x_{4}\zeta _{7}$.

\item[(4.13b)] $H^{\ast }(PE_{7};\mathbb{F}_{2})=[\func{Im}\pi ^{\ast
}\otimes \Delta _{\mathbb{F}_{2}}(\varsigma _{1},\varsigma _{2t-1})_{t\in
\{3,5,9\}}]\otimes \Lambda _{\mathbb{F}_{2}}(\varsigma _{2s-1})_{s\in
\{8,12,14\}}$ with

\quad i) $\func{Im}\pi ^{\ast }=\frac{\mathbb{F}_{2}[\omega
,x_{3},x_{5},x_{9}]}{\left\langle
x_{1}^{2},x_{3}^{2},x_{5}^{2},x_{9}^{2}\right\rangle }$; ii) $\delta
_{3}(\varsigma _{1})=\omega $, $\delta _{3}(\varsigma _{2t-1})=x_{t}$, $%
\delta _{3}(\varsigma _{2s-1}))=0$,

where in view of v) of Theorem 4.5 and in the order of $r=1,3,5,8,9,12,14$

$\quad \varsigma _{2r-1}:=\iota ,\zeta _{5},\zeta _{9},\zeta
_{15}+x_{3}\zeta _{9},\zeta _{17},\zeta _{23}+x_{3}\zeta _{17},\zeta
_{27}+x_{5}\zeta _{17}$.
\end{enumerate}

In what follows we put $c_{\{1,4\}}=\delta _{3}(\varsigma _{1}\varsigma
_{7})\in H^{9}(PE_{6};\mathbb{F}_{3})$. For a multi--index $I\subseteq
\{1,3,5,9\}$ define the elements

\begin{quote}
$c_{I}=\delta _{2}(\varsigma _{I})\in H^{\ast }(PE_{7};\mathbb{F}_{2})$ with 
$\varsigma _{I}=\underset{s\in I}{\Pi }\varsigma _{2s-1}$,.
\end{quote}

\noindent \textbf{Lemma 4.10.} \textsl{The cohomology} $\overline{H}^{\ast
}(PG;\mathbb{F}_{p})$, \textsl{together with} $\func{Im}\delta _{p}$\textsl{,%
} \textsl{are given by}

\begin{quote}
\textsl{i)} $\overline{H}^{\ast }(PE_{6};\mathbb{F}_{3})\cong \Lambda _{%
\mathbb{F}_{3}}(\omega ^{8}\varsigma _{1},x_{4}^{2}\varsigma _{7},\varsigma
_{3},\varsigma _{9},\varsigma _{11},\varsigma _{15})$

$\quad \func{Im}\delta _{3}=\frac{\mathbb{F}_{3}[\omega
,x_{4},c_{\{1,4\}}]^{+}}{\left\langle \omega
^{9},x_{4}^{3},c_{\{1,4\}}^{2},x_{1}^{8}x_{4}^{2}c_{\{1,4\}}\right\rangle }%
\otimes \Lambda _{\mathbb{F}_{3}}(\varsigma _{3},\varsigma _{9},\varsigma
_{11},\varsigma _{15})$

\textsl{ii)} $\overline{H}^{\ast }(PE_{7};\mathbb{F}_{2})\cong \Lambda _{%
\mathbb{F}_{2}}(\omega \varsigma _{1},x_{3}\varsigma _{5},x_{5}\varsigma
_{9},x_{9}\varsigma _{17},\varsigma _{15},\varsigma _{23},\varsigma _{27})$

\quad $\func{Im}\delta _{2}=\frac{\mathbb{F}_{2}[\omega
,x_{3},x_{5},x_{9},c_{I}]^{+}}{\left\langle \omega
^{2},x_{3}^{2},x_{5}^{2},x_{9}^{2},D_{I},R_{I},S_{I,J}\right\rangle }\otimes
\Lambda _{\mathbb{F}_{2}}(\varsigma _{15},\varsigma _{23},\varsigma _{27})$%
\textsl{\ with }$\left\vert I\right\vert ,\left\vert J\right\vert \geq 2$,
\end{quote}

\noindent \textsl{where in the presentation of }$\func{Im}\delta _{2}$%
\textsl{\ the relations }$D_{I},R_{I},S_{I,J}$ \textsl{are, respectively, }

\begin{enumerate}
\item[(4.14)] $\underset{t\in I}{\Sigma }x_{t}c_{I_{t}}=0$\textsl{,} $\quad (%
\underset{t\in I}{\Pi }x_{t})c_{I}=0$\textsl{, }$c_{I}c_{J}+\underset{t\in I}%
{\Sigma }x_{t}\underset{s\in I_{t}\cap J}{\Pi }\varsigma
_{2s-1}^{2}c_{\left\langle I_{t},J\right\rangle }=0$\textsl{,}
\end{enumerate}

\noindent \textsl{and where} $(\varsigma _{1}^{2},\varsigma
_{5}^{2},\varsigma _{9}^{2},\varsigma _{17}^{2})=(\omega ,x_{5},x_{9},0)$%
\textsl{, }$I_{t}$ \textsl{is the complement of} $t\in I$\textsl{,} $%
\left\langle I,J\right\rangle $ \textsl{denotes the complement of} \textsl{%
the intersection }$I\cap J$\textsl{\ in the union} $I\cup J$\textsl{.}

\bigskip

\noindent \textbf{Proof.} The results in i) and ii) come from the same
calculation. We may therefore focus on the relatively nontrivial case ii).
In the presentation (4.13b) of $H^{\ast }(PE_{7};\mathbb{F}_{2})$ the first
factor $\func{Im}\pi ^{\ast }\otimes \Delta _{\mathbb{F}_{2}}(\varsigma
_{1},\varsigma _{2t-1})$ is the Koszul complex studied in \cite[Theorem 2.1]%
{DZ2}, while the differential $\delta _{2}$ acts trivially on the second
factor $\Lambda _{\mathbb{F}_{2}}(\varsigma _{2s-1})$. One gets by \cite[%
Theorem 2.1]{DZ2} and the K\"{u}nneth formula that

\begin{enumerate}
\item[c)] $\overline{H}^{\ast }(PE_{7};\mathbb{F}_{2})=\Delta _{\mathbb{F}%
_{2}}(x_{1}\varsigma _{1},x_{3}\varsigma _{5},x_{5}\varsigma
_{9},x_{9}\varsigma _{17})\otimes \Lambda _{\mathbb{F}_{2}}(\varsigma
_{15},\varsigma _{23},\varsigma _{27})$,

\item[d)] $\func{Im}\delta _{2}=\frac{\func{Im}\pi ^{\ast }\{1,c_{I}\}^{+}}{%
\left\langle D_{J},R_{K}\right\rangle }\otimes \Lambda _{\mathbb{F}%
_{2}}(\varsigma _{15},\varsigma _{23},\varsigma _{27})$,
\end{enumerate}

\noindent where $I,J,K\subseteq \{1,3,5,9\}\QTR{sl}{\ }$with $\left\vert
I\right\vert ,\left\vert J\right\vert ,\left\vert K\right\vert \geq 2$. By
v) of Theorem 4.5 the formula c) of $\overline{H}^{\ast }(PE_{7};\mathbb{F}%
_{2})$ is identical to the one stated in ii). To modify the additive
presentation of $\func{Im}\delta _{2}$ in d) into its ring presentation in
ii) one needs to clarify the multiplicative rule among the classes $c_{I}$%
's. This brings us the relations of the type $S_{I,J}$ which are obtained by
the following calculation:

\begin{quote}
$c_{I}c_{J}=\delta _{2}(\varsigma _{I})\delta _{2}(\varsigma _{J})=\delta
_{2}(\delta _{2}(\varsigma _{I})\varsigma _{J})$ (since $\delta _{2}^{2}=0$)

$=\delta _{2}(\underset{t\in I}{\Sigma }x_{t}\varsigma _{I_{t}}\varsigma
_{J})$ (since $\delta _{2}(\varsigma _{I})=\underset{t\in I}{\Sigma }%
-x_{t}\varsigma _{I_{t}}$)

$=\delta _{2}(\underset{t\in I}{\Sigma }x_{t}\underset{s\in I_{t}\cap J}{\Pi 
}\zeta _{2s-1}^{2}\zeta _{\left\langle I_{t},J\right\rangle })$ (with $%
\underset{s\in I\cap J}{\Pi }\varsigma _{2s-1}^{2}=1$ if $I\cap J=\emptyset $%
)

$=\underset{t\in I}{\Sigma }x_{t}\underset{s\in I_{t}\cap J}{\Pi }\varsigma
_{2s-1}^{2}c_{\left\langle I_{t},J\right\rangle }$ (since $\delta
_{2}(x_{t})=0$, $\delta _{2}(\varsigma _{2s-1}^{2})=0$).$\square $
\end{quote}

The presentations of the Bockstein cohomology\textsl{\ }$\overline{H}^{\ast
}(PG;\mathbb{F}_{p})$ in Lemma 4.10 provide us with crucial information on
reduction $r_{p}:H^{\ast }(PG)\rightarrow H^{\ast }(PG;\mathbb{F}_{p})$ for $%
(G,p)=(E_{6},3)$ and $(E_{7},2)$. Indeed, in view of the decomposition (4.8)
for $X=PG$ we write $r_{p}^{0}$ and $r_{p}^{1}$ for the restrictions of $%
r_{p}$ on $\mathcal{F}(PG)$ and $\sigma _{p}(PG)$, respectively. With

\begin{quote}
$\dim \overline{H}^{\ast }(PE_{6};\mathbb{F}_{3})=2^{6}$ and $\dim \overline{%
H}^{\ast }(PE_{7};\mathbb{F}_{2})=2^{7}$
\end{quote}

\noindent by Lemma 4.10 the result \cite[Theorem 2.7]{DZ2} implies that

\begin{enumerate}
\item[(4.15)] $r_{p}^{1}$ is an isomorphism $\sigma _{p}(PG)\cong \func{Im}%
\delta _{p}$;

\item[(4.16)] $r_{p}^{0}$ induces an isomorphism $\mathcal{F}(PG)/p\cdot 
\mathcal{F}(PG)\rightarrow \overline{H}^{\ast }(PG;\mathbb{F}_{p})$.
\end{enumerate}

To apply the exact sequence (1.8) to compute $H^{\ast }(PG)$ we need
information on the cohomology of the group $G$. In term of the set $S(G)$ of
primary characteristic polynomials for $G=E_{6}$ and $E_{7}$ over $\mathbb{Z}
$ given in Table 5.4 (see also Example 3.6) let $D(G)$ be the degree set of
elements in $S(G)$.\textbf{\ }That is

\begin{quote}
$D(E_{6})=\{2,5,6,8,9,12\}$; $D(E_{7})=\{2,6,8,10,12,14,18\}$
\end{quote}

\noindent For each $s\in D(G)$ we set $\gamma _{2s-1}=\varphi (P)\in
E_{3}^{\ast ,1}(G)$, where $P\in S(G)$ with $\deg P=s$. To save notations we
maintain $\gamma _{2s-1}$ for $\kappa (\gamma _{2s-1})\in H^{\ast }(G)$. By 
\cite[Theorem 1.9]{DZ2} the integral cohomology $H^{\ast }(G)$ has the
following presentations

\begin{enumerate}
\item[(4.17a)] $H^{\ast }(E_{6})=\Delta (\gamma _{3})\otimes \Lambda (\gamma
_{9},\gamma _{11},\gamma _{15},\gamma _{17},\gamma _{23})\oplus \sigma
_{2}(E_{6})\oplus \sigma _{3}(E_{6})$ with

$\qquad \sigma _{2}(E_{6})=\mathbb{F}_{2}[x_{3}]^{+}/\left\langle
x_{3}^{2}\right\rangle \otimes \Delta (\gamma _{3})\otimes \Lambda (\gamma
_{9},\gamma _{15},\gamma _{17},\gamma _{23})$\textsl{,}

where $\gamma _{3}^{2}=x_{3}$\textsl{, }$x_{3}\gamma _{11}=0$\textsl{, }$%
x_{4}\gamma _{23}=0$\textsl{.}

\item[(4.17b)] $H^{\ast }(E_{7})=\Delta (\gamma _{3})\otimes \Lambda _{%
\mathbb{Z}}(\gamma _{11},\gamma _{15},\gamma _{19},\gamma _{23},\gamma
_{27},\gamma _{35})\underset{p=2,3}{\oplus }\sigma _{p}(E_{7})$ with

$\qquad \sigma _{3}(E_{7})=\frac{\mathbb{F}_{3}[x_{4}]^{+}}{\left\langle
x_{4}^{3}\right\rangle }\otimes \Lambda (\gamma _{3},\gamma _{11},\gamma
_{15},\gamma _{19},\gamma _{27},\gamma _{35})$\textsl{,}

where $\gamma _{3}^{2}=x_{3}$\textsl{, }$x_{4}\gamma _{23}=0$, $r\in
\{11,19,35\}$,
\end{enumerate}

\noindent where the classes $x_{i}$'s are the special Schubert classes on $%
E_{n}/T$, $n=6,7$, specified in (5.1), and where the formulae of the ideals $%
\sigma _{3}(E_{6})$ and $\sigma _{2}(E_{7})$ are not needed in sequel, hence
are omitted.

With $\gamma _{2s-1}=\kappa \varphi (P)$ one computes $\theta (\gamma
_{2s-1})$ by the formula (3.17) to yield

\begin{enumerate}
\item[(4.18)] $\theta (\gamma _{2s-1})=0$ with the only exceptions: $\theta
(\gamma _{2s-1})=\omega ^{8}$ or $\omega $ for\textsl{\ }$(G,s)=(E_{6},9)$
or $(E_{7},2)$ (see the contents of Table 5.5).
\end{enumerate}

\noindent In view of the presentations (4.13a) and (4.13b) introduce the
elements

\begin{quote}
$\mathcal{C}_{\{1,4\}}=\beta _{3}(\varsigma _{1}\varsigma _{7})$ $\in \sigma
_{3}(PE_{6})$, $\mathcal{C}_{K}=\beta _{2}(\underset{t\in K}{\Pi }\varsigma
_{2t-1})\in \sigma _{2}(E_{7})$,
\end{quote}

\noindent where $K\subseteq \{1,3,5,9\}$. For $I,J\subset \{1,3,5,9\}$\ with 
$\left\vert I\right\vert ,\left\vert J\right\vert \geq 2$ let $\mathcal{D}%
_{I},\mathcal{R}_{I},\mathcal{S}_{I,J}\in \sigma _{2}(E_{7})$ be
respectively the elements\ obtained by substituting in the polynomials $%
D_{I},R_{I},S_{I,J}\in \mathbb{F}_{2}[\omega ,x_{3},x_{5},x_{9},c_{I}]^{+}$
in (4.14) the classes $c_{I}$\ by $C_{I}$.

\bigskip

\noindent \textbf{Lemma 4.11.} \textsl{For each }$s\in D(G)$\textsl{\ there
exists a class }$\rho _{2s-1}\in \func{Im}\kappa $\textsl{\ satisfying}

\begin{enumerate}
\item[i)] $C^{\ast }(\rho _{2s-1})=a_{s}\cdot \gamma _{2s-1}$\textsl{, where 
}$a_{s}=1$\textsl{\ with the only exceptions: }

$a_{9}=3$\textsl{\ for }$E_{6}$\textsl{; }$a_{2}=2$\textsl{\ for }$E_{7}$%
\textsl{.}

\item[ii)] $\mathcal{F}(PG)=\left\{ 
\begin{tabular}{l}
$\Delta (\rho _{3})\otimes \Lambda (\rho _{9},\rho _{11},\rho _{15},\rho
_{17},\rho _{23})$ \textsl{with} $\rho _{3}^{2}=x_{3}$ \textsl{for} $G=E_{6}$%
; \\ 
$\Lambda (\rho _{3},\rho _{11},\rho _{15},\rho _{19},\rho _{23},\rho
_{27},\rho _{35})$ \textsl{for} $G=E_{7}$.%
\end{tabular}%
\right. $

\item[iii)] $\sigma _{p}(PG)=\left\{ 
\begin{tabular}{l}
$\frac{\mathbb{F}_{3}[\omega ,x_{4},\mathcal{C}_{\{1,4\}}]^{+}\otimes
\Lambda (\rho _{3},\rho _{9},\rho _{11},\rho _{15},\rho _{17})}{\left\langle
\omega ^{9},x_{4}^{3},\omega \cdot \rho _{17},\mathcal{C}_{\{1,4\}}^{2},%
\omega ^{8}x_{4}^{2}\mathcal{C}_{\{1,4\}}\right\rangle }$ \textsl{for} $%
(G,p)=(E_{6},3);$ \\ 
$\frac{\mathbb{F}_{2}[\omega ,x_{3},x_{5},x_{9},\mathcal{C}_{I}]^{+}}{%
\left\langle x_{1}^{2},x_{3}^{2},x_{5}^{2},x_{9}^{2},\mathcal{D}_{I},%
\mathcal{R}_{I},\mathcal{S}_{I,J}\right\rangle }\otimes \Lambda (\rho
_{15},\rho _{23},\rho _{27})$ \textsl{for} $(G,p)=(E_{7},2)$.%
\end{tabular}%
\right. $
\end{enumerate}

\noindent \textbf{Proof.} The classes $\rho _{2s-1}\in \func{Im}\kappa $ are
specified by $\theta (a_{s}\gamma _{2s-1})=0$ by (4.18), as well as the
exact sequence (1.8). This shows i).

For ii) the degree of the covering $c:G\rightarrow PG$ satisfies that $\deg
c=a_{D(D)}$ ($=3$ or $2$ for $E_{6}$ or $E_{7}$) by i). The same argument as
that used in the proof of iii) of Lemma 4.6 shows that the set $\{1,\rho
_{I}\mid $ $I\subseteq D(G)\}$ of monomials is a basis of $\mathcal{F}(PG)$.
Since $\rho _{2s-1}\in \func{Im}\kappa $ and since the reduction $r_{p}$
preserves $\func{Im}\kappa $ we get from (4.16) and the formula of $%
\overline{H}^{\ast }(PG;\mathbb{F}_{p})$ in Lemma 4.10 the relations

\begin{enumerate}
\item[(4.19a)] $r_{3}(\rho _{2s-1})\in \left\{ \omega ^{8}\varsigma
_{1},x_{4}^{2}\varsigma _{7},\varsigma _{3},\varsigma _{9},\varsigma
_{11},\varsigma _{15}\right\} $ \textsl{for} $PE_{6}$;

\item[(4.19b)] $r_{2}(\rho _{2s-1})\in \left\{ \omega \varsigma
_{1},x_{3}\varsigma _{5},x_{5}\varsigma _{9},x_{9}\varsigma _{17},\varsigma
_{15},\varsigma _{23},\varsigma _{27}\right\} $ \textsl{for} $PE_{7}$.
\end{enumerate}

\noindent These determine the classes $r_{p}(\rho _{2s-1})\in $ $H^{\ast
}(PG;\mathbb{F}_{p})$ for the degree reason.

To complete the proof of ii) it remains to show that $\rho _{2s-1}^{2}=0$
with the only exception $\rho _{3}^{2}=x_{3}$ when $G=E_{6}$. For $PE_{7}$
this comes from $\rho _{2s-1}^{2}\in \sigma _{2}(PE_{7})$, the injectivity
of $r_{2}$ on $\sigma _{2}(PE_{7})$ by (4.15), as well as the relation
(4.19b). For $PE_{6}$ we shall show in the proof of Theorem 4.12 that the
map $C^{\ast }$ restricts to a ring monomorphism $\sigma
_{2}(PE_{6})\rightarrow \sigma _{2}(E_{6}\times S^{1})$. One obtains i) from

\begin{quote}
$C^{\ast }(\rho _{2s-1})\equiv \gamma _{2s-1}\func{mod}2$, $C^{\ast
}(x_{3})\equiv x_{3}\func{mod}2$ (see ii) of Lemma 4.2),
\end{quote}

\noindent as well as the relations $\gamma _{3}^{2}=x_{6},\gamma
_{2s-1}^{2}=0$, $s\geq 3$, on the ring $H^{\ast }(E_{6})$ by (4.19a).

Finally, granted with the isomorphisms $r_{p}^{1}$ in (4.15), the relations
(4.19a) and (4.19b) suffice to translate the presentations of $\func{Im}%
\delta _{p}$ in Lemma 4.10 into the formulae for the ideals $\sigma _{p}(PG)$
stated in iii).$\square $

\bigskip

We present the integral cohomology $H^{\ast }(PG)$ with $G=E_{6}$\ and $%
E_{7} $ by the set $\{\rho _{2s-1}\}\subset \func{Im}\kappa $ of $1$--forms
specified by i) of Lemma 4.11, the subring $\func{Im}\pi ^{\ast }$ given by
(4.2), together the torsion classes $\mathcal{C}_{K}$.

\bigskip

\noindent \textbf{Theorem 4.12. }\textsl{The rings }$H^{\ast }(PG)$\textsl{\
with }$G=E_{6}$\textsl{\ and }$E_{7}$ \textsl{are}

\begin{enumerate}
\item[i)] $H^{\ast }(PE_{6})=\Delta (\rho _{3})\otimes \Lambda (\rho
_{9},\rho _{11},\rho _{15},\rho _{17},\rho _{23})\underset{p=2,3}{\oplus }%
\sigma _{p}(PE_{6})$\textsl{\ with}

$\quad \sigma _{2}(PE_{6})=\mathbb{F}_{2}[x_{3}]^{+}/\left\langle
x_{3}^{2}\right\rangle \otimes \Delta (\rho _{3})\otimes \Lambda (\rho
_{9},\rho _{15},\rho _{17},\rho _{23})$\textsl{,}

$\quad \sigma _{3}(PE_{6})=\frac{\mathbb{F}_{3}[\omega ,x_{4},\mathcal{C}%
_{\{1,4\}}]^{+}\otimes \Lambda (\rho _{3},\rho _{9},\rho _{11},\rho
_{15},\rho _{17})}{\left\langle \omega ^{9},x_{4}^{3},\omega \cdot \rho
_{17},\mathcal{C}_{\{1,4\}}^{2},\omega ^{8}x_{4}^{2}\mathcal{C}%
_{\{1,4\}}\right\rangle }$\textsl{,}

\textsl{that are subject the relations}

$\quad \rho _{3}^{2}=x_{3}$\textsl{, }$x_{3}\rho _{11}=0$\textsl{, }$%
x_{4}\rho _{23}=0$\textsl{, }$\omega \rho _{23}=x_{4}^{2}\mathcal{C}%
_{\{1,4\}}$\textsl{,} $\mathcal{C}_{\{1,4\}}\rho _{23}=0$.

\item[ii)] $H^{\ast }(PE_{7})=\Lambda (\rho _{3},\rho _{11},\rho _{15},\rho
_{19},\rho _{23},\rho _{27},\rho _{35})\underset{p=2,3}{\oplus }\sigma
_{p}(PE_{7})$ \textsl{with}

$\quad \sigma _{2}(PE_{7})=\frac{\mathbb{F}_{2}[\omega ,x_{3},x_{5},x_{9},%
\mathcal{C}_{I}]^{+}}{\left\langle x_{1}^{2},x_{3}^{2},x_{5}^{2},x_{9}^{2},%
\mathcal{D}_{I},\mathcal{R}_{I},\mathcal{S}_{I,J}\right\rangle }\otimes
\Lambda (\rho _{15},\rho _{23},\rho _{27})$;

$\quad \sigma _{3}(PE_{7})=\frac{\mathbb{F}_{3}[x_{4}]^{+}}{\left\langle
x_{4}^{3}\right\rangle }\otimes \Lambda (\rho _{3},\rho _{11},\rho
_{15},\rho _{19},\rho _{27},\rho _{35})$\textsl{,}

\textsl{that are subject to the relations, where} $K\subseteq
\{1,3,5,9\},s\in \{2,6,10,18\}$\textsl{,}

$\quad x_{4}\varrho _{23}=0$\textsl{,} $\rho _{2s-1}C_{K}=0$ \textsl{if} $%
s\in K,$ $\rho _{2s-1}C_{K}=x_{\frac{s}{2}}C_{K\cup \{s\}}$ \textsl{if} $%
s\notin K$\textsl{.}
\end{enumerate}

\noindent \textbf{Proof. }The presentations of $\mathcal{F}(PG)$, together
with the ideals $\sigma _{3}(PE_{6})$ and $\sigma _{2}(PE_{7})$, have been
shown by Lemma 4.11. Moreover, with $\left\langle \omega \right\rangle \in
\sigma _{3}(PE_{6})$ by $3\omega =0$ ($\left\langle \omega \right\rangle \in
\sigma _{2}(PE_{7})$ by $2\omega =0$) formula (4.8) yields that

\begin{quote}
$H^{\ast }(PE_{6})_{\left\langle \omega \right\rangle }=\mathcal{F}(PE_{6})%
\underset{p\neq 3}{\oplus }\sigma _{p}(PE_{6})\oplus \sigma
_{3}(PE_{6})/\left\langle \omega \right\rangle $

(resp. $H^{\ast }(PE_{7})_{\left\langle \omega \right\rangle }=\mathcal{F}%
(PE_{7})\underset{p\neq 2}{\oplus }\sigma _{p}(PE_{7})\oplus \sigma
_{2}(PE_{6})/\left\langle \omega \right\rangle $).
\end{quote}

\noindent In particular, the map $C^{\ast }$ in (1.8) restricts to the
monomorphisms

\begin{quote}
$\sigma _{p}(PE_{6})\rightarrow \sigma _{p}(E_{6}\times S^{1})$, $p\neq 3$ ($%
\sigma _{p}(PE_{7})\rightarrow \sigma _{p}(E_{7}\times S^{1})$, $p\neq 2$).
\end{quote}

\noindent It implies, in addition to $\sigma _{p}(PG)=0$ for all $p\neq 2,3$%
, that

\begin{quote}
$\sigma _{2}(PE_{6})\cong \sigma _{2}(E_{6})$ (resp. $\sigma
_{3}(PE_{7})\cong \sigma _{3}(E_{7})$) under $C^{\ast }$.
\end{quote}

\noindent One obtains from (4.17a) (resp. (4.17b)) the presentation of $%
\sigma _{2}(PE_{6})$ (resp. $\sigma _{3}(PE_{7})$), together with the
relation $x_{3}\rho _{11}=0$ (resp. $x_{4}\varrho _{23}=0$), as that stated
in the theorem.

Finally, granted with the isomorphism $r_{p}^{1}$ in (4.15), Theorem 4.5, as
well as the relations (4.19a) and (4.19b) that characterize the reduction $%
r_{p}^{1}:$ $\mathcal{F}(PG)\rightarrow H^{\ast }(PG;\mathbb{F}_{p})$, one
verifies the following relations

\begin{quote}
a) $x_{4}\rho _{23}=0$\textsl{, }$\omega \rho _{23}=x_{4}^{2}\mathcal{C}%
_{\{1,4\}}$\textsl{,} $\mathcal{C}_{\{1,4\}}\rho _{23}=0$;

b) $x_{4}\varrho _{23}=0$\textsl{,} $\rho _{2s-1}C_{K}=0$ if $s\in K,$ $=x_{%
\frac{s}{2}}C_{I\cup \{s\}}$ if $s\notin K$.
\end{quote}

\noindent that characterize the action of $\mathcal{F}(PG)$ on $\sigma
_{p}(PG)$, respectively for $(G,p)=(E_{6},3)$ or $(E_{7},2)$. As example,
combining (4.15) and iv) of Theorem 4.5, the relations in a) are verified by
the following calculations in the ring $H^{\ast }(PE_{6};\mathbb{F}_{3})$

\begin{quote}
$r_{3}(\omega \rho _{23})=\omega x_{4}^{2}\zeta _{7}=x_{4}^{2}(\omega \zeta
_{7}-\iota x_{4})=x_{4}^{2}c_{\{1,4\}}$;

$r_{3}(x_{4}\rho _{23})=x_{4}^{3}\zeta _{7}=0$; $r_{3}(\mathcal{C}%
_{\{1,4\}}\rho _{23})=(\omega \zeta _{7}-\iota x_{4})x_{4}^{2}\zeta _{7}=0$,
\end{quote}

\noindent where $r_{3}(\rho _{23})=x_{4}^{2}\zeta _{7}$ by (4.19a).$\square $

\subsection{Historical remarks}

The rings $H^{\ast }(PG;\mathbb{F}_{p})$ concerned by Theorem 4.5 have been
previously computed by Baum and Browder \cite{BB} for $G=SU(n),Sp(n)$, and
by Toda, Kono and Ishitoya \cite{IKT,Ko} for $G=E_{6},E_{7}$. However, with
respect to our explicitly constructed generators on $H^{\ast }(PG;\mathbb{F}%
_{p})$ the Bockstein operator $\beta _{p}$ can be effectively calculated,
which played a crucial role to decide the torsion ideals $\sigma _{p}(PG)$
of the integral cohomology $H^{\ast }(PG)$.

In \cite[Theorem A]{R} Ruiz stated a presentation of the integral cohomology
ring of the projective complex Stiefel manifold $Y_{n,n-m}$. It implies when
taking $m=0$ that the map $h$ in Theorem 1.3 is an isomorphism for the group 
$PSU(n)$. That is, the relations of the form $\omega \theta (\gamma _{I})$
with $\left\vert I\right\vert \geq 2$ are absent. Computation in Example 4.9
indicates that these missing relations are highly nontrival.

The formula (2.9) for the Borel transgression $\tau $ is an essential
ingredient for the computation throughout Section 4. On the other hand in 
\cite[formula (4)]{K} Ka\v{c} stated a formula for the differential $d_{2}$
on $E_{2}^{\ast ,\ast }(G;\mathbb{F}_{p})$ which implies that the map $\tau $
is an isomorphism for any semi--simple Lie group $G$ and characteristic $p$.
These suggest to us that a proof of Theorem 2.4 is unavoidable.

According Grothendieck \cite{G} the subring $\func{Im}\pi ^{\ast }$ $\subset
H^{\ast }(PG)$ is the \textsl{Chow ring} $A^{\ast }(PG^{c})$ of the
reductive algebraic group $PG^{c}$ corresponding to $PG$. In this regard the
formulae (4.2) presents the Chow rings $A^{\ast }(PG^{c})$ for $G=SU(n)$, $%
Sp(n)$, $E_{6}$, $E_{7}$ by explicit Schubert classes on $G/T$.

It has been known for a long time that for cohomology with coefficients in a
field $\mathbb{F}$ one has $E_{3}^{\ast ,\ast }(G;\mathbb{F})=H^{\ast }(G;%
\mathbb{F})$ \cite{K,Re}. For the integral cohomology the maps $\pi ^{\ast }$
and $\kappa $ in (3.5) and (3.6) provide a direct passage from $E_{3}^{\ast
,\ast }(G)$ to $H^{\ast }(G)$ useful to show the much stronger relation

\begin{enumerate}
\item[(4.20)] $E_{3}^{\ast ,\ast }(G)=H^{\ast }(G)$.
\end{enumerate}

\noindent Indeed, for the $1$--connected Lie groups $G$ the relation (4.20)
was conjecture by Marlin \cite{M} and has been confirmed by the authors in 
\cite[Theorem 3.6]{DZ2}, while our proofs of Theorem 4.7 and Theorem 4.12
indicate that the relation (4.20) holds more generally by all compact Lie
groups.

\section{Schubert calculus}

This section generates and records the intermediate data facilitating the
computation in Section 4. It serves also the purpose to illustrate how the
construction and computation with the cohomology of compact Lie groups $G$
can be boiling down to computing with certain polynomials in the Schubert
classes on $G/T$.

\subsection{ Schubert presentation of the ring $H^{\ast }(G/T)$}

Let $m=n-1,n,6$ or $7$ in accordance to $G=SU(n),Sp(n),E_{6}$ or $E_{7}$. In
view of the Schubert basis $\left\{ \omega _{1},\cdots ,\omega _{m}\right\} $
on $H^{2}(G/T)$ (see Lemma 2.3) we define the polynomials $c_{r}(G)\in
H^{\ast }(BT)=\mathbb{Z}\left[ \omega _{1},\cdots ,\omega _{m}\right] $ to
be the $r^{th}$ elementary symmetric polynomials on the set $\Omega (G)$
specified below

\begin{quote}
$\Omega (SU(n))=\left\{ \omega _{{\footnotesize 1}},\omega _{{\footnotesize k%
}}-\omega _{{\footnotesize k-1}},-\omega _{{\footnotesize n-1}}\mid 2\leq
k\leq n-1\right\} $;

$\Omega (Sp(n))=\left\{ \pm \omega _{{\footnotesize 1}},\pm (\omega _{%
{\footnotesize k}}-\omega _{{\footnotesize k-1}})\mid 2\leq k\leq n\right\} $%
;

$\Omega (E_{6})=\left\{ \omega _{{\footnotesize 6}},\omega _{{\footnotesize 5%
}}-\omega _{{\footnotesize 6}},\omega _{{\footnotesize 4}}-\omega _{%
{\footnotesize 5}},\omega _{{\footnotesize 2}}+\omega _{{\footnotesize 3}%
}-\omega _{{\footnotesize 4}},\omega _{{\footnotesize 1}}+\omega _{%
{\footnotesize 2}}-\omega _{{\footnotesize 3}},\omega _{{\footnotesize 2}%
}-\omega _{{\footnotesize 1}}\right\} $;

${\small \Omega (E}_{7}{\small )=}\left\{ \omega _{{\footnotesize 7}},\omega
_{{\footnotesize 6}}-\omega _{{\footnotesize 7}},\omega _{{\footnotesize 5}%
}-\omega _{{\footnotesize 6}},\omega _{{\footnotesize 4}}-\omega _{%
{\footnotesize 5}},\omega _{{\footnotesize 2}}+\omega _{{\footnotesize 3}%
}-\omega _{{\footnotesize 4}},\omega _{1}+\omega _{2}-\omega _{3},\omega _{%
{\footnotesize 2}}-\omega _{{\footnotesize 1}}\right\} $.
\end{quote}

\noindent Using the \textsl{Weyl coordinates} $\sigma _{\lbrack i_{1},\cdots
,i_{k}]}$ for the Schubert classes on $G/T$ (\cite[Definition 2]{DZ2}) we
introduce also the \textsl{special Schubert classes} $x_{r}$ on $G/T$ by

\begin{enumerate}
\item[(5.1)] $(x_{3},x_{4}):=(\sigma _{\lbrack 5,4,2]},\sigma _{\lbrack
6,5,4,2]})$ for $G=E_{6}$;

$(x_{3},x_{4},x_{5},x_{9}):=(\sigma _{\lbrack 5,4,2]},\sigma _{\lbrack
6,5,4,2]},\sigma _{\lbrack {7,6,5,4,2}]},\sigma _{\lbrack
1,5,4,3,7,6,5,4,2]})$ for $G=E_{7}$.
\end{enumerate}

\noindent With these notation we have by \cite{DZ1} the following
presentations of the cohomologies $H^{\ast }(G/T)$ in term of explicit
generators and relations, where in the cases of $G=E_{6}$ and $E_{7}$ and in
comparison with Theorem 2.6, we use $x_{\deg y_{i}}$ in place of $y_{i}$,
and write $R_{\deg h_{i}}$, $R_{\deg f_{j}}$, $R_{\deg g_{j}}$ instead of $%
h_{i},f_{j},g_{j}$. In view of the map $\psi :G/T\rightarrow BT$ in (3.6) we
shall also reserve the nation $c_{r}(G)$ for the class $\psi ^{\ast
}(c_{r}(G))\in H^{\ast }(G/T)$.

\bigskip

\noindent \textbf{Theorem 5.1. }\textsl{The ring }$H^{\ast }(G/T)$ \textsl{%
has the following presentations}

\begin{enumerate}
\item[i)] $H^{\ast }(SU(n)/T){\small =}\mathbb{Z}\left[ \omega _{1},\cdots
,\omega _{n-1}\right] /\left\langle c_{2},\cdots ,c_{n}\right\rangle $%
\textsl{,} $c_{r}=c_{r}({\small SU(n)})$,

\item[ii)] $H^{\ast }(Sp(n)/T)=\mathbb{Z}\left[ \omega _{1},\cdots ,\omega
_{n}\right] /\left\langle c_{2},c_{4},\cdots ,c_{2n}\right\rangle $\textsl{,}
$c_{2r}=c_{2r}(Sp(n))$.

\item[iii)] $H^{\ast }(E_{6}/T)=\mathbb{Z}[\omega _{1},\cdots ,\omega
_{6},x_{3},x_{4}]/\left\langle
R_{2},R_{3},R_{4},R_{5},R_{6},R_{8},R_{9},R_{12}\right\rangle $\textsl{,
where}

$R_{2}=4\omega _{2}^{2}-c_{2}$;

$R_{3}=2x_{3}+2\omega _{2}^{3}-c_{3}$;

$R_{4}=3x_{4}+\omega _{2}^{4}-c_{4}$;

$R_{5}=2\omega _{2}^{2}x_{3}-\omega _{2}c_{4}+c_{5}$;

$R_{6}=x_{3}^{2}-\omega _{2}c_{5}+2c_{6}$;

$R_{8}=x_{4}(c_{4}-\omega _{2}^{4})-2c_{5}x_{3}-\omega _{2}^{2}c_{6}+\omega
_{2}^{3}c_{5};$

$R_{9}=2x_{3}c_{6}-\omega _{2}^{3}c_{6}$;

$R_{12}=x_{4}^{3}-c_{6}^{2}$.

\item[iv)] $H^{\ast }(E_{7}/T)=\mathbb{Z}[\omega _{1},\cdots ,\omega
_{7},x_{3},x_{4},x_{5},x_{9}]/\left\langle R_{t}\right\rangle $\textsl{,
where }$t\in \{2,3,4,5,6,8,$ $9,10,12,14,18\}$\textsl{, and where}

$R_{2}=4\omega _{2}^{2}-c_{2}$;

$R_{3}=2x_{3}+2\omega _{2}^{3}-c_{3}$;

$R_{4}=3x_{4}+\omega _{2}^{4}-c_{4}$;

$R_{5}=2x_{5}-2\omega _{2}^{2}x_{3}+\omega _{2}c_{4}-c_{5}$;

$R_{6}=x_{3}^{2}-\omega _{2}c_{5}+2c_{6}$;

$R_{8}=3x_{4}^{2}-x_{5}(2\omega _{2}^{3}-c_{3})-2x_{3}c_{5}+2\omega
_{2}c_{7}-\omega _{2}^{2}c_{6}+\omega _{2}^{3}c_{5}$;

$R_{9}=2{x_{9}}+{x_{4}(2\omega _{2}^{2}x_{3}-\omega _{2}c_{4}+c_{5})}-2{%
x_{3}c_{6}}-{\omega _{2}^{2}c_{7}}+{\omega _{2}^{3}c_{6}}$;

$R_{10}=x_{5}^{2}-2x_{3}c_{7}+\omega _{2}^{3}c_{7}$;

$R_{12}=x_{4}^{3}-4x_{5}c_{7}-c_{6}^{2}+(2\omega
_{2}^{3}-c_{3})(x_{9}+x_{4}x_{5})+2\omega _{2}x_{5}c_{6}+3\omega
_{2}x_{4}c_{7}+c_{5}c_{7}$;

$R_{14}=c_{7}^{2}-(2\omega _{2}^{2}x_{3}-\omega
_{2}c_{4}+c_{5})x_{9}+2x_{3}x_{4}c_{7}-\omega _{2}^{3}x_{4}c_{7}$;

$R_{18}=x_{9}^{2}+2x_{5}c_{6}c_{7}-x_{4}c_{7}^{2}-(2\omega
_{2}^{2}x_{3}-\omega _{2}c_{4}+c_{5})x_{4}x_{9}-(2\omega
_{2}^{3}-c_{3})x_{5}^{3}$

$\qquad \quad -5\omega _{2}x_{5}^{2}c_{7}$.$\square $
\end{enumerate}

As in Section 4.1 we take a set $\Omega =\{\phi _{1},\cdots ,\phi _{m}\}$ of
fundamental dominant weights as a basis for the unit lattice $\Lambda _{e}$
of the adjoint Lie group $PG$, and let $\{t_{1},\cdots ,t_{m}\}$ be the
corresponding basis on the group $H^{1}(T)$. Granted with the Cartan
matrices of simple Lie groups given in \cite[p.59]{H} the formula (2.9) for
the transgression $\tau $ in the fibration $\pi :$ $PG\rightarrow G/T$
yields that

\bigskip

\noindent \textbf{Lemma 5.2.} \textsl{In the order of }$%
G=SU(n),Sp(n),E_{6},E_{7}$\textsl{\ the following relations hold in the
quotient group} $H^{2}(G/T)/\func{Im}\tau $

\begin{quote}
i) $\omega _{k}=k\omega _{1}$\textsl{, }$1\leq k\leq n-1$\textsl{, }$n\omega
_{1}=0$\textsl{;}

ii) $\omega _{k}=k\omega _{1}$\textsl{,}$1\leq k\leq n$\textsl{, }$2\omega
_{1}=0$\textsl{;}

iii) $\omega _{2}=\omega _{4}=0$\textsl{, }$\omega _{1}=\omega _{5}=2\omega
_{3}=2\omega _{6}$\textsl{, }$3\omega _{1}=0$\textsl{;}

iv) $\omega _{1}=\omega _{3}=\omega _{4}=\omega _{6}=0$\textsl{, }$\omega
_{5}=\omega _{7}=\omega _{2}$\textsl{, }$2\omega _{2}=0$\textsl{.}
\end{quote}

\noindent \textsl{Consequently}

\begin{quote}
a) $c_{r}(SU(n))\mid _{\tau (t_{1})=\cdots =\tau (t_{n-1})=0\QTR{sl}{\ }}=%
\binom{n}{r}\omega _{1}^{r}$\textsl{, }$2\leq r\leq n$\textsl{;}

b) $c_{2r}(Sp(n))\mid _{\tau (t_{1})=\cdots =\tau (t_{n})=0\QTR{sl}{\ }}=%
\binom{n}{r}\omega _{1}^{2r}$\textsl{, }$1\leq r\leq n$\textsl{;}

c) $c_{r}(E_{6})\mid _{\tau (t_{1})=\cdots =\tau (t_{6})=0\QTR{sl}{\ }%
}=(-1)^{r}\binom{6}{r}\omega _{1}^{r}$\textsl{, }$1\leq r\leq 6$\textsl{;}

d) $c_{r}(E_{7})\mid _{\tau (t_{1})=\cdots =\tau (t_{7})=0\QTR{sl}{\ }}=%
\binom{7}{r}\omega _{2}^{r}$\textsl{, }$1\leq r\leq 7$.$\square $
\end{quote}

\subsection{Computing with the $\mathbb{F}_{p}$--characteristic polynomials}

Let $(G,p)=(SU(n),p)$, $(Sp(n),2)$, $(E_{6},3)$ and $(E_{7},2)$. From
Theorem 5.1 one deduces a set $S_{p}(G):=\left\{ \delta _{1},\cdots ,\delta
_{m}\right\} $ of\textsl{\ }primary characteristic polynomials for $G$ over $%
\mathbb{F}_{p}$ (see Example 3.7) as that presented in the following table,
where $(H)_{p}\in \mathbb{F}_{p}[\omega _{1},\cdots ,\omega _{m}]$ denotes $H%
\func{mod}p$, $H\in \mathbb{Z}[\omega _{1},\cdots ,\omega _{m}]$.

\begin{center}
\textbf{Table 5.1. }A set $S_{p}(G)$ of primary characteristic polynomials
for $G$ over $\mathbb{F}_{p}$

\begin{tabular}{l|l}
\hline\hline
$(G,p)$ & $S_{p}(G)=\left\{ \delta _{1},\cdots ,\delta _{m}\right\} $ \\ 
\hline
$(SU(n),p)$ & $\left( c_{k}\right) _{p},2\leq k\leq n$ \\ \hline
$(Sp(n),2)$ & $\left( c_{2k}\right) _{2},1\leq k\leq n$ \\ \hline
$(E_{6},3)$ & 
\begin{tabular}{l}
$(\omega
_{2}^{2}-c_{2})_{3},(c_{2}^{2}-c_{4})_{3},(c_{5}+c_{2}c_{3})_{3},\left(
c_{6}-c_{2}c_{4}-c_{3}^{2}\right) _{3}$ \\ 
$(-c_{3}c_{5}-c_{2}c_{6})_{3},(c_{6}c_{3})_{3}$%
\end{tabular}
\\ \hline
$(E_{7},2)$ & 
\begin{tabular}{l}
$\left( c_{2}\right) _{2},\left( c_{3}\right) _{2},\left( c_{5}+\omega
_{2}c_{4}\right) _{2},\left( c_{4}^{2}+\omega _{2}^{2}c_{6}+\omega
_{2}^{3}c_{5}+\omega _{2}^{8}\right) _{2},$ \\ 
$\left( \omega _{2}^{2}c_{7}+\omega _{2}^{3}c_{6}\right) _{2},\left(
c_{6}^{2}+c_{4}^{3}\right) _{2},\left( c_{7}^{2}+c_{4}^{2}c_{6}+\omega
_{2}^{2}c_{6}^{2}\right) _{2}$%
\end{tabular}
\\ \hline\hline
\end{tabular}
\end{center}

Let $\tau ^{\prime }$ be the transgression in the fibration $\pi ^{\prime }$
contained in the diagram (3.8). Granted with the class $\varpi =\tau
^{\prime }(t_{0})$ determined in Lemma 4.2, as well as the results of Lemma
5.2, formula (3.17) is applicable to evaluate the derivative $\partial
P/\partial \varpi $ for $P\in S_{p}(G)$, as that presented in the following
table:

\begin{center}
\textbf{Table 5.2.} The derivative $\partial P/\partial \varpi $ for $P\in
S_{p}(G)$

\begin{tabular}{l|l}
\hline\hline
$(G,p)$ & $\left\{ \partial P/\partial \varpi \mid P\in S_{p}(G)\right\} $
\\ \hline
$(SU(n),p)$ & $\left( \binom{n}{k}\omega _{1}^{k-1}\right) _{p}$, $2\leq
k\leq n$ \\ \hline
$(Sp(n),2)$ & $\left( \binom{n}{k}\omega _{1}^{2k-1}\right) _{2}$, $1\leq
k\leq n$ \\ \hline
$(E_{6},3)$ & $0,0,0,0,0,\left( \omega _{1}^{8}\right) _{3}$ \\ \hline
$(E_{7},2)$ & $\left( \omega _{2}\right) _{2},\left( \omega _{2}^{2}\right)
_{2},0,0,{0,}0,\left( \omega _{2}^{13}\right) _{2}$ \\ \hline\hline
\end{tabular}%
.
\end{center}

\noindent Consequently, one obtains $\theta _{1}(\varphi _{p}(P))\in
E_{3}^{\ast ,0}(PG;\mathbb{F}_{p})$, $P\in S_{p}(G)$, as that stated in the
formula (4.5).

By the algorithm given in the proof of Lemma 3.11, for each $P\in S_{p}(G)$
with $\theta _{1}(\varphi _{p}(P))=0$ one can construct a polynomial $%
P^{\prime }\in \left\langle \func{Im}\widetilde{\tau }_{p}\right\rangle \cap
\ker f_{p}$ satisfying the relation $C^{\ast }(\varphi _{p}(P^{\prime
}))=\varphi _{p}(P)$. Explicitly, a set $S_{p}(PG)$ of the polynomials $%
P^{\prime }$ so obtained is given in the following table:

\begin{center}
\textbf{Table 5.3.} A set $S_{p}(PG)$ of primary characteristic polynomials
over $\mathbb{F}_{p}$

\begin{tabular}{l|l}
\hline\hline
$(G,p)$ & $S_{p}(PG)\subset \left\langle \func{Im}\widetilde{\tau }%
_{p}\right\rangle \cap \ker f_{p}$ \\ \hline
$(SU(n),p)$ & 
\begin{tabular}{l}
$\left( c_{k}\right) _{p}$ for $2\leq k<p^{r}$, $(c_{k}-t_{n,k}c_{p^{r}}%
\omega _{1}^{k-p^{r}})_{p}$ for $k\geq p^{r}$ \\ 
where $n=p^{r}n^{\prime }$ with $(n\prime ,p)=1$;%
\end{tabular}
\\ \hline
$(Sp(n),2)$ & 
\begin{tabular}{l}
$\left( c_{2k}\right) _{2}$ for $2\leq k<2^{r}$, $(c_{2k}-t_{n,2}c_{2^{r+1}}%
\omega _{1}^{2(k-2^{r})})_{2}$ for $k\geq 2^{r}$ \\ 
where $n=2^{r}(2b+1)$;%
\end{tabular}
\\ \hline
$(E_{6},3)$ & $%
\begin{tabular}{l}
$\left( \omega _{2}^{2}-c_{2}\right) _{3},\left( c_{2}^{2}-c_{4}\right)
_{3},\left( c_{5}+c_{2}c_{3}\right) _{3},$ \\ 
$\left( c_{6}-c_{2}c_{4}-c_{3}^{2}\right) _{3},\left(
-c_{3}c_{5}-c_{2}c_{6}\right) _{3}$;%
\end{tabular}%
$ \\ \hline
$(E_{7},2)$ & 
\begin{tabular}{l}
$\left( c_{3}-c_{2}\omega _{2}\right) _{2},\left( c_{5}+\omega
_{2}c_{4}\right) _{2},\left( c_{4}^{2}+\omega _{2}^{2}c_{6}+\omega
_{2}^{3}c_{5}+\omega _{2}^{8}\right) _{2}$ \\ 
$\left( \omega _{2}^{2}c_{7}+\omega _{2}^{3}c_{6}\right) _{2},\left(
c_{6}^{2}+c_{4}^{3}\right) _{2},\left( c_{7}^{2}+c_{4}^{2}c_{6}+\omega
_{2}^{2}c_{6}^{2}-c_{2}\omega _{2}^{12}\right) _{2}$%
\end{tabular}
\\ \hline\hline
\end{tabular}
\end{center}

\noindent where $t_{n,k}>0$ is an integer with $t_{n,k}\binom{n}{p^{r}}%
\equiv \binom{n}{k}\func{mod}p$.

By formula (4.7) the cohomology $H^{\ast }(PG;\mathbb{F}_{p})$ has the
presentation

\begin{enumerate}
\item[(5.2)] $H^{\ast }(PG;\mathbb{F}_{p})=\pi ^{\ast }E_{3}^{\ast ,0}(PG;%
\mathbb{F}_{p})\otimes \Delta (\iota ,\zeta _{2\deg P-1})_{P\in S_{p}(PG)}$,
\end{enumerate}

\noindent where $\zeta _{2\deg P-1}=\kappa \circ \varphi _{p}(P)$, $P\in
S_{p}(PG)$. In the statement and the proof of the following result, we note
by the proof of Lemma 3.7 that, with $\zeta _{2\deg P-1}\in \func{Im}\kappa $%
,

\begin{enumerate}
\item[(5.3)] $\beta _{p}(\zeta _{2\deg P-1})\in \pi ^{\ast }(E_{3}^{\ast
,0}(PG))\subset H^{\ast }(PG)$,
\end{enumerate}

\noindent where the ring $E_{3}^{\ast ,0}(PG)$ has been determined by (4.2).

\bigskip

\noindent \textbf{Lemma 5.3. }\textsl{With respect to (5.2)}\textbf{\ }%
\textsl{the Bockstein} $\beta _{p}:H^{\ast }(PG;\mathbb{F}_{p})\rightarrow
H^{\ast }(PG)$ \textsl{satisfies that}

\begin{quote}
\textsl{i) for }$(G,p)=(SU(n),p)$\textsl{\ with }$n=p^{r}n^{\prime }$, $%
(n\prime ,p)=1$

$\qquad \beta _{p}(\zeta _{2s-1})=-p^{r-t-1}\omega ^{p^{t}}$ \textsl{if} $%
s=p^{t}$ \textsl{with} $t<r$\textsl{,} $0$ \textsl{otherwise;}

\textsl{ii) for }$(G,p)=(Sp(n),2)$\textsl{\ with }$n=2^{r}(2b+1)$\textsl{:}

$\qquad \beta _{2}(\zeta _{4s-1})=\omega ^{2^{r}}$ \textsl{if} $s=2^{r-1}$%
\textsl{,} $0$ \textsl{if} $s\neq 2^{r-1}$\textsl{;}

\textsl{iii) for }$(G,p)=(E_{6},3)$\textsl{\ and in the order of }$%
s=2,4,5,6,8$

$\qquad \beta _{3}(\zeta _{2s-1})=0,-x_{4},0,0,-x_{4}^{2}$\textsl{;}

\textsl{iv) for }$(G,p)=(E_{7},2)$\textsl{\ and in the order of }$%
s=3,5,8,9,12,14$

$\qquad \beta _{2}(\zeta
_{2s-1})=x_{3},x_{5},x_{3}x_{5},x_{9},x_{3}x_{9},x_{5}x_{9}.$
\end{quote}

\noindent \textbf{Proof. }For the case $(SU(n),p)$ with $n=p^{r}n^{\prime }$
and $(n\prime ,p)=1$ (resp. $(Sp(n),2)$ with $n=2^{r}(2b+1)$), one has by
(5.3) that

\begin{quote}
$\beta _{p}(\zeta _{2s-1})\in \mathbb{Z}[\omega ]^{+}/\left\langle
p^{r}\omega ,p^{r-1}\omega ^{p},\cdots ,\omega ^{p^{r}}\right\rangle $

(resp. $\beta _{2}(\zeta _{4s-1})\in \mathbb{Z}[\omega ]/\left\langle
2\omega ,\omega ^{2^{r+1}}\right\rangle $),
\end{quote}

\noindent where $\omega =\pi ^{\ast }(\omega _{1})$ (see iii) of Theorem
3.9). By the degree reason one gets

\begin{quote}
$\beta _{p}(\zeta _{2s-1})=0$ if $s\geq p^{r}$ (resp. $\beta _{2}(\zeta
_{4s-1})=0$ if $s\geq 2^{r}$).
\end{quote}

\noindent In the remaining cases $2\leq s<p^{r}$ (resp. $1\leq s<2^{r}$) an
integral lift of the characteristic polynomial $\left( c_{s}\right) _{p}$
(resp. $\left( c_{2s}\right) _{2}$) of the class $\zeta _{2s-1}$ (resp. $%
\zeta _{4s-1}$) is easily seen to be

\begin{quote}
$c_{s}-$ $\binom{n}{s}\omega _{1}^{s}\in \left\langle \func{Im}\widetilde{%
\tau }\right\rangle $ (resp. $c_{2s}-$ $\binom{n}{k}\omega _{1}^{2s}\in
\left\langle \func{Im}\widetilde{\tau }\right\rangle $).
\end{quote}

\noindent The formula in Lemma 3.7 then yields that

\begin{quote}
$\beta _{p}(\zeta _{2s-1})=\pi ^{\ast }\frac{1}{p}f(c_{s}-\binom{n}{s}\omega
_{1}^{s})=\pi ^{\ast }(-\frac{1}{p}\binom{n}{s}\omega _{1}^{s})=-\frac{1}{p}%
\binom{n}{s}\omega ^{s}$

(resp. $\beta _{2}(\zeta _{4s-1})=\pi ^{\ast }(\frac{1}{2}f(c_{2s}-\binom{n}{%
s}\omega _{1}^{2s})=\pi ^{\ast }(\frac{1}{2}\binom{n}{s}\omega _{1}^{2s})=%
\frac{1}{2}\binom{n}{s}\omega ^{2s}$),
\end{quote}

\noindent where the second equality follows from $f(c_{s})=0$ (resp. $%
f(c_{2s})=0$) by i) (resp. ii)) of Theorem 5.1. The relation i) (resp. ii))
of the present lemma comes now from $p^{r-t}\omega ^{p^{t}}=0$, $0\leq t\leq
r$ (resp. $2\omega ,\omega ^{2^{r+1}}=0$), and the property (5.6) of the
binomial coefficients $\binom{n}{s}$.

Turning to the cases $(G,p)=(E_{6},3)$ or $(E_{7},2)$ for each $P=(H)_{p}\in
S_{p}(PG)$ given by Table 5.3 the enclosed polynomial $H$ satisfies $H\in
\left\langle \func{Im}\widetilde{\tau }\right\rangle $ by the relation c) or
d) of Lemma 5.2. Therefore, by Lemma 3.7

\begin{quote}
$\beta _{p}(\zeta _{2\deg P-1})=\pi ^{\ast }(\frac{1}{p}f(H)\mid _{\tau
(t_{1})=\cdots =\tau (t_{n})=0})$.
\end{quote}

\noindent This formula, together with the relations $R_{i}$'s in the
presentations iii) and iv) of Theorem 5.1, suffices to yield the results in
iii) and iv). As examples, when $(G,p)=(E_{6},3)$ it gives rise to

\begin{quote}
$\beta _{3}(\zeta _{3})=\pi ^{\ast }\frac{\omega _{2}^{2}-c_{2}}{3}=\pi
^{\ast }\omega _{2}^{2}=0$;

$\beta _{3}(\zeta _{7})=\pi ^{\ast }\frac{c_{2}^{2}-c_{4}}{3}=\pi ^{\ast
}(5\omega _{2}^{2}-x_{4})=-x_{4}$;

$\beta _{3}(\zeta _{9})=\pi ^{\ast }\frac{c_{5}+c_{2}c_{3}}{3}=\pi ^{\ast
}(\omega _{2}^{5}+\omega _{2}x_{4}+\omega _{2}^{2}c_{3})=0$;

$\beta _{3}(\zeta _{11})=\pi ^{\ast }\frac{c_{2}c_{4}+c_{3}^{2}-c_{6}}{3}%
=\pi ^{\ast }(5\omega _{2}^{2}x_{4}+\omega _{2}^{6}+\omega _{2}c_{5}+\omega
_{2}^{3}c_{3}-3c_{6})=0$;

$\beta _{3}(\zeta _{15})=\pi ^{\ast }(-\frac{c_{3}c_{5}+c_{2}c_{6}}{3})=\pi
^{\ast }(-4x_{4}^{2}+c_{3}c_{5}+4\omega _{2}^{3}c_{5})=-x_{4}^{2}$,
\end{quote}

\noindent where, in each of the above equations, the second equality is
deduced from the relations $R_{i}$'s in iii) of Theorem 5.1, where the last
one is obtained by the formula c) of Lemma 5.2.$\square $

\bigskip

For $G=SU(n),Sp(n),E_{6}$ or $E_{7}$ consider the fibration $G/T\overset{%
\psi }{\hookrightarrow }BT\overset{Bi}{\rightarrow }BG$ induced by the
inclusion $i:T\rightarrow G$ of a maximal torus $T$, see (3.6). It was shown
in \cite[Lemma 5.4]{DZ3} that there exists a complex bundle $\xi $ on $BT$
so that

\begin{quote}
$C(\xi )=1+c_{1}+c_{2}+\cdots $ for $G=SU(n),E_{6}$ or $E_{7}$;

$P(\xi )=1+c_{2}+c_{4}+\cdots $ for $G=Sp(n)$,
\end{quote}

\noindent where $C(\xi )$ and $P(\xi )$ are the total Chern class and the
total Pontrjagin class of $\xi $, respectively. It follows from the
Wu--formula \cite[p.94]{MS} that the following relation holds on the ring $%
H^{\ast }(BT;\mathbb{F}_{2})$

\begin{enumerate}
\item[(5.4)] $Sq^{2s-2}c_{s}=c_{s-1}c_{s}+c_{s-2}c_{s+1}+\cdots +c_{2s-1}$
for $G=SU(n),E_{6}$ or $E_{7}$;

$Sq^{4s-2}c_{2s}=c_{2s-2}c_{2s}+c_{2s-4}c_{2s+2}+\cdots +c_{4s-2}$ for $%
G=Sp(n)$.
\end{enumerate}

\noindent \textbf{Lemma 5.4.} \textsl{With respect to the presentation (5.2)}
\textsl{the Steenrod operators }$\delta _{2}$ \textsl{and} $Sq^{2r}$ \textsl{%
on} $H^{\ast }(PG;\mathbb{F}_{p})$ \textsl{satisfy the following relations}

\begin{quote}
\textsl{i) }$G=SU(n)$\textsl{\ with }$n=2^{r}(2b+1):$

$\quad \delta _{2}\zeta _{2s-1}=\omega _{1}^{2^{r-1}}$ \textsl{if} $%
s=2^{r-1} $, $0$ \textsl{if} $s\neq 2^{r-1}$\textsl{; }

$\quad Sq^{2s-2}\zeta _{2s-1}=\zeta _{4s-3}$ \textsl{for }$2s-1\leq 2^{r-1}$%
\textsl{;}

\textsl{ii) }$G=Sp(n)$\textsl{\ with }$n=2^{r}(2b+1):$

$\quad \delta _{2}\zeta _{4s-1}=\omega _{1}^{2^{r}}$ \textsl{if} $s=2^{r-1}$%
, $0$ \textsl{if} $s\neq 2^{r-1}$\textsl{, }

$\quad Sq^{4s-2}\zeta _{4s-1}=\zeta _{8s-3}$\textsl{\ for }$4s-1\leq 2^{r}$;

\textsl{iii) }$G=E_{7}$\textsl{\ and in accordance to} $s=3,5,8,9,12,14%
\QTR{sl}{:}$

$\quad \delta _{2}\zeta
_{2s-1}=y_{3},y_{5},y_{3}y_{5},y_{9},y_{3}y_{9},y_{5}y_{9}$\textsl{; }$%
Sq^{2s-2}\zeta _{2s-1}=\zeta _{9},\zeta _{17},0,0,0,0.$
\end{quote}

\noindent \textbf{Proof.} The results on $\delta _{2}\zeta _{2s-1}$ comes
directly from the relation $\delta _{2}=r_{2}\circ \beta _{2}$, in which the
action of $\beta _{2}$ on $\zeta _{2s-1}$ has been decided in Lemma 5.3,
while the reduction $r_{2}$ is transparent in view of the presentations of
the rings $E_{3}^{\ast ,0}(PG)$ and $E_{3}^{\ast ,0}(PG;\mathbb{F}_{2})$ in
(4.2) and Lemma 4.3.

For $G=SU(n)$\textsl{\ }(resp. $Sp(n)$) with\textsl{\ }$n=2^{r}(2b+1)$ we
get from Lemma 3.8 and the formula (5.4) that

\begin{quote}
$Sq^{2s-2}\zeta _{2s-1}=\kappa \varphi _{2}(Sq^{2s-2}c_{s})=\kappa \varphi
_{2}(c_{2s-1})=\zeta _{4s-3}$

(resp. $Sq^{4s-2}\zeta _{4s-1}=\varphi _{2}(Sq^{4s-2}c_{2s})=\varphi
_{2}(c_{4s-2})=\zeta _{4s-3}$),
\end{quote}

\noindent where $2s-1\leq 2^{r-1}$ (resp. $4s-1\leq 2^{r}$), the second
equality comes from $\varphi _{2}(c_{s-1-t}c_{s+t})=0$ for all $0\leq t\leq
s-2$ by Lemma 3.4, and where the last equality follows from the fact that,
with $2s-1<2^{r-1}$, $c_{2s-1}$ is a characteristic polynomial of the class $%
\zeta _{4s-3}\in H^{\ast }(PSU(n);\mathbb{F}_{2})$, see Table (5.3). This
completes the proof of i) (resp. ii)).

Finally, for the case $G=E_{7}$, granted with the set $S_{2}(PE_{7})$ of
explicitly polynomials given by Table 5.3 and the Wu--formula (5.4), it is
straightforward to show the formulae in iii) by the method entailed above.
Alternatively, see \cite[Section 5.2]{DZ3} for the relevant computation.$%
\square $

\subsection{The integral characteristic polynomials}

In view of the presentation of the rings $H^{\ast }(G/T)$ in Theorem 5.1,
one formulates the set $S(G)$ of\textsl{\ }characteristic\textsl{\ }%
polynomials for the $1$--connected groups $G=SU(n),Sp(n),E_{6},E_{7}$ (see
Example 3.6) as that tabulated below

\begin{center}
\textbf{Table 5.4.} The set $S(G)$ of primary characteristic polynomials
over $\mathbb{Z}$

\begin{tabular}{l|l}
\hline\hline
$G$ & $S(G)$ \\ \hline
$SU(n)$ & $c_{k},2\leq k\leq n$ \\ \hline
$Sp(n)$ & $c_{2k},1\leq k\leq n$ \\ \hline
$E_{6}$ & $R_{2},R_{5},2R_{6}-x_{3}R_{3},R_{8},R_{9},3R_{12}-x_{4}^{2}R_{4}$
\\ \hline
$E_{7}$ & $%
R_{2},2R_{6}-x_{3}R_{3},R_{8},2R_{10}-x_{5}R_{5},3R_{12}-x_{4}^{2}R_{4},R_{14},2R_{18}-x_{9}R_{9} 
$ \\ \hline\hline
\end{tabular}
\end{center}

\noindent In (4.9), (4.17a) and (4.17b), the rings $H^{\ast }(G)$ are
presented by the primary $1$--forms $\gamma _{2\deg P-1}:=\kappa \varphi
(P)\in E_{3}^{\ast ,1}(G)$ with $P\in S(G)$.

Let $\tau ^{\prime }$ be the transgression in the fibration $\pi ^{\prime }$
in the diagram (3.8). Granted with the class $\varpi =\tau ^{\prime }(t_{0})$
determined by Lemma 4.2, as well as the results of Lemma 5.2, formula (3.17)
is ready to apply to evaluate the derivation $\partial P/\partial \varpi $
for $P\in S(G)$. The results are so obtained presented in the following
table.

\begin{center}
\textbf{Table 5.5.} The derivation $\partial P/\partial \varpi $ for $P\in
S(G)$

\begin{tabular}{l|l}
\hline\hline
$G$ & $\left\{ \partial P/\partial \varpi \mid P\in S(G)\right\} $ \\ \hline
$SU(n)$ & $\binom{n}{s}\omega _{1}^{s-1}$, $2\leq s\leq n$ \\ \hline
$Sp(n)$ & $\binom{n}{s}\omega _{1}^{2s-1},$ $1\leq s\leq n$ \\ \hline
$E_{6}$ & $0,0,0,0,\omega _{1}^{8},0$ \\ \hline
$E_{7}$ & $\omega _{2},\omega _{2}^{2}y_{3},\omega _{2}^{2}y_{5},0,\omega
_{2}^{2}(y_{9}+y_{4}y_{5}+y_{4}\omega _{2}^{5}),\omega _{2}^{9}(y_{4}+\omega
_{2}^{4}),0$ \\ \hline\hline
\end{tabular}
\end{center}

\subsection{Arithmetic properties of binomial coefficients}

We begin with a brief account for three arithmetic properties of the
binomial coefficients $\binom{n}{k}=\frac{n!}{k!(n-k)!}$ which are required
by computing with cohomology of the adjoint Lie group $PSU(n)$, where $1\leq
k\leq n$.

Consider $b_{n,k}=g.c.d.\{\binom{n}{1},\cdots ,\binom{n}{k}\}$, $1\leq k\leq
n$. Since $b_{n,k+1}\mid b_{n,k}$ with $b_{n,1}=n,b_{n,n}=1$ one defines the
integers $a_{n,k}:=\frac{b_{n,k-1}}{b_{n,k}}$, $k\geq 2$. The following
result has been shown in \cite{DLin}, and is required by proving i) of Lemma
4.6.

\begin{enumerate}
\item[(5.5)] If $n>2$ is an integer with the prime factorization $%
n=p_{1}^{r_{1}}\cdots p_{t}^{r_{t}}$, then $a_{n,k}=p_{i}$ or $1$ in
accordance to $k\in Q_{p_{i}}(n)$ or $k\in Q_{0}(n)$, where $1\leq i\leq t$.
\end{enumerate}

For an integer $m$ denote by $ord_{p}m$ the biggest integer $a$ so that $m$
is divisible by the power $p^{a}$. The following result is needed by proving
i) of Lemma 5.3. Assume that $n=p^{r}n^{\prime }$ with $(n\prime ,p)=1$.

\begin{enumerate}
\item[(5.6)] For any $p^{t}\leq s<p^{t+1}$ with $t+1<r$ one has $ord_{p}%
\binom{n}{s}\geq r-t$, while the equality holds iff $s=p^{t}$.
\end{enumerate}

The remaining part of this section is devoted to a proof of Proposition 4.8.
The calculation will be based on the following result

\begin{enumerate}
\item[(5.7)] For any $1\leq s\leq r$\ there is a sequence $\{h_{1},\cdots
,h_{s}\}$\ of integers such that

$\quad \binom{p^{r}}{p^{s}}-p^{r-s}=h_{1}\binom{p^{r}}{p^{s-1}}+h_{2}\binom{%
p^{r}}{p^{s-2}}+\cdots +h_{s}\binom{p^{r}}{1}$.
\end{enumerate}

The center of the special unitary group $SU(n)$ is the cyclic group $\mathbb{%
Z}_{n}$ generated by the diagonal matrix $diag\{e^{\frac{2\pi }{n}i},\cdots
,e^{\frac{2\pi }{n}i}\}\in SU(n)$. The total space of the circle bundle $C$
on $PSU(n)$ (see (1.4)) is the unitary group $U(n)$, whose maximal torus $T$
is $\{diag\{e^{i\theta _{1}},\cdots ,e^{i\theta _{n}}\}\mid $ $\theta
_{i}\in \lbrack 0,2\pi ]\}$. In the terminologies of Section 3.4 we can
furnish the group $H^{1}(T)$ with the basis $\left\{ t_{1},\cdots
,t_{n-1},t_{0}\right\} $ in which the subset $\left\{ t_{1},\cdots
,t_{n-1}\right\} $ corresponds to a set of fundamental dominant weights of
the group $PSU(n)$. Consequently, in the short exact sequence (3.9)
associated to the bundle $C$

\begin{quote}
$0\rightarrow E_{2}^{\ast ,k}(PSU(n))\overset{C^{\ast }}{\rightarrow }%
E_{2}^{\ast ,k}(U(n))\overset{\overline{\theta }}{\rightarrow }E_{2}^{\ast
,k-1}(PSU(n))\rightarrow 0$
\end{quote}

\noindent one has

\begin{quote}
$E_{{\small 2}}^{{\small \ast ,k}}(PSU(n))=H^{\ast }(\frac{{\small U(n)}}{T}%
)\otimes \Lambda (t_{{\small 1}},\cdots ,t_{{\small n-1}})$;

(resp. $E_{2}^{\ast ,k}(U(n))=H^{\ast }(\frac{{\small U(n)}}{T})\otimes
\Lambda (t_{1},\cdots ,t_{n-1},t_{0})$),
\end{quote}

\noindent while by ii) of Remark2.5 the $d_{2}$ actions on $E_{{\small 2}}^{%
{\small \ast ,k}}$ are determined by

\begin{quote}
$\tau (t_{{\small 1}})=2\omega _{1}-\omega _{2};\tau (t_{{\small 2}%
})=-\omega _{1}+2\omega _{2}-\omega _{3};\cdots ;$

$\tau (t_{{\small n-1}})=-\omega _{n-2}+2\omega _{n-1}$ and $\tau ^{\prime
}(t_{{\small 0}})=\omega _{1}$
\end{quote}

\noindent where $\tau $ and $\tau ^{\prime }$ are the transgression in $\pi $
and $\pi ^{\prime }$ (see in the diagram (3.8)), respectively. It follows
that if $\left\{ c_{2},\cdots ,c_{n}\right\} $ is the set of primary
characteristic polynomials for the group $SU(n)$ over $\mathbb{Z}$ given by
Table 5.4, then

\begin{quote}
$H^{\ast }(U(n))=E_{3}^{\ast ,\ast }(U(n))=\Lambda (\gamma _{1},\gamma
_{3},\cdots ,\gamma _{2n-1})$,
\end{quote}

\noindent where, in terms of Lemmas 3.2 and 3.3,

\begin{quote}
$\gamma _{1}=-(n-1)\otimes t_{1}-(n-2)\otimes t_{2}-\cdots -1\otimes
t_{n-1}+n\otimes t_{0}$,

$\gamma _{2r-1}=\left[ \widetilde{c}_{r}\right] $.
\end{quote}

\noindent In the following result we assume that $n=p^{r}$ with $p$ a prime.

\bigskip

\noindent \textbf{Lemma 5.5.} \textsl{For each} $p^{s}\in Q_{p}(n)$ \textsl{%
the class }$\gamma _{2p^{s}-1}$\textsl{\ has a characteristic polynomial of
the form }$\alpha _{p^{s}}+p^{r-s}\omega _{1}^{p^{s}}$\textsl{\ with }$%
\alpha _{p^{s}}\in \left\langle \func{Im}\widetilde{\tau }\right\rangle $%
\textsl{. As results}

\begin{quote}
\textsl{i) }$\gamma _{2p^{s}-1}=\left[ \widetilde{\alpha }%
_{p^{s}}+p^{r-s}\omega _{1}^{p^{s}-1}\otimes t_{0}\right] $\textsl{\ with} $%
\widetilde{\alpha }_{p^{s}}\in E_{2}^{\ast ,1}(PSU(n));$

\textsl{ii)} $d_{2}(\widetilde{\alpha }_{p^{s}})=-p^{r-s}\omega _{1}^{p^{s}}$%
.
\end{quote}

\noindent \textbf{Proof.} For a given $p^{s}\in Q_{p}(n)$ let $%
\{h_{1},\cdots ,h_{s}\}$ be the sequence of integers satisfying (5.7).
Consider the polynomial

\begin{quote}
$\alpha _{p^{s}}=c_{p^{s}}-h_{1}\omega ^{p^{s}-p^{s-1}}c_{p^{s-1}}-\cdots
-h_{s}\omega ^{p^{s}-1}c_{1}-p^{r-s}\omega _{1}^{p^{s}}\in \mathbb{Z}[\omega
_{1},\cdots ,\omega _{n-1}]$
\end{quote}

\noindent From $\alpha _{p^{s}}\mid _{\tau (t_{1})=\cdots =\tau (t_{n-1})=0%
\QTR{sl}{\ }}=0$ by (5.7) and a) of Lemma 5.2 one finds that $\alpha
_{p^{s}}\in \left\langle \func{Im}\widetilde{\tau }\right\rangle $.
Moreover, with respect to the surjection

\begin{quote}
$f:$ $\mathbb{Z}\left[ \omega _{1},\cdots ,\omega _{n-1}\right] \rightarrow
H^{\ast }(U(n)/T)${\small \ }(see Section 3.2)
\end{quote}

\noindent the polynomial $\alpha _{p^{s}}+p^{r-s}\omega _{1}^{p^{2}}\in 
\mathbb{Z}\left[ \omega _{1},\cdots ,\omega _{n-1}\right] $ has a lift to $%
E_{2}^{\ast ,1}(U(n))$ (see (3.1)) with the form

\begin{center}
$\widetilde{c}_{p^{s}}-h_{1}\omega ^{p^{s}-p^{s-1}}\widetilde{c}%
_{p^{s-1}}-h_{2}\omega ^{p^{s}-p^{s-2}}\widetilde{c}_{p^{s-2}}-\cdots
-h_{s}\omega ^{p^{s}-1}\widetilde{c}_{1}$

$=\widetilde{c}_{p^{s}}-d_{2}(h_{1}\omega ^{p^{s}-p^{s-1}-1}\widetilde{c}%
_{p^{s-1}}+h_{2}\omega ^{p^{s}-p^{s-2}-1}\widetilde{c}_{p^{s-2}}+\cdots
+h_{s}\omega ^{p^{s}-2}\widetilde{c}_{1})\otimes t_{0})$.
\end{center}

\noindent These shows that $\gamma _{2p^{s}-1}=\left[ \widetilde{\alpha }%
_{p^{s}}+p^{r-s}\omega _{1}^{p^{s}-1}\otimes t_{0}\right] $\textsl{\ }with $%
\widetilde{\alpha }_{p^{s}}\in E_{2}^{\ast ,1}(PSU(n))$.$\square $

\bigskip

By ii) of Lemma 5.5 the $1$--forms $p\cdot \widetilde{\alpha }%
_{p^{s}}-\omega ^{p^{s}-p^{s-1}}\widetilde{\alpha }_{p^{s-1}}\in E_{2}^{\ast
,1}(PSU(n))$ with $1\leq s\leq r$ are $d_{2}$ closed and therefore, define
the cohomology classes

\begin{enumerate}
\item[(5.8)] $\rho _{2p^{s}-1}:=\left[ p\widetilde{\alpha }_{p^{s}}-\omega
^{p^{s}-p^{s-1}}\widetilde{\alpha }_{p^{s-1}}\right] \in E_{3}^{\ast
,1}(PSU(n))$, $1\leq s\leq r$.
\end{enumerate}

\noindent Moreover, the computation

\begin{quote}
$C^{\ast }(\rho _{2p^{s}-1})-p\cdot \gamma _{2p^{s}-1}$

$=\left[ p\widetilde{\alpha }_{p^{s}}-\omega ^{p^{s}-p^{s-1}}\widetilde{%
\alpha }_{p^{s-1}}]-[p\widetilde{\alpha }_{p^{s}}+p^{r-s+1}\omega
_{1}^{p^{s}-1}\otimes t_{0}\right] $

$=\left[ -\omega ^{p^{s}-p^{s-1}}\widetilde{\alpha }_{p^{s-1}}-p^{r-s+1}%
\omega _{1}^{p^{s}-1}\otimes t_{0}\right] $

$=\left[ d_{2}(-\omega ^{p^{s}-p^{s-1}-1}\widetilde{\alpha }%
_{p^{s-1}}\otimes t_{0})\right] =0$
\end{quote}

\noindent indicates that $C^{\ast }(\rho _{2p^{s}-1})=p\cdot \gamma
_{2p^{s}-1}$, where the second equality comes from i) of Lemma 5.5, while
the fourth one follows from ii) of Lemma 5.5 and $d_{2}(t_{0})=\omega _{1}$.
Granted with the classes $\rho _{2p^{s}-1}$\ ($1\leq s\leq r$) defined by
(5.8) we are ready to show Proposition 4.8.

\bigskip

\noindent \textbf{Proof of Proposition 4.8.} In term of the formulae (4.9)
and (4.3) of the rings $H^{\ast }(U(n))$ and $H^{\ast }(U(n);\mathbb{F}_{p})$
the reduction $r_{p}$ on $H^{\ast }(U(n))$ satisfies the relation $%
r_{p}(\gamma _{I})=\xi _{I}$. Moreover, if $I=(p^{i_{1}},\cdots
,p^{i_{k}})\in Q_{p}(n)$\textsl{\ }with\textsl{\ }$1\leq i_{1}<\cdots
<i_{k}<r$ one has $\theta (\xi _{I})=0$ by (4.5). The commutative diagram
induced by $r_{p}$

\begin{quote}
\begin{tabular}{lll}
$H^{\ast }(U(n))$ & $\overset{\theta }{\rightarrow }$ & $H^{\ast }(PSU(n))$
\\ 
$r_{p}\downarrow $ &  & $\allowbreak r_{p}\downarrow $ \\ 
$H^{\ast }(U(n);\mathbb{F}_{p})$ & $\overset{\theta }{\rightarrow }$ & $%
H^{\ast }(PSU(n);\mathbb{F}_{p})$%
\end{tabular}%
.
\end{quote}

\noindent now concludes that $\allowbreak r_{p}(\theta (\gamma _{I}))=0$.
This shows assertion i) of Proposition 4.8.

With $\gamma _{2p^{s}-1}=\left[ \widetilde{\alpha }_{p^{s}}+p^{r-s}\omega
_{1}^{p^{s}-1}\otimes t_{0}\right] $ by i) of Lemma 5.5 one gets by (3.10)
that

\begin{quote}
$\theta (\gamma _{2p^{s}-1})=\left[ \overline{\theta }(\widetilde{\alpha }%
_{p^{s}}+p^{r-s}\omega _{1}^{p^{s}-1}\otimes t_{0})\right] =p^{r-s}\omega
^{p^{s}-1}$,
\end{quote}

\noindent where we have made use of the relations $\omega _{1}=\varpi $ by
Lemma 4.2, and $\pi ^{\ast }(\varpi )=\omega $ by iii) of Theorem 3.9. This
shows formula ii) of Proposition 4.8.

Finally, for the formula iii) of Proposition 4.8 one gets from (3.10), the
decomposition $\gamma _{I}=\gamma _{I^{e}}[\widetilde{\alpha }%
_{p^{i_{k}}}+p^{r-i_{k}}\omega _{1}^{p^{i_{k}}-1}\otimes t_{0}]$ by i) of
Lemma 5.5, as well as the relation $t_{0}^{2}=0$ on $E_{2}^{\ast ,\ast }$,
the following relation in $E_{3}^{\ast ,\ast }(PSU(n))$

\begin{enumerate}
\item[(5.9)] $\overline{\theta }(\gamma _{I})=\overline{\theta }(\gamma
_{I^{e}})\widetilde{a}_{p^{i_{k}}}+p^{r-i_{k}}\omega ^{p^{i_{k}}-1}%
\widetilde{a}_{I^{e}}$,
\end{enumerate}

\noindent where $\widetilde{a}_{K}=\underset{s\in K}{\Pi }\widetilde{a}_{s}$%
. Granted with the fact that the class $\overline{\theta }(\gamma _{I^{e}})$
is divisible by $p$ the formula iii) is verified by the following calculation

\begin{quote}
$\overline{\theta }(\gamma _{I})=(\frac{1}{p}\theta (\gamma _{I^{e}}))\cdot p%
\widetilde{a}_{p^{i_{k}}}+p^{r-i_{k}}\omega ^{p^{i_{k}}-1}\widetilde{a}%
_{I^{e}}$

$=\frac{1}{p}\theta (\gamma _{I^{e}})(\rho _{2p^{i_{k}}-1}+\omega
^{p^{i_{k}}-p^{i_{k}-1}}\widetilde{a}_{p^{i_{k}-1}-1})+p^{r-i_{k}}\omega
^{p^{i_{k}}-1}\widetilde{a}_{I^{e}}$ (by (5.8))

$=\frac{1}{p}\theta (\gamma _{I^{e}})\rho _{2p^{i_{k}}-1}+\frac{1}{p}\omega
^{p^{i_{k}}-p^{i_{k}-1}}(\theta (\gamma _{I^{e}})\widetilde{a}%
_{p^{i_{k}-1}-1}+p^{r-i_{k}+1}\omega ^{p^{i_{k}-1}-1}\widetilde{a}_{I^{e}})$

$=\frac{1}{p}\theta (\gamma _{I^{e}})\rho _{2p^{i_{k}}-1}+\frac{1}{p}\omega
^{p^{i_{k}}-p^{i_{k}-1}}\theta (\gamma _{I^{\partial }})$ (by (5.9)).$%
\square $
\end{quote}


\begin{thebibliography}{99}
\bibitem{BB} P. F. Baum, W. Browder, The cohomology of quotients of
classical groups. Topology 3(1965), 305--336.

\bibitem{Ch} C. Chevalley, Sur les D\'{e}compositions Cellulaires des
Espaces G/B, in Algebraic groups and their generalizations: Classical
methods, W. Haboush ed. Proc. Symp. in Pure Math. 56 (part 1) (1994), 1-26.

\bibitem{D} J. Dieudonn\'{e}, A History of Algebraic and Differential
Topology 1900--1960, Boston; Basel; Birkh\"{a}user, 1989.

\bibitem{DLin} H. Duan and X. Lin, Topology of unitary groups and the prime
orders of binomial coefficients, arXiv:1502.00401.

\bibitem{DL} H. Duan, SL. Liu, The isomorphism type of the centralizer of an
element in a Lie group, Journal of algebra, 376(2013), 25-45.

\bibitem{DZ1} H. Duan, Xuezhi Zhao, Schubert presentation of the integral
cohomology ring of the flag manifolds $G/T$, LMS J. Comput. Math. Vol.18,
no.1(2015), 489-506.

\bibitem{DZ2} H. Duan, Xuezhi Zhao, Schubert calculus and cohomology of Lie
groups. Part I. 1--connected Lie groups, math.AT (math.AG). arXiv:0711.2541.

\bibitem{DZ3} H. Duan, Xuezhi Zhao, Schubert calculus and the Hopf algebra
structures of exceptional Lie groups, Forum. Math. Vol,25, no.1(2014),
113--140.

\bibitem{G} A. Grothendieck, Torsion homologique et sections rationnelles,
Sem. C. Chevalley, ENS 1958, expos\'{e} 5, Secreatariat Math. IHP, Paris,
1958.

\bibitem{H} J. E. Humphreys, Introduction to Lie algebras and representation
theory, Graduated Texts in Math. 9, Springer-Verlag New York, 1972.

\bibitem{IKT} K. Ishitoya, A. Kono, H. Toda, Hopf algebra structure of mod 2
cohomology of simple Lie groups. Publ. Res. Inst. Math. Sci. 12 (1976/77),
no. 1, 141--167.

\bibitem{K} V.G. Ka\v{c}, Torsion in cohomology of compact Lie groups and
Chow rings of reductive algebraic groups, Invent. Math. 80(1985), no. 1,
69--79.

\bibitem{Ko} A. Kono, Hopf algebra structure of simple Lie groups, J. Math.
Kyoto Univ. 17 (1977), no. 2, 259--298.

\bibitem{Ma} W. S. Massey, On the cohomology ring of a sphere bundle, J.
Math. Mech. 7 1958 265-289.

\bibitem{Mc} J. McCleary, A user's guide to spectral sequences, Second
edition. Cambridge Studies in Advanced Mathematics, 58. Cambridge University
Press, Cambridge, 2001.

\bibitem{M} R. Marlin, Une conjecture sur les anneaux de Chow $A(G,\mathbb{Z}%
)$ renforc\'{e}e par un calcul formel, Effective methods in algebraic
geometry (Castiglioncello, 1990), 299--311, Progr. Math., 94, Birkh\"{a}user
Boston, Boston, MA, 1991.

\bibitem{MS} J. W. Milnor and J. D. Stasheff, Characteristic classes. Annals
of Mathematics Studies, No. 76. Princeton University Press, Princeton, N.
J.; University of Tokyo Press, Tokyo, 1974.

\bibitem{MT} M. Mimura and H. Toda, Topology of Lie groups. I, II.,
Translations of Mathematical Monographs, 91, American Mathematical Society,
Providence, RI, 1991.

\bibitem{Re} M. Reeder, On the cohomology of compact Lie groups. Enseign.
Math. (2) 41 (1995), no. 3-4, 181--200.

\bibitem{R} C. A. Ruiz, The cohomology of the complex projective Stiefel
manifold, Trans. Amer. Math. Soc. 146(1969), 541--547.

\bibitem{S} M. R. Sepanski, Compact Lie groups, Graduate Texts in
Mathematics, 235. Springer, New York, 2007.

\bibitem{SE} N. E. Steenrod and D. B. A. Epstein, Cohomology Operations,
Ann. of Math. Stud., Princeton Univ. Press, Princeton, NJ, 1962.

\bibitem{W} A. Weil, Foundations of algebraic geometry, American
Mathematical Society, Providence, R.I. 1962.
\end{thebibliography}
\end{document}